\documentclass[12pt,a4paper,leqno]{article}
\usepackage{etoolbox} \usepackage{comment}
\usepackage{amsfonts}
\usepackage{multirow}
\usepackage{float}
\usepackage[makeroom]{cancel}
\usepackage[a4paper,includeheadfoot,margin=2.54cm]{geometry}
\usepackage{times}
\usepackage[amsthm,thmmarks]{ntheorem}
\usepackage{amsmath,amsfonts,amssymb,mathtools}
\usepackage{amsmath,xypic}
\usepackage{pdflscape}
\usepackage[export]{adjustbox} 
\usepackage[cal=boondoxo,calscaled=.96,scr=rsfs]{mathalpha} \usepackage{multicol}
\usepackage{hyphenat}
\usepackage[usenames,dvipsnames]{color}
\usepackage{rotating}
\usepackage{tikz-cd}
\usepackage{longtable}
\usepackage{caption}
\usepackage[most]{tcolorbox}

\usepackage{epigraph}
\usepackage{enumitem}
\setlist[enumerate]{listparindent=0.5in}
\usepackage{mathrsfs}
\newcommand{\be}{\begin{equation}}
\newcommand{\ee}{\end{equation}}
\newcommand{\bes}{\begin{equation*}}
\newcommand{\ees}{\end{equation*}}
\newcommand{\bea}{\begin{eqnarray}}
\newcommand{\eea}{\end{eqnarray}}
\newcommand{\beas}{\begin{eqnarray}}
\newcommand{\eeas}{\end{eqnarray}}
\newcommand{\ben}{\begin{note}}
\newcommand{\een}{\end{note}}
\newcommand{\bexl}{\vskip0.1em\noindent\hrulefill\vskip1em\begin{ExerciseList}}
\newcommand{\eexl}{\end{ExerciseList}\hrulefill}

\newcommand{\bthm}{\begin{theorem}}
\newcommand{\ethm}{\end{theorem}}
\newcommand{\bpro}{\begin{prop}}
\newcommand{\epro}{\end{prop}}
\newcommand{\bcor}{\begin{corollary}}
\newcommand{\ecor}{\end{corollary}}
\newcommand{\bcon}{\begin{conjecture}}
\newcommand{\econ}{\end{conjecture}}
\newcommand{\bp}{\begin{proof}}
\newcommand{\ep}{\end{proof}}
\newcommand{\blem}{\begin{lemma}}
\newcommand{\elem}{\end{lemma}}
\newcommand{\bn}{\begin{note}}
\newcommand{\en}{\end{note}}
\newcommand{\benum}{\begin{enumerate}}
\newcommand{\eenum}{\end{enumerate}}
\newcommand{\bed}{\begin{defn}}
\newcommand{\eed}{\end{defn}}
\newcommand{\brem}{\begin{remark}}
\newcommand{\erem}{\end{remark}}

\newcommand{\btik}{\begin{tikzpicture}\begin{axis}[scale=0.5,axis y line=center, axis x line=middle]}
\newcommand{\etik}{\end{axis}\end{tikzpicture}}

\let\into=\hookrightarrow
\let\mapsto=\longmapsto

\newcommand{\upperRomannumeral}[1]{\uppercase\expandafter{\romannumeral#1}}

 \usepackage{graphicx}
\usepackage[pagewise]{lineno}
\usepackage{stackengine}
\usepackage{xr}
\newcommand{\construntilts}[2]{\cite[{#1\ref{U-#2}}]{joshi-untilts}}
\newcommand{\constrone}[2]{\cite[{#1\ref{I-#2}}]{joshi-teich}}
\newcommand{\constrtwo}[2]{\cite[{#1\ref{II-#2}}]{joshi-teich-estimates}}
\newcommand{\constrtwoh}[2]{\cite[{#1\ref{II5-#2}}]{joshi-teich-def}}
\newcommand{\constrthr}[2]{\cite[{#1\ref{III-#2}}]{joshi-teich-rosetta}}

\usepackage{stmaryrd}
\usepackage[normalem]{ulem} \usepackage[cal=boondoxo,calscaled=.96,scr=rsfs]{mathalpha}
\DeclareMathAlphabet{\mathscrbf}{OMS}{mdugm}{b}{n}

\usepackage{titlesec}
\newtoggle{arxiv}
\toggletrue{arxiv}
\usepackage{natbib}
	\let\cite=\citep
\usepackage{fancyhdr}
	\usepackage[colorlinks,citecolor=blue]{hyperref}
\usepackage[capitalise,noabbrev]{cleveref}
\usepackage{fancyhdr}
\usepackage{cleveref}
\newtoggle{draft}
\togglefalse{draft}

\theoremstyle{theorem}
\newtheorem{theorem}[equation]{Theorem}      \newtheorem{lemma}[equation]{Lemma}          \newtheorem{corollary}[equation]{Corollary}  \newtheorem{proposition}[equation]{Proposition}

\newtheorem{conj}[equation]{Conjecture}

\newtheorem{defn}[equation]{Definition}

\theoremsymbol{\ensuremath{\bullet}}

\newtheorem{remark}[equation]{Remark}
\theoremsymbol{\ensuremath{\bullet}}

\numberwithin{equation}{section}

\let\oldproofname=\proofname
\renewcommand{\proofname}{{\bfseries\sffamily\textup{\oldproofname}}}

\titleformat{\subsection}[runin]{\normalfont\bfseries}{\thesubsection}{.5em}{}[{\ \ }]
\titlespacing{\subsection}{0pt}{1.5ex plus .1ex minus .2ex}{0pt}
\titleformat{\subsubsection}[runin]{\normalfont\bfseries}{\thesubsubsection}{.5em}{}[{\ \ }]
\titlespacing{\subsubsection}{0pt}{1.5ex plus .1ex minus .2ex}{0pt}

\newcommand{\nwss}{\numberwithin{equation}{subsection}}

\crefname{equation}{equation}{Equation}
\crefname{section}{§}{§§}
\crefname{subsection}{§}{§§}
\crefname{subsubsection}{§}{§§}

\usepackage{microtype}
\usepackage[margin=15pt,font=small,labelfont={bf,sf}]{caption}
\usepackage{sectsty}
\allsectionsfont{\sffamily}

\newtcolorbox[auto counter,number within=subsection,
crefname={Table}{Table}]{boxedcontent}[2][]{colback=white,coltitle=black,colframe=white!25!brown,fonttitle=\bfseries,
	title=Table \thetcbcounter: #2,#1}

\let\into=\hookrightarrow
\let\isom=\simeq

\newcommand{\A}{\mathscr{A}}
\newcommand{\abs}[1]{\left\vert#1\right\vert}

\newcommand{\bF}{{\bar{F}}}
\newcommand{\bQ}{{\bar{\Q}}}

\newcommand{\C}{{\mathbb C}}

\newcommand{\F}{{\mathbb F}}

\newcommand{\gal}{{\rm Gal}}

\newcommand{\mydot}{{\scriptscriptstyle{\bullet}}}
\newcommand{\N}{\mathscr{N}}

\newcommand{\Q}{{\mathbb Q}}
\newcommand{\R}{{\mathbb R}}

\newcommand{\Spec}{{\rm Spec}}

\newcommand{\Z}{{\mathbb Z}}

\renewcommand{\O}{{\mathscr O}}
\renewcommand{\P}{{\mathbb P}}
\renewcommand{\wp}{{\mathfrak p}}

\newcommand{\fa}{{\mathscr{A}}}
\newcommand{\invlim}{\varprojlim}

\renewcommand{\bpro}{\begin{proposition}}
	\renewcommand{\epro}{\end{proposition}}
\renewcommand{\bcon}{\begin{conj}}
	\renewcommand{\econ}{\end{conj}}

\title{Construction of Arithmetic Teichmuller Spaces IV:\\ Proof of the $abc$-conjecture 
\\	\textcolor{red}{{\large\bf Preliminary version for comments}}}
\author{Kirti Joshi}

\newcommand{\Address}{\bigskip\noindent{\footnotesize\textsc{{Math. department, University of Arizona, 617 N Santa Rita, Tucson
		85721-0089, USA.}}\par\nopagebreak 
\noindent\textit{Email:}	\texttt{kirti@math.arizona.edu}}}

\begin{document}
	\maketitle
\epigraphwidth0.65\textwidth

\lhead{}

\iftoggle{draft}{\pagewiselinenumbers}{\relax}
\newcommand{\needsproof}{\textcolor{red}{{\large THIS NEEDS A PROOF\ }}}
\newcommand{\act}{\curvearrowright}
\newcommand{\lmp}{{\Pi\act\Ot}}
\newcommand{\lmpi}{{\lmp}_{\int}}
\newcommand{\lmpf}{\lmp_F}
\newcommand{\Om}{\O^{\times\mu}}
\newcommand{\Omf}{\O^{\times\mu}_{\bF}}
\renewcommand{\N}{\mathbb{N}}
\newcommand{\yoga}{Yoga}
\newcommand{\gl}[1]{{\rm GL}(#1)}
\newcommand{\bK}{\overline{K}}
\newcommand{\reptrip}{\rho:G_K\to\gl V}
\newcommand{\reptripp}[1]{\rho\circ\alpha:G_{\ifstrempty{#1}{K}{{#1}}}\to\gl V}
\newcommand{\benumlab}{\begin{enumerate}[label={{\bf(\arabic{*})}}]}
\newcommand{\benumlabresume}{\begin{enumerate}[resume,label={{\bf(\arabic{*})}}]}
\newcommand{\ord}{\mathop{\rm ord}\nolimits}	
\newcommand{\kcs}{K^\circledast}
\newcommand{\lcs}{L^\circledast}
\renewcommand{\A}{\mathbb{A}}
\newcommand{\bfq}{\bar{\mathbb{F}}_q}
\newcommand{\tripod}{\P^1-\{0,1728,\infty\}}

\newcommand{\vseq}[2]{{#1}_1,\ldots,{#1}_{#2}}
\newcommand{\anab}[4]{\left({#1},\{#3 \}\right)\anabelmap\left({#2},\{#4 \}\right)}

\newcommand{\gln}{{\rm GL}_n}
\newcommand{\glo}[1]{{\rm GL}_1(#1)}
\newcommand{\glt}[1]{{\rm GL_2}(#1)}

\newcommand{\iut}{\cite{mochizuki-iut1, mochizuki-iut2, mochizuki-iut3,mochizuki-iut4}}
\newcommand{\topics}{\cite{mochizuki-topics1,mochizuki-topics2,mochizuki-topics3}}

\newcommand{\linv}{\mathfrak{L}}
\newcommand{\bedef}{\begin{defn}}
\newcommand{\eedef}{\end{defn}}
\renewcommand{\act}[1][]{\overset{#1}{\curvearrowright}}
\newcommand{\bfx}{\overline{F(X)}}
\newcommand{\anabelmap}{\leftrightsquigarrow}
\newcommand{\ban}[1][G]{\mathscr{B}({#1})}
\newcommand{\pit}{\Pi^{temp}}
 
 \newcommand{\bL}{\overline{L}}
 \newcommand{\bkm}{\bK_M}
 \newcommand{\vbk}{v_{\bK}}
 \newcommand{\vbkm}{v_{\bkm}}
\newcommand{\ocs}{\O^\circledast}
\newcommand{\ot}{\O^\triangleright}
\newcommand{\ocsk}{\ocs_K}
\newcommand{\otk}{\ot_K}
\newcommand{\ok}{\O_K}
\newcommand{\oko}{\O_K^1}
\newcommand{\oks}{\ok^*}
\newcommand{\Qpb}{\overline{\Q}_p}
\newcommand{\Qpbh}{\widehat{\overline{\Q}}_p}
\newcommand{\tr}{\triangleright}
\newcommand{\ocpt}{\O_{\C_p}^\tr}
\newcommand{\ocpf}{\O_{\C_p}^\flat}
\newcommand{\sG}{\mathscr{G}}
\newcommand{\sX}{\mathscr{X}}
\newcommand{\sxfe}{\sX_{F,E}}
\newcommand{\sxfep}{\sX_{F,E'}}
\newcommand{\sxcpte}{\sX_{\cpt,E}}
\newcommand{\sxcptep}{\sX_{\cpt,E'}}
\newcommand{\loglt}{\log_{\sG}}
\newcommand{\fc}{\mathfrak{t}}
\newcommand{\ku}{K_u}
\newcommand{\kup}{\ku'}
\newcommand{\kt}{\tilde{K}}
\newcommand{\sGpf}{\sG(\O_K)^{pf}}
\newcommand{\hgm}{\widehat{\mathbb{G}}_m}
\newcommand{\bE}{\overline{E}}
\newcommand{\sY}{\mathscr{Y}}
\newcommand{\syfe}{\mathscr{Y}_{F,E}}
\newcommand{\syfep}{\mathscr{Y}_{F,E'}}
\newcommand{\syfqp}[1]{\mathscr{Y}_{\cptl{#1},\Q_p}}
\newcommand{\syfqpe}[1]{\mathscr{Y}_{{#1},E}}
\newcommand{\syfqpep}[1]{\mathscr{Y}_{{#1},E'}}
\newcommand{\fJ}{\mathfrak{J}}
\newcommand{\fM}{\mathfrak{M}}
\newcommand{\locvar}{local arithmetic-geometric anabelian variation of fundamental group of $X/E$ at $\wp$}
\newcommand{\fjxep}{\fJ(X,E,\wp)}
\newcommand{\fjxe}{\fJ(X,E)}
\newcommand{\fpc}[1]{\widehat{{\overline{\F_p(({#1}))}}}}
\newcommand{\cpt}{\C_p^\flat}
\newcommand{\cpti}{\C_{p_i}^\flat}
\newcommand{\cptl}[1]{\C_{p,{#1}}^\flat}
\newcommand{\fja}[1]{\fJ^{\rm arith}({#1})}
\newcommand{\ainfe}{A_{\inf,E}(\O_F)}
\renewcommand{\ainfe}{W_{\O_E}(\O_F)}
\newcommand{\gmh}{\widehat{\mathbb{G}}_m}
\newcommand{\sE}{\mathscr{E}}
\newcommand{\bpi}{B^{\varphi=\pi}}
\newcommand{\bpip}{B^{\varphi=p}}
\newcommand{\onto}{\twoheadrightarrow}

\newcommand{\cpmax}{{\C_p^{\rm max}}}
\newcommand{\xan}{X^{an}}
\newcommand{\yan}{Y^{an}}
\newcommand{\bPi}{\overline{\Pi}}
\newcommand{\bPit}{\bPi^{\rm{\scriptscriptstyle temp}}}
\newcommand{\Pit}{\Pi^{\rm{\scriptscriptstyle temp}}}
\renewcommand{\pit}[1]{\Pi^{\scriptscriptstyle temp}_{#1}}
\newcommand{\pitk}[2]{\Pi^{\scriptscriptstyle temp}_{#1;#2}}
\newcommand{\pio}[1]{\pi_1({#1})}
\newcommand{\fTeich}{\widetilde{\fJ(X/L)}}
\newcommand{\ssep}{\S\,} \newcommand{\vphi}{\varphi}
\newcommand{\sgt}{\widetilde{\sG}}
\newcommand{\sxqp}{\mathscr{X}_{\cpt,\Q_p}}
\newcommand{\syQp}{\mathscr{Y}_{\cpt,\Q_p}}

\setcounter{tocdepth}{2}

\newcommand{\mywork}[1]{\textcolor{teal}{#1}}

\togglefalse{draft}
\newcommand{\FF}{\cite{fargues-fontaine}}
\iftoggle{draft}{\pagewiselinenumbers}{\relax}

\newcommand{\attportion}{Sections~\cref{se:number-field-case}, \cref{se:construct-att}, \cref{se:relation-to-iut}, \cref{se:self-similarity} and \cref{se:applications-elliptic}}

\newcommand{\inithtdata}{\cite[{\ssep\ref{III-ss:theta-data-fixing}, \ssep\ref{III-ss:theta-data-fixing2}}]{joshi-teich-rosetta}} 
\newcommand{\initassumptions}{Initial Theta Data \cite[{\ssep\ref{III-ss:elliptic-curve-assumptions}, \ssep\ref{III-ss:theta-data-fixing}, \ssep\ref{III-ss:theta-data-fixing2}}]{joshi-teich-rosetta}}
\newcommand{\assumptions}{\cref{sss:degrees}, \cref{sss:field-assump}}

\newcommand{\Pib}{\overline{\Pi}}
\newcommand{\four}{Sections~\cref{se:grothendieck-conj}, \cref{se:untilts-of-Pi}, and \cref{se:riemann-surfaces}}

\newcommand{\ENDDOCUMENT}{\phantomsection\addcontentsline{toc}{section}{References}\bibliography{../../master/masterofallbibs.bib}\Address\end{document}}

\begin{thebibliography}{73}
\providecommand{\natexlab}[1]{#1}
\providecommand{\url}[1]{\texttt{#1}}
\expandafter\ifx\csname urlstyle\endcsname\relax
  \providecommand{\doi}[1]{doi: #1}\else
  \providecommand{\doi}{doi: \begingroup \urlstyle{rm}\Url}\fi

\bibitem[Amor\'{o}s et~al.(2000)Amor\'{o}s, Bogomolov, Katzarkov, and
  Pantev]{bogomolov00}
J.~Amor\'{o}s, F.~Bogomolov, L.~Katzarkov, and T.~Pantev.
\newblock Symplectic {L}efschetz fibrations with arbitrary fundamental groups.
\newblock \emph{J. Differential Geom.}, 54\penalty0 (3):\penalty0 489--545,
  2000.
\newblock URL \url{http://projecteuclid.org/euclid.jdg/1214339791}.
\newblock With an appendix by Ivan Smith.

\bibitem[Andr\'{e}(2003)]{andre03}
Yves Andr\'{e}.
\newblock On a geometric description of $\gal({\bQ}_p/{\Q}_p)$ and a $p$-adic
  avatar of $\widehat{GT}$.
\newblock \emph{Duke Math. Journal}, 119\penalty0 (1):\penalty0 1--39, 2003.

\bibitem[Beauville(2002)]{beauville02}
Arnaud Beauville.
\newblock The {S}zpiro inequality for higher genus fibrations.
\newblock In \emph{Algebraic geometry}, pages 61--63. de Gruyter, Berlin, 2002.

\bibitem[Bombieri and Gubler(2006)]{bombieri-gubler}
Enrico Bombieri and Walter Gubler.
\newblock \emph{Heights in diophantine geometry}.
\newblock Cambridge University Press, Cambridge, UK ; New York, 2006.

\bibitem[Castelvecchi(2015)]{castelvecchi-1}
Davide Castelvecchi.
\newblock The biggest mystery in mathematics: {S}hinichi {M}ochizuki and the
  impenetrable proof.
\newblock \emph{Nature 526}, 526:\penalty0 178--181, 2015.
\newblock \doi{https://doi.org/10.1038/526178a}.

\bibitem[Castelvecchi(2020)]{castelvecchi-2}
Davide Castelvecchi.
\newblock Mathematical proof that rocked number theory will be published.
\newblock \emph{Nature}, 580:\penalty0 177--177, 2020.
\newblock \doi{https://doi.org/10.1038/d41586-020-00998-2}.

\bibitem[Colmez(2002)]{colmez2000}
Pierre Colmez.
\newblock Espaces de banach de dimension finie.
\newblock \emph{J. Inst. Math. Jussieu}, 1\penalty0 (3):\penalty0 331--439,
  2002.

\bibitem[Dupuy and Hilado(2020)]{dupuy2020probabilistic}
Taylor Dupuy and Anton Hilado.
\newblock Probabilistic {S}zpiro, {B}aby {S}zpiro, and {E}xplicit {S}zpiro from
  {M}ochizuki's corollary 3.12.
\newblock \emph{arXiv preprint arXiv:2004.13108}, 2020.
\newblock URL \url{https://arxiv.org/abs/2004.13108}.

\bibitem[Dusart(2016)]{dusart2016}
Pierre Dusart.
\newblock Explicit estimates of some functions over primes.
\newblock \emph{The Ramanujan Journal}, 45:\penalty0 227--251, 2016.
\newblock \doi{https://doi.org/10.1007/s11139-016-9839-4}.

\bibitem[{E}lkies(1991)]{elkies1991}
Noam~{D} {E}lkies.
\newblock Abc implies mordell.
\newblock \emph{International Mathematics Research Notices}, 1991\penalty0
  (7):\penalty0 99--109, 1991.
\newblock \doi{10.1155/S1073792891000144}.

\bibitem[Faltings(1983)]{faltings1983}
Gerd Faltings.
\newblock Endlichkeitss\"{a}tze f\"{u}r abelsche {V}ariet\"{a}ten \"{u}ber
  {Z}ahlk\"{o}rpern.
\newblock \emph{Invent. Math.}, 73\penalty0 (3):\penalty0 349--366, 1983.
\newblock URL \url{https://doi.org/10.1007/BF01388432}.

\bibitem[Fargues and Fontaine(2018)]{fargues-fontaine}
Laurent Fargues and Jean-Marc Fontaine.
\newblock Courbes et fibr\'{e}s vectoriels en th\'{e}orie de {H}odge
  {$p$}-adique.
\newblock \emph{Ast\'{e}risque}, \penalty0 (406):\penalty0 xiii+382, 2018.
\newblock ISSN 0303-1179.
\newblock With a preface by Pierre Colmez.

\bibitem[Imayoshi and Taniguchi(1992)]{imayoshi-book}
Yoichi Imayoshi and Masahiko Taniguchi.
\newblock \emph{An Introduction to Teichmuller Spaces}.
\newblock Springer Japan, Tokyo, 1st ed. 1992. edition, 1992.

\bibitem[Joshi(2019)]{joshi-formal-groups}
Kirti Joshi.
\newblock Mochizuki's anabelian variation of ring structures and formal groups.
\newblock 2019.
\newblock URL \url{https://arxiv.org/abs/1906.06840}.

\bibitem[Joshi(2020{\natexlab{a}})]{joshi-anabelomorphy}
Kirti Joshi.
\newblock On {M}ochizuki's idea of anabelomorphy and applications.
\newblock 2020{\natexlab{a}}.
\newblock URL \url{https://arxiv.org/abs/2003.01890}.

\bibitem[Joshi(2020{\natexlab{b}})]{joshi-gconj}
Kirti Joshi.
\newblock The absolute {G}rothendieck conjecture is false for
  {F}argues-{F}ontaine curves.
\newblock 2020{\natexlab{b}}.
\newblock URL \url{https://arxiv.org/abs/2008.01228}.
\newblock Preprint.

\bibitem[Joshi(2021{\natexlab{a}})]{joshi-teich}
Kirti Joshi.
\newblock Construction of {A}rithmetic {T}eichmuller {S}paces {I}.
\newblock 2021{\natexlab{a}}.
\newblock URL \url{https://arxiv.org/abs/2106.11452}.

\bibitem[Joshi(2021{\natexlab{b}})]{joshi-teich-summary-comments}
Kirti Joshi.
\newblock Comments on {A}rithmetic {T}eichmuller {T}heory.
\newblock 2021{\natexlab{b}}.
\newblock URL \url{https://arxiv.org/abs/2111.06771}.

\bibitem[Joshi(2022)]{joshi-untilts}
Kirti Joshi.
\newblock Untilts of fundamental groups: construction of labeled isomorphs of
  fundamental groups ({A}rithmetic {H}olomorphic {S}tructures).
\newblock 2022.
\newblock URL \url{https://arxiv.org/abs/2210.11635}.

\bibitem[Joshi(2023{\natexlab{a}})]{joshi-teich-def}
Kirti Joshi.
\newblock Construction of {A}rithmetic {T}eichmuller {S}paces
  {II}$(\frac{1}{2})$: {D}eformations of {N}umber {F}ields.
\newblock 2023{\natexlab{a}}.
\newblock URL \url{https://arxiv.org/abs/2305.10398}.

\bibitem[Joshi(2023{\natexlab{b}})]{joshi-teich-estimates}
Kirti Joshi.
\newblock Construction of {A}rithmetic {T}eichmuller {S}paces {II}: Proof a
  local prototype of {M}ochizuki's {C}orollary 3.12.
\newblock 2023{\natexlab{b}}.
\newblock URL \url{https://arxiv.org/abs/2303.01662}.

\bibitem[Joshi(2023{\natexlab{c}})]{joshi-teich-quest}
Kirti Joshi.
\newblock Mochizuki's {C}orollary 3.12 and my quest for its proof.
\newblock 2023{\natexlab{c}}.
\newblock URL \url{https://www.math.arizona.edu/~kirti/joshi-teich-quest.pdf}.

\bibitem[Joshi(2024{\natexlab{a}})]{joshi-teich-abc}
Kirti Joshi.
\newblock Construction of {A}rithmetic {T}eichmuller {S}paces {IV}: Proof of
  the $abc$-conjecture.
\newblock \emph{Preprint}, 2024{\natexlab{a}}.
\newblock URL \url{https://www.arxiv.org}.

\bibitem[Joshi(2024{\natexlab{b}})]{joshi-teich-abc-conj}
Kirti Joshi.
\newblock Construction of {A}rithmetic {T}eichmuller {S}paces {IV}: A proof of
  the $abc$-conjecture.
\newblock \emph{Available at arxiv.org}, 2024{\natexlab{b}}.
\newblock URL \url{https://arxiv.org/pdf/2403.10430.pdf}.

\bibitem[Joshi(2024{\natexlab{c}})]{joshi-teich-rosetta}
Kirti Joshi.
\newblock Construction of {A}rithmetic {T}eichmuller {S}paces {III}: A proof of
  {M}ochizuki’s corollary 3.12 and a {R}osetta {S}tone.
\newblock 2024{\natexlab{c}}.
\newblock URL \url{https://arxiv.org/pdf/2401.13508.pdf}.

\bibitem[Joshi(2025{\natexlab{a}})]{joshi-expository}
Kirti Joshi.
\newblock A worked example illustrating my construction of {Arithmetic
  Teichmuller Spaces}.
\newblock 2025{\natexlab{a}}.

\bibitem[Joshi(2025{\natexlab{b}})]{joshi-final}
Kirti Joshi.
\newblock Final report on the status of the {M}ochizuki-{S}cholze-{S}tix
  {C}ontroversy.
\newblock 2025{\natexlab{b}}.
\newblock URL
  \url{http://www.math.arizona.edu/~kirti/Final-Mochizuki-Scholze-Stix-Controversy.pdf}.

\bibitem[Joshi(June, 2024)]{joshi-report}
Kirti Joshi.
\newblock The status of the {S}cholze-{S}tix {R}eport and an {a}nalysis of the
  {M}ochizuki-{S}cholze-{S}tix {C}ontroversy.
\newblock June, 2024.
\newblock URL
  \url{https://math.arizona.edu/\~kirti/report-on-scholze-stix-mochizuki-controversy.pdf}.

\bibitem[Joshi(November 2022)]{joshi-blogpost}
Kirti Joshi.
\newblock
  \href{https://thehighergeometer.wordpress.com/2022/11/25/a-study-in-basepoints-guest-post-by-kirti-joshi/}{A
  study in basepoints: guest post by {K}irti {J}oshi}, November 2022.
\newblock URL
  \url{https://thehighergeometer.wordpress.com/2022/11/25/a-study-in-basepoints-guest-post-by-kirti-joshi/}.

\bibitem[Joshi(October 2020)]{joshi-untilts-2020}
Kirti Joshi.
\newblock Untilts of fundamental groups: construction of labeled isomorphs of
  fundamental groups ({\bf 2020 version}).
\newblock October 2020.
\newblock URL \url{https://arxiv.org/abs/2010.05748}.

\bibitem[Joshi and Pauly(2015)]{joshi16}
Kirti Joshi and Christian Pauly.
\newblock Hitchin-{M}ochizuki morphism, opers and {F}robenius-destabilized
  vector bundles over curves.
\newblock \emph{Adv. Math.}, 274:\penalty0 39--75, 2015.
\newblock \doi{10.1016/j.aim.2015.01.004}.
\newblock URL \url{http://dx.doi.org/10.1016/j.aim.2015.01.004}.

\bibitem[Kaplansky(1942)]{kaplansky42}
Irving Kaplansky.
\newblock Maximal fields with valuations.
\newblock \emph{Duke Math. J.}, 9:\penalty0 303--321, 1942.
\newblock URL \url{http://projecteuclid.org/euclid.dmj/1077493226}.

\bibitem[Kedlaya and Liu(2015)]{kedlaya-liu15}
Kiran~S. Kedlaya and Ruochuan Liu.
\newblock Relative {$p$}-adic {H}odge theory: foundations.
\newblock \emph{Ast\'{e}risque}, 371:\penalty0 239, 2015.
\newblock ISSN 0303-1179.

\bibitem[Kedlaya and Temkin(2018)]{kedlaya18}
Kiran~S. Kedlaya and Michael Temkin.
\newblock Endomorphisms of power series fields and residue fields of
  {F}argues-{F}ontaine curves.
\newblock \emph{Proc. Amer. Math. Soc.}, 146\penalty0 (2):\penalty0 489--495,
  2018.
\newblock URL \url{https://doi.org/10.1090/proc/13818}.

\bibitem[Kim(1997)]{kim97}
Minhyong Kim.
\newblock Geometric height inequalities and the {K}odaira-{S}pencer map.
\newblock \emph{Compositio Math.}, 105\penalty0 (1):\penalty0 43--54, 1997.

\bibitem[Lang(2002)]{lang-algebra}
Serge Lang.
\newblock \emph{Algebra}, volume 211 of \emph{Graduate Texts in Mathematics}.
\newblock Springer-Verlag, New York, third edition, 2002.
\newblock URL \url{https://doi.org/10.1007/978-1-4613-0041-0}.

\bibitem[{M}asser(1990)]{masser1990}
{D}.~{W}. {M}asser.
\newblock Note on a conjecture of {S}zpiro.
\newblock Number 183, pages 19--23. 1990.
\newblock S\'{e}minaire sur les Pinceaux de Courbes Elliptiques (Paris, 1988).

\bibitem[Matignon and Reversat(1984)]{matignon84}
Michel Matignon and Marc Reversat.
\newblock Sous-corps ferm\'{e}s d'un corps valu\'{e}.
\newblock \emph{J. Algebra}, 90\penalty0 (2):\penalty0 491--515, 1984.
\newblock URL \url{https://doi.org/10.1016/0021-8693(84)90186-8}.

\bibitem[Mirzakhani(2007)]{mirzakhani07}
Maryam Mirzakhani.
\newblock Simple geodesics and {W}eil-{P}etersson volumes of moduli spaces of
  bordered {R}iemann surfaces.
\newblock \emph{Invent. {M}ath.}, 167:\penalty0 179--222, 2007.
\newblock \doi{https://doi.org/10.1007/s00222-006-0013-2}.

\bibitem[Mochizuki(1997)]{mochizuki-local-gro}
Shinichi Mochizuki.
\newblock A version of the {G}rothendieck conjecture for $p$-adic local fields.
\newblock \emph{International Journal of Math.}, 8:\penalty0 499--506, 1997.
\newblock URL
  \url{http://www.kurims.kyoto-u.ac.jp/~motizuki/A%20Version%20of%20the%20Grothendieck%20Conjecture%20for%20p-adic%20Local%20Fields.pdf}.

\bibitem[Mochizuki(2010)]{mochizuki-general-pos}
Shinichi Mochizuki.
\newblock Arithmetic elliptic curves in general position.
\newblock \emph{Math. J. Okayama Univ.}, 52:\penalty0 1--28, 2010.
\newblock ISSN 0030-1566.

\bibitem[Mochizuki(2012)]{mochizuki-topics1}
Shinichi Mochizuki.
\newblock Topics in absolute anabelian geometry {I}: generalities.
\newblock \emph{J. Math. Sci. Univ. Tokyo}, 19\penalty0 (2):\penalty0 139--242,
  2012.

\bibitem[Mochizuki(2013)]{mochizuki-topics2}
Shinichi Mochizuki.
\newblock Topics in absolute anabelian geometry {II}: decomposition groups and
  endomorphisms.
\newblock \emph{J. Math. Sci. Univ. Tokyo}, 20\penalty0 (2):\penalty0 171--269,
  2013.

\bibitem[Mochizuki(2015)]{mochizuki-topics3}
Shinichi Mochizuki.
\newblock Topics in absolute anabelian geometry {III}: global reconstruction
  algorithms.
\newblock \emph{J. Math. Sci. Univ. Tokyo}, 22\penalty0 (4):\penalty0
  939--1156, 2015.

\bibitem[Mochizuki(2016)]{mochizuki-bogomolov}
Shinichi Mochizuki.
\newblock Bogomolov's proof of the geometric version of the {S}zpiro conjecture
  from the point of view of {I}nter-{U}niversal {T}eichm\"{u}ller {T}heory.
\newblock \emph{Res. Math. Sci.}, 3:\penalty0 Paper No. 6, 21, 2016.
\newblock ISSN 2522-0144.
\newblock \doi{10.1186/s40687-016-0057-x}.
\newblock URL \url{https://doi.org/10.1186/s40687-016-0057-x}.

\bibitem[Mochizuki(2018)]{mochizuki-scholze-stix-report}
Shinichi Mochizuki.
\newblock Comments on the manuscript by {S}cholze-{S}tix concerning
  {I}nter-{U}niversal {T}eichmuller {T}heory ({IUTCH}).
\newblock 2018.
\newblock URL \url{https://www.kurims.kyoto-u.ac.jp/~motizuki/Cmt2018-05.pdf}.

\bibitem[Mochizuki(2020)]{mochizuki-gaussian}
Shinichi Mochizuki.
\newblock The mathematics of mutually alien copies: from {G}aussian integrals
  to {I}nter-{U}niversal {T}eichmüller {T}heory.
\newblock 2020.
\newblock URL
  \url{https://www.kurims.kyoto-u.ac.jp/~motizuki/Alien\%20Copies,\%20Gaussians,\%20and\%20Inter-universal\%20Teichmuller\%20Theory.pdf}.

\bibitem[Mochizuki(2021{\natexlab{a}})]{mochizuki-iut1}
Shinichi Mochizuki.
\newblock Inter-{U}niversal {T}eichmuller theory {I}: construction of {H}odge
  {T}heaters.
\newblock \emph{Publ. Res. Inst. Math. Sci.}, 57\penalty0 (1/2):\penalty0
  3--207, 2021{\natexlab{a}}.
\newblock URL \url{https://ems.press/journals/prims/articles/201525}.

\bibitem[Mochizuki(2021{\natexlab{b}})]{mochizuki-iut2}
Shinichi Mochizuki.
\newblock Inter-{U}niversal {T}eichmuller {T}heory {II}: {H}odge-{A}rakelov
  {T}heoretic {E}valuations.
\newblock \emph{Publ. Res. Inst. Math. Sci.}, 57\penalty0 (1/2):\penalty0
  209--401, 2021{\natexlab{b}}.
\newblock URL \url{https://ems.press/journals/prims/articles/201526}.

\bibitem[Mochizuki(2021{\natexlab{c}})]{mochizuki-iut3}
Shinichi Mochizuki.
\newblock Inter-{U}niversal {T}eichmuller {T}heory {III}: canonical splittings
  of the log-theta lattice.
\newblock \emph{Publ. Res. Inst. Math. Sci.}, 57\penalty0 (1/2):\penalty0
  403--626, 2021{\natexlab{c}}.
\newblock URL \url{https://ems.press/journals/prims/articles/201527}.

\bibitem[Mochizuki(2021{\natexlab{d}})]{mochizuki-iut4}
Shinichi Mochizuki.
\newblock Inter-{U}niversal {T}eichmuller {T}heory {IV}: Log-volume
  computations and set theoretic foundations.
\newblock \emph{Publ. Res. Inst. Math. Sci.}, 57\penalty0 (1/2):\penalty0
  627--723, 2021{\natexlab{d}}.
\newblock URL \url{https://ems.press/journals/prims/articles/201528}.

\bibitem[Mochizuki(2022)]{mochizuki-essential-logic}
Shinichi Mochizuki.
\newblock On the essential logical structure of {I}nter-{U}niversal
  {T}eichmuller {T}heory in terms of logical and “$\wedge$”/logical or
  “$\vee$” relations: Report on the occasion of the publication of the four
  main papers on {I}nter-{U}niversal {T}eichmuller {T}heory.
\newblock 2022.
\newblock URL
  \url{https://www.kurims.kyoto-u.ac.jp/~motizuki/Essential\%20Logical\%20Structure\%20of\%20Inter-universal\%20Teichmuller\%20Theory.pdf}.

\bibitem[Oesterle(1988)]{oesterle1988}
J.~Oesterle.
\newblock Nouvelles approches du ``th\'eor\`eme'' de {F}ermat.
\newblock \emph{Asterisque}, 161:\penalty0 165--186, 1988.

\bibitem[Poonen(1993)]{poonen93}
Bjorn Poonen.
\newblock Maximally complete fields.
\newblock \emph{Enseign. Math. (2)}, 39\penalty0 (1-2):\penalty0 87--106, 1993.

\bibitem[Schmidt(1933)]{schmidt33}
F.K. Schmidt.
\newblock Mehrfach perfekte korper.
\newblock \emph{Math. Annalen}, 108\penalty0 (1):\penalty0 1--25, 1933.

\bibitem[Scholze(2012)]{scholze12-perfectoid-ihes}
Peter Scholze.
\newblock Perfectoid spaces.
\newblock \emph{Publ. Math. Inst. Hautes \'{E}tudes Sci.}, 116:\penalty0
  245--313, 2012.
\newblock URL \url{https://doi.org/10.1007/s10240-012-0042-x}.

\bibitem[Scholze(2017)]{scholze-diamonds}
Peter Scholze.
\newblock \emph{\'Etale cohomology of Diamonds}.
\newblock 2017.
\newblock URL \url{https://arxiv.org/abs/1709.07343}.

\bibitem[Scholze(2021)]{scholze-review}
Peter Scholze.
\newblock Review of {M}ochizuki's paper: {I}nter-{U}niversal {T}eichmüller
  {T}heory. {I},{II},{III},{IV}.
\newblock \emph{zbMath Open (formerly Zentralblatt Math)}, 2021.
\newblock URL \url{https://zbmath.org/pdf/07317908.pdf}.

\bibitem[Scholze and Stix(2018)]{scholze-stix}
Peter Scholze and Jakob Stix.
\newblock Why $abc$ is still a conjecture.
\newblock 2018.
\newblock URL
  \url{https://www.math.uni-bonn.de/people/scholze/WhyABCisStillaConjecture.pdf}.

\bibitem[Serre(1979)]{serre1979-local-fields}
J.-P. Serre.
\newblock \emph{Local fields}, volume~67 of \emph{Graduate {T}exts in
  {M}athematics}.
\newblock Springer, Berlin, Springer-{V}erlag, 1979.

\bibitem[Silverman(1985)]{silverman-arithmetic}
Joseph Silverman.
\newblock \emph{The arithmetic of elliptic curves}, volume 106 of
  \emph{Graduate {T}ext in {Mathematics}}.
\newblock Springer-{V}erlag, Berlin, 1985.

\bibitem[{S}ilverman(1986)]{silverman1986}
Joseph~{H}. {S}ilverman.
\newblock \emph{The Theory of Height Functions}, page 151–166.
\newblock Springer New York, New York, NY, 1986.
\newblock \doi{10.1007/978-1-4613-8655-1_6}.

\bibitem[Stewart and Tijdeman(1986)]{stewart1986}
C.~L. Stewart and R.~Tijdeman.
\newblock On the {O}esterl\'{e}-{M}asser conjecture.
\newblock \emph{Monatsh. Math.}, 102\penalty0 (3):\penalty0 251--257, 1986.
\newblock URL \url{https://doi.org/10.1007/BF01294603}.

\bibitem[Szpiro(1979)]{szpiro-book1979}
Lucien Szpiro.
\newblock \emph{Lectures on equations defining space curves (Notes by N. Mohan
  Kumar)}.
\newblock Tata Institute of Fundamental Research; Narosa Pub. House, Bombay,
  New Delhi, 1979.

\bibitem[Szpiro(1981)]{szpiro81}
Lucien Szpiro.
\newblock S\'eminaire sur les pinceaux de courbes de genre au moins deux.
\newblock \emph{Ast\'erisque}, 86, 1981.

\bibitem[Szpiro(1990)]{szpiro1991}
Lucien Szpiro.
\newblock Discriminant et conducteur des courbes elliptiques, in s\'minaire sur
  les pinceaux de courbes elliptiques (\`a la recherche de mordell effectif),.
\newblock \emph{Ast\'erisque, no. 183 (1990), pp. 7-18}, 183:\penalty0 7--18,
  1990.
\newblock URL \url{http://www.numdam.org/item/AST_1990__183__7_0/}.

\bibitem[Tan(2018)]{fucheng}
Fucheng Tan.
\newblock Note on {I}nter-{U}niversal {T}eichmuller {T}heory.
\newblock 2018.

\bibitem[{v}an {F}rankenhuysen(2002)]{frankenhuysen2002}
Machiel {v}an {F}rankenhuysen.
\newblock The abc conjecture implies vojta's height inequality for curves.
\newblock \emph{Journal of number theory}, 95\penalty0 (2):\penalty0 289--302,
  2002.
\newblock \doi{10.1016/S0022-314X(01)92769-6}.

\bibitem[Vojta(1998)]{vojta1998}
Paul Vojta.
\newblock A more general abc conjecture.
\newblock \emph{International Mathematics Research Notices}, 1998\penalty0
  (21):\penalty0 1103--1116, 1998.
\newblock \doi{10.1155/S1073792898000658}.

\bibitem[Wada(2021)]{wada2016}
Yuki Wada.
\newblock Near miss {\it abc}-triples in compactly bounded subsets.
\newblock In \emph{Inter-universal {T}eichm\"uller {T}heory {S}ummit 2016},
  volume B84 of \emph{RIMS K\^oky\^uroku Bessatsu}, pages 193--227. Res. Inst.
  Math. Sci. (RIMS), Kyoto, 2021.

\bibitem[Wright(2019)]{wright19}
Alex Wright.
\newblock A tour through {M}irzakhani's work on moduli spaces of {R}iemann
  surfaces, 2019.
\newblock URL \url{https://arxiv.org/abs/1905.01753}.

\bibitem[Yamashita(2019)]{yamashita}
Go~Yamashita.
\newblock A proof of the abc-conjecture.
\newblock \emph{RIMS--Kokyuroku Bessatsu}, 2019.

\bibitem[Zhang(2001)]{zhang01}
Shouwu Zhang.
\newblock Geometry of algebraic points.
\newblock In \emph{First {I}nternational {C}ongress of {C}hinese
  {M}athematicians ({B}eijing, 1998)}, volume~20 of \emph{AMS/IP Stud. Adv.
  Math.}, pages 185--198. Amer. Math. Soc., Providence, RI, 2001.
\newblock URL \url{https://doi.org/10.2298/pim140921001z}.

\end{thebibliography}

\begin{abstract}
This is a continuation of my work on the theory of Arithmetic Teichmuller Spaces. The theory developed in the present series of papers leads to the proof $abc$-conjecture following Mochizuki's rubric in \iut.
\end{abstract}

\numberwithin{equation}{subsection}
\newcommand{\omu}{\O_{\bQ_p}^{\mu}}
\newcommand{\lmod}{L_{\rm mod}}
\newcommand{\ttheta}{{{\widetilde\Theta}}}
\newcommand{\tthetaj}[1]{\ttheta_{\scriptscriptstyle{Joshi},#1}}
\newcommand{\tthetam}[1]{\ttheta_{\scriptscriptstyle{Mochizuki},#1}}

\newcommand{\moccor}{\cite[Corollary 3.12]{mochizuki-iut3}}
\newcommand{\thetam}{{\ttheta}_{\scriptscriptstyle{Mochizuki}}}
\newcommand{\thetamp}{{\ttheta}_{\scriptscriptstyle{Mochizuki},p}}
\newcommand{\thetampi}{{\ttheta}_{\scriptscriptstyle{Mochizuki},p}^{\bsI}}
\newcommand{\thetaminfty}{{\ttheta}_{\scriptscriptstyle{Mochizuki},\infty}^{\bsI}}
\newcommand{\thetaj}{{{\ttheta}_{\scriptscriptstyle{Joshi}}}}
\newcommand{\thetajp}{{{\ttheta}_{\scriptscriptstyle{Joshi},p}}}
\newcommand{\sM}{\mathscr{M}}
\newcommand{\sMb}{\overline{\mathscr{M}}}
\newcommand{\pib}{\overline{\Pi}}
\newcommand{\bN}{\mathbb{N}}
\newcommand{\sD}{\mathscr{D}}
\newcommand{\sF}{\mathscr{F}}
\newcommand{\sL}{\mathscr{L}}

\newcommand{\bdrp}{B_{dR}^+}
\newcommand{\bdr}{B_{dR}}

\newcommand{\iutthr}{\cite{mochizuki-iut1,mochizuki-iut2, mochizuki-iut3}}

\newcommand{\thetajpp}{{\widehat{\Theta}}_{\scriptscriptstyle{Joshi},p}}
\newcommand{\thetaja}{{\widehat{\Theta}}_{\scriptscriptstyle{Joshi}}}
\newcommand{\thetajpph}{{\widehat{\widehat{\Theta}}}_{\scriptscriptstyle{Joshi},p}}
\newcommand{\thetajppa}{{\widehat{\widehat{\Theta}}}_{\scriptscriptstyle{Joshi}}}

\newcommand{\ells}{{\ell^*}}
\newcommand{\sP}{\mathscr{P}}
\newcommand{\spc}{\sP^{Teich}}
\newcommand{\pitop}[1]{\pi_1^{top}(#1)}
\newcommand{\sppi}{\sP(\Pi\hookleftarrow\pib)}
\newcommand{\sppim}{\sP_{\scriptscriptstyle{Mochizuki}}(\Pi\hookleftarrow\pib)}
\newcommand{\sppij}{\sP_{\scriptscriptstyle{Joshi}}(\Pi\hookleftarrow\pib)}

\newcommand{\bcris}{B_{cris}}

\newcommand{\hol}[3]{\mathfrak{hol}_{#1}(#3)_{#2}}

\newcommand{\xdia}{X^\lozenge}
\newcommand{\holdia}{{\mathfrak{Hol}(X/E)}}
\newcommand{\spd}{{\rm Spd}}
\newcommand{\perf}{{\rm Perf}}
\newcommand{\perfld}{{\rm PerfFld}}
\newcommand{\spa}{{\rm Spa}}
\newcommand{\bnddsub}{{\tiny\circ}}
\newcommand{\xad}{X^{ad}}
\newcommand{\sppimold}{\sP_{\scriptscriptstyle{Mochizuki}}(\Pi\hookleftarrow\pib)}
\renewcommand{\sppim}{\sP_{\scriptscriptstyle{Mochizuki}}'(\Pi\hookleftarrow\pib)}

\newcommand{\chxfe}{{\rm Ch^1(\sxfe)}}
\newcommand{\divxfe}{{\rm Div}(\sxfe)}

\newcommand{\qprojqp}[1]{\mathcal{QProj}_{#1}}
\newcommand{\projqp}[1]{\mathcal{Proj}_{#1}}

\newcommand{\tsigma}{\tilde{\sigma}}
\newcommand{\ismg}{{\mathbf{Ism}(G)}}
\renewcommand{\tripod}{\P^1-\{0,1,\infty\}}
\newcommand{\okbt}{\O_{\bE;K}^\triangleright}
\newcommand{\Ob}{\overline{\O}} 

\renewcommand{\deg}{{\rm deg}}
\newcommand{\udeg}[1]{\underline{{\rm deg}}(#1)}
\newcommand{\divp}{Div_{+}}
\newcommand{\divt}{\widetilde{{\rm Div_{+}}}}
\newcommand{\pdivt}{\widetilde{{\rm PDiv}}}
\newcommand{\fjxlp}{\fJ(X,L_\wp)}

\newcommand{\bbvl}{\mathbb{V}_L}
\newcommand{\bbvlp}{\mathbb{V}_{L'}}
\newcommand{\bbvlmod}{\mathbb{V}_{\lmod}}
\newcommand{\bbvlmodgood}{\mathbb{V}_{\lmod}^{good}}
\newcommand{\bbvlgood}{\mathbb{V}_{L}^{good}}
\newcommand{\bbvlpgood}{\mathbb{V}_{L'}^{good}}
\newcommand{\bbvlmodoss}{\mathbb{V}_{\lmod}^{odd,ss}}
\newcommand{\bbvloss}{\mathbb{V}_{L}^{odd,ss}}
\newcommand{\bbvlposs}{\mathbb{V}_{L'}^{odd,ss}}
\newcommand{\ubblv}{\underline{\mathbb{V}}}
\newcommand{\ubblvgood}{\underline{\mathbb{V}}^{good}}
\newcommand{\ubblvoss}{\underline{\mathbb{V}}^{odd,ss}}
\newcommand{\ubblvossp}{\underline{\mathbb{V}}^{odd,ss}_p}

\newcommand{\ul}[1]{\underline{#1}}

\newcommand{\breveB}{{\breve{B}}}
\newcommand{\breveBb}{{\breve{\mathbf{B}}}}
\newcommand{\bigB}{\bigotimes_{E'\anabelmap E} \left(\oplus_{j=1}^\ells B_{E'} \right)}
\newcommand{\tSigma}{\widetilde{\Sigma}}
\newcommand{\thetajh}{{\widehat{\Theta}}_{\scriptscriptstyle{Joshi}}}

\newcommand{\sI}{\mathcal{I}}
\newcommand{\utheta}{U_\Theta}

\newcommand{\logsh}[1]{\mathscr{I}({#1})}
\newcommand{\logshq}[1]{\mathscr{I}^{\Q_p}({#1})}
\newcommand{\logshe}{\logsh{E}}
\newcommand{\logsheq}{\logshq{E}}
\newcommand{\logshep}{\logsh{E'}}
\newcommand{\logshepq}{\logshq{E'}}

\newcommand{\bS}{\mathbb{S}}

\newcommand{\sIp}{\sI_p}
\newcommand{\sIpq}{{\sIp}^{\Q}}
\newcommand{\sIm}{\sI_{\scriptscriptstyle{Mochizuki}}}
\newcommand{\sImq}{\sI_{\scriptscriptstyle{Mochizuki}}^{\Q}}
\newcommand{\mulog}[2]{\mu_{#1}^{\log}({#2})}
\newcommand{\EB}{{}^EB} 

\newcommand{\vnon}{\V^{non}}
\newcommand{\varc}{\V^{arc}}

\newcommand{\dmod}{d_{mod}}
\newcommand{\dmods}{d^*_{mod}}
\newcommand{\emod}{e_{mod}}
\newcommand{\emods}{e^*_{mod}}
\newcommand{\ltpd}{L_{tpd}}
\newcommand{\ff}{\mathfrak{f}}
\newcommand{\fd}{\mathfrak{d}}
\newcommand{\fq}{\mathfrak{q}}
\newcommand{\fQ}{\mathfrak{Q}}
\newcommand{\fT}{\mathfrak{Tate}}
\newcommand{\fl}{\mathfrak{l}}
\newcommand{\dtpd}{\fd_{\ltpd}}
\newcommand{\ftpd}{\ff_{\ltpd}}
\newcommand{\etap}{\eta_{prm}}
\newcommand{\V}{\mathbb{V}}
\newcommand{\vq}{\V_{\Q}}
\newcommand{\vnl}{\V_L^{non}}
\newcommand{\vlnon}{\vnl}
\newcommand{\vlarch}{\mathbb{V}_L^{arc}}
\newcommand{\vlp}{\V_{L'}}
\newcommand{\vlpp}{\V_{L',p}}
\newcommand{\vnlp}{\V_{L'}^{non}}
\newcommand{\vlnonp}{\vnlp}
\newcommand{\vlarchp}{\mathbb{V}_L^{'arc}}
\newcommand{\fA}{\mathfrak{A}}

\newcommand{\bsI}{\boldsymbol{\mathscrbf{I}}}
\newcommand{\thetajti}{{\thetaj}^{\tilde{\mathscrbf{I}}}}
\newcommand{\thetaji}{{\thetaj}^{\mathscrbf{I}}}
\newcommand{\thetamti}{{\thetam}^{\tilde{\mathscrbf{I}}}}
\newcommand{\thetami}{{\thetam}^{\mathscrbf{I}}}
\newcommand{\flog}{\mathfrak{log}}

\newcommand{\lmodv}{L_{mod,v}}
\newcommand{\vmdst}{\V_M^{dst}}
\newcommand{\vqdst}{\V_{\Q}^{dst}}
\newcommand{\vmnon}{\vnon_M}
\newcommand{\supp}{{\rm Supp}}

\newpage
{\sffamily
	\tableofcontents
}

\newpage

\definecolor{darkmidnightblue}{rgb}{0.0, 0.2, 0.4}
\definecolor{carmine}{rgb}{0.59, 0.0, 0.09}
\newtcbox{\mybox}[1][red]{on line,
	arc=0pt,outer arc=0pt,colback=white,colframe=darkmidnightblue,
	boxsep=0pt,left=0pt,right=0pt,top=2pt,bottom=2pt,
	boxrule=0pt,leftrule=1pt, rightrule=1pt,bottomrule=1pt,toprule=1pt}
\newcommand{\tmb}[1]{\mybox{#1}}
\newcommand{\present}{present series of papers (\cite{joshi-anabelomorphy,joshi-untilts,joshi-teich-estimates,joshi-teich-def,joshi-teich-rosetta}  and \cite{joshi-anabelomorphy,joshi-formal-groups})}

\newcommand{\abc}{$abc$-conjecture}
\epigraph{Discovery is the privilege of the child: the child who has no fear of being once again wrong, of looking like an idiot, of not being serious, of not doing things like everyone else.}{Alexander Grothendieck}
\section{Introduction}\label{se:intro}\nwss

\newcommand{\vl}{\V_L}
\newcommand{\lvbh}{\widehat{\overline{L}}_v}
\newcommand{\lvbht}{{\widehat{\overline{L}}_v}^\flat}
\newcommand{\lvbhmax}{\widehat{\overline{L}}_v^{max}}
\newcommand{\lvbhtmax}{{\widehat{\overline{L}}_v}^{max\flat}}
\newcommand{\lvbhreal}{\widehat{\overline{L}}_v^{\R}}
\newcommand{\lvbhtreal}{{\widehat{\overline{L}}_v}^{\R\flat}}

\newcommand{\arith}[1]{\mathfrak{arith}(#1)}
\newcommand{\arithl}{\arith{L}}
\newcommand{\adel}[1]{\mathfrak{adel}(#1)}
\newcommand{\adell}{\adel{L}}
\newcommand{\by}{{\bf y}}

\newcommand{\syflv}{\sY_{\lvbht,L_v}}
\newcommand{\syflvmax}{\sY_{\lvbhtmax,L_v}}
\newcommand{\sxflv}{\sX_{\lvbht,L_v}}
\newcommand{\sxflvmax}{\sX_{\lvbhtmax,L_v}}
\newcommand{\bsY}{\mathscrbf{Y}}
\newcommand{\bsX}{\mathscrbf{X}}
\newcommand{\yadl}{\bsY_L}
\newcommand{\yadlp}{\bsY_{L'}}
\newcommand{\yadlmax}{\bsY_L^{max}}
\newcommand{\yadq}{\bsY_{\Q}}
\newcommand{\yadqmax}{\bsY_{\Q}^{max}}
\newcommand{\xadl}{\bsX_L}
\newcommand{\xadq}{\bsX_{\Q}}
\newcommand{\xadlmax}{\bsX_L^{max}}
\newcommand{\xadqmax}{\bsX_{\Q}^{max}}
\newcommand{\yadlpoint}{\{(L_\wp\into K_{\wp}, K_\wp\isom \cpt)\}_{\wp\in\vlnon}}
\newcommand{\xadlpoint}{\{(L_\wp\into K_{\wp}, K_\wp\isom \cpt)\}_{\wp\in\vlnon}}

\newcommand{\bv}{\bar{v}}
\newcommand{\sGal}[1]{\mathbf{G}_{#1}}

\newcommand{\fet}{\mathcal{F\hat{e}t}}
\newcommand{\sA}{\mathcal{A}}

\newcommand{\logp}{\log^+}
\newcommand{\sFrob}{\mathcal{Frob}}
\newcommand{\pow}[2]{#1\llbracket#2\rrbracket}

\newcommand{\ainf}{A_{\textrm{inf}}}

\newcommand{\xs}{X^S}
\newcommand{\sZ}{\mathscr{Z}}
\newcommand{\bvphi}{\boldsymbol{\varphi}}

\subsection{The goal of this paper} 
[\textcolor{red}{This is still a work in progress!}] This paper completes (in \Cref{th:main-dioph-thm}) the remarkable proof of the \abc\ announced by Shinichi Mochizuki in \iut\  (see \cref{ss:outline} for an outline). 

Unfortunately, \iut\ is quite opaque and labyrinthine. In the decade between the release of Mochizuki's preprints in 2012 and its publication in 2021, Mochizuki's claims  have remained poorly understood; and following the release of \cite{scholze-stix} (and \cite{scholze-review}), many came to the conclusion that Mochizuki's theory was vacuous. 

My work, \cite{joshi-formal-groups}, and the \present, commenced amidst this backdrop, provides an independent construction of the  Theory of Arithmetic Teichmuller Spaces which includes, as a special case, Mochizuki's Theory and shows that many impressions readers may previously  have had about Mochizuki's Theory must be reconsidered (\cite{joshi-report}). 

The most important discovery of \iutthr\ is that the arithmetic of a fixed number field is itself deformable, there is a Frobenius morphism of a number field (known as the $\flog$-Link in \iutthr), and, in fact, there is a Teichmuller Theory of a fixed Number Field\footnote{I use the capitalization Number Field here in the same sense as the usage of the capitalization Banach Space in \cite{colmez2000} i.e. because \cite{joshi-teich-def} provides a category equipped with a functor to infinite dimensional $\Q_p$-Banach spaces for each prime  number $p$.}. Each of these assertions, especially the assertion that a number field is deformable (just as a Riemann surface is deformable), has turned out to be an extremely subtle and challenging one, and to be sure, Mochizuki does not provide any quantification of this, even though his theory rests upon it. 

One of the original contributions of the present series of papers is to provide a precise quantification, of the idea that an apparently rigid object such as a number field has topological deformations, and also of the other assertions mentioned above. The above assertions are now robustly established in \cite{joshi-teich-def}, where the notion of a deformation of  number field is captured in the notion of an \textit{arithmeticoid}. Each distinct arithmeticoid presents a quantifiably different avatar of a fixed number field, there exists topologically non-isomorphic arithmeticoids and the deformability is established by constructing, in analogy with classical Teichmuller Theory,  a metrisable topological space which parameterizes arithmeticoids (so it makes sense to talk of proximate arithmeticoids). 

At any rate, the importance of Mochizuki's assertion should not be underestimated: it turns on its head the absolute rigidity of Arithmetic which has dominated number theory since the introduction of number fields at the hands of Gauss, Kummer and Dedekind and it also brings the analogy between Riemann Surfaces and number fields to yet another level by bringing it on par with Teichmuller's refinement of Riemann's Moduli Theory. 

In particular, the point of view of \iutthr\ is that one can now treat the Arithmetic of a fixed number field as a dynamical variable and average quantities of Diophantine interest over aggregates of Arithmetics of a fixed number field. The breath-taking novelty and the scope of this idea should not be lost on us. 

As an exhibit of what this does for us, Mochizuki offers a proof of the \abc\ in \cite{mochizuki-iut4}. While the proof  of the \abc\ (see \cref{ss:outline} for an outline) described in this paper  is based on Mochizuki's proof in \cite{mochizuki-iut4}, my approach to it is based on the \present\ and at this juncture, readers need not have any familiarity with \iut\ to understand the  claims established in the \present. However, for readers who wish to compare the two theories, I have provided  detailed references to \iut\  and a `Rosetta Stone' for transliterating the assertions of the two theories. 

I have also organized the proof (given in \cite{mochizuki-iut4}) from my point of view for greater transparency,  filled in details and provided additional proofs (for \cite{mochizuki-iut4} and \cite{mochizuki-general-pos}) which I felt were necessary for completeness. My approach to the proof of the \abc\  demonstrates  quantitatively the idea of  ``comparing and averaging over distinct arithmeticoids or avatars of the familiar arithmetic of integers'' (see \cite{joshi-teich-def}, \cite{joshi-teich-rosetta}). 

Mochizuki's vision (in \iutthr) is quite different from mine: he has viewed his proof as changing of mathematical universes and comparing arithmetic in  independent set theories (\cite[Page 26]{mochizuki-iut1}) and hence the name \textit{Inter-Universal Teichmuller Theory}. [Readers should bear in mind that unfortunately, Mochizuki does not offer a precise quantification of this assertion in his papers.] My theory, on the other hand, is rooted in a fixed mathematical universe,  within which one can quantify precisely the existence of topologically distinct arithmetics of integers (pluralized `arithmetics' is intended) and the inter-relationships between them which bears all the features of classical Teichmuller Theory (hence its name: \textit{The Theory of Arithmetic Teichmuller Spaces}). My work allows one to precisely quantify Mochizuki's ideas as well. 

The theory presented in my work so closely  mirrors Classical Teichmuller Theory that it is tempting to dream that it will allow us to weave, in Arithmetic Geometry, a tapestry of ideas and theorems, as rich as the one that has emerged in the last hundred years with the invention of classical Teichmuller Theory.

\subsection{Why is this paper needed?}\label{ss:why-this-is-needed}
My role in this complicated saga (see \cite{castelvecchi-1}, \cite{castelvecchi-2}) surrounding \iut\ has been to address the foundational and existential questions which have engulfed \iutthr. This has allowed me to establish a complete and independent proof of \moccor\ which is central to the claims of \cite{mochizuki-iut4}. 
 
It is clear to me that there is a valid basis for concern about  \iut, Mochizuki's methodology, and the approach  he has adopted to persuade his readership (e.g. \cite{mochizuki-essential-logic}). However, it must also be emphasized that limited experience, among the arithmetic geometry and the anabelian geometry communities, with  classical Teichmuller Theory has played a central role in the unfortunate debacle that has surrounded \iut. 

Let me discuss some important examples of issues with \iut. The central difficulty one faces in \iut\ is how to disambiguate objects within the isomorphism class of, say, a smooth hyperbolic curve over a $p$-adic field. In the archimedean  case, Classical Teichmuller Theory of Riemann surfaces solves this problem by using quasi-conformal mappings between Riemann surfaces in place of conformal mappings used in Riemann's moduli theory. While Mochizuki constructed his theory in analogy with classical Teichmuller Theory (this is why Teichmuller Theory appears in the title of Mochizuki's papers), unfortunately \iut\ does not provide any way of disambiguating objects  used in his theory. This difficulty was first highlighted by \cite{scholze-stix} (\textit{to be clear, my work shows that the theory does need distinct objects}). They declared the impossibility of disambiguation, but my work shows that this is incorrect. Mochizuki, on the other hand, dismissed the disambiguation problem altogether in \cite{mochizuki-essential-logic}. [My discussion of \cite{scholze-stix}, \cite{scholze-review} is in \cref{the scholze-stix-report}.]

Here is another important example of inadequacies of \iut. Some anabelian geometers (Mochizuki included) have asserted that the collection of algebraically closed complete valued fields $\{\C_p \}_{p\in\mathbb{V}_{\Q}}$ is adequate for establishing \iutthr. This claim is incorrect because
\benumlab
\item Mochizuki's Key Principle of Inter-Universality \cite[\ssep I3]{mochizuki-iut1} (which according to Mochizuki lies at the foundation of \iut) requires one to work with arbitrary geometric base-points for tempered fundamental groups, 
\item   the definition of geometric base-points for tempered fundamental groups given in \cite[Section 4.1]{andre03} (Andre's paper lays the foundation of the theory of tempered fundamental groups), and
\item  \constrone{Lemma }{lem:perfectoid},
\eenum
taken together, imply that arbitrary algebraically closed perfectoid fields naturally enter \iut. On the other hand, as discussed in \cite{joshi-report}, \cite{joshi-teich-rosetta}, arbitrary algebraically closed perfectoid fields cannot be recovered using Mochizuki's group-theoretic methods of \topics\ (on which \iut\ is based) and this leads to a whole range of mathematical issues in \iut. 

In late 2019--early 2020, I came to the conclusion (which is documented in \cite{joshi-untilts-2020}--current version is \cite{joshi-untilts}) that the introduction of algebraically closed perfectoid fields provides a  canonical solution to the disambiguation problem which is fully  compatible with Mochizuki's anabelian ideas, his Key Principle of Inter-Universality  and also fully compatible with our understanding of the role of analytic functions in the classical case, thus making it clear that there  exists a canonical Teichmuller Theory, in a stronger sense than has been asserted by Mochizuki, even in the $p$-adic context.

A more central and rather troublesome difficulty  encountered with Mochizuki's  approach using Frobenioids (or monoidal data) is this: it is exceedingly difficult to mathematically establish, using only the methods of \topics\ and \iut, that modifications to a Frobenioid or a monoid necessarily arise from an intrinsic variation of ring structures. Much of the criticism, by Scholze-Stix, myself and others, of \iut\ has been about the lack of precise answers to such technical questions. On the other hand, my work makes it clear from the very beginning how modifications of monoidal data lead to distinct arithmetic holomorphic structures (see \cite{joshi-formal-groups}, \cite{joshi-teich}, \cite{joshi-teich-def}).

Another important question regarding \iutthr\ is what is the significance of Mochizuki's Three Indeterminacies  (\cite[Pages 576--577]{mochizuki-iut3})  which are central to the formulation of \cite[Theorem 3.11]{mochizuki-iut3}  and hence to \moccor.  My work provides (see \constrthr{\ssep}{ss:Mochizuki-indeterminacies}) a canonical geometric description of Mochizuki's Three Fundamental Indeterminacies Ind1, Ind2, Ind3--in particular this eliminates  doubts about their existence asserted in \cite{scholze-stix}. 

There are more subtle phenomena which were discovered or claimed by Mochizuki (e.g. Galois action is used to change arithmetic holomorphic structures in \cite[Theorem 3.11]{mochizuki-iut3}), but whose proofs are now available \textit{only} through the \present. Notably, by \cite{joshi-teich-def}, there exists a Teichmuller Theory of number fields as asserted by Mochizuki.

Hence I believe that, without addressing the concerns raised about \iutthr\ and without their resolution as presented in my work, \iut\ and its proof of the \abc\ unfortunately stands  incomplete. 

The \present\ detail a solution to the issues raised about \iut, which is fully compatible with Mochizuki's key ideas and especially, fully compatible with his Key Principle of Inter-Universality \cite[\ssep I3]{mochizuki-iut1}  (for a detailed discussion see the `Rosetta Stone' provided in \cite[{\ssep\ref{III-se:intro-rosetta-stone}}]{joshi-teich-rosetta}). The solution which I have found, independently of Mochizuki's ideas, goes well beyond the original confines of the theory presented in \iutthr, and extends to higher genus and even higher dimensions, and provides us with a new and radical way of thinking about arithmetic of number fields. 

The following table provides a convenient reference point for the key concepts and properties which need to be established  in \iutthr\ to arrive at complete and satisfactory proofs of Mochizuki's claims (and hence \cite{mochizuki-iut4}): 

\begin{boxedcontent}[label={tab:quickref},grow to left by=0.5cm, grow to right by=0.5cm]{Quick Reference}{}
\begin{center}
	\begin{table}[H]
		\renewcommand{\arraystretch}{1.5}
		\resizebox{1.01\textwidth}{!}{\begin{tabular}{|c|c|}
				\hline 
				{\bf Concepts of Mochizuki's Theory which need to be established} & {\bf Reference in the present series of papers} \\
				\hline 
				Anabelomorphisms (I introduced this term) & \cite[{\ssep\ref{0-se:anabelomorphy}}]{joshi-anabelomorphy}  \\ 
				\hline 
				Mochizuki's Arithmetic Holomorphic Structures	& \cite[{Def. \ref{U-def:arith-hol-strs}}]{joshi-untilts}, (also \cite{joshi-untilts-2020}) \\ 
				\hline 
				Mochizuki's $\Theta_{gau}$-Links	& \cite[{\ssep\ref{II-se:ansatz}}]{joshi-teich-estimates} (local), \cite[{\ssep\ref{III-se:adelic-ansatz}}]{joshi-teich-rosetta} (global)  \\ 
				\hline 
				Mochizuki's Log-Links (Mochizuki's proxy for global Frobenius)	& \cite[{Thm. \ref{II-th:flog-links}}]{joshi-teich-estimates} (local), \cite[{Cor. \ref{II5-cor:global-frobenius}}]{joshi-teich-def} (global)  \\ 
				\hline 
				Mochizuki's Theta-Pilot Objects	& \cite[{Thm. \ref{II-thm:theta-pilot-object-appears}}]{joshi-teich-estimates} (local), \cite[{Def. \ref{III-def:additive-theta-pilot-object}}]{joshi-teich-rosetta} (global)  \\ 
				\hline
				Mochizuki's Theta values locus $\thetami$	& \cite[{\ssep\ref{II-se:construction-ttheta}}]{joshi-teich-estimates} (local), \cite[{\ssep\ref{III-ss:def-thetam}}]{joshi-teich-rosetta} (global)  \\ 
				\hline 
				Mochizuki's Three Indeterminacies	&  \cite[{\ssep\ref{III-ss:Mochizuki-indeterminacies}}]{joshi-teich-rosetta} (global), \cite[{\ssep\ref{U-ss:replacements-iut}}]{joshi-untilts}, (also \cite{joshi-teich})   \\ 
				\hline 
				Mochizuki's Prime-Strips	&  \cite[{\ssep\ref{III-th:gluing-prime-strips}}]{joshi-teich-rosetta} (global)  \\ 
				\hline 
				Mochizuki's Hodge Theaters	&  \cite[{\ssep\ref{III-ss:hodge-theater}}]{joshi-teich-rosetta} (global)  \\ 
				\hline 
				Mochizuki's Frobenioids	&  \cite[{\ssep\ref{III-th:gluing-frobenioids}}]{joshi-teich-rosetta} (global)  \\ 
				\hline 
			\end{tabular} 
		}
	\end{table}
\end{center}
\end{boxedcontent}

\subsection{Archival Note} [A detailed time-line of events is given in \cite{joshi-report}.]
Starting with my sabbatical visit to RIMS, Kyoto (Spring 2018), Mochizuki and I actively discussed matters related to \iut\ during 2018--2020. Unfortunately, Mochizuki stopped corresponding with me after the release of \cite{joshi-untilts-2020}. I reached out to him on a number of occasions during 2020--2024 but my emails went unanswered.

After  \cite{joshi-teich-rosetta} was released on the arxiv, I wrote to Mochizuki suggesting that it would be best that we write this paper as a joint paper, revisiting \cite{mochizuki-iut4} and settle, once and for all, the issues which have surrounded his remarkable proof (once again, I offered to explain my work in Zoom meetings/lectures). But neither Mochizuki nor the other anabelian geometers who have been involved with \iut,  responded to my e-mails.
\subsection{Suggested reading order and logical dependency of various results}
This paper depends on results contained in other papers in this series (with links to the arxiv versions):
\begin{description}
\item[{\href{https://arxiv.org/pdf/1906.06840}{Formal Groups} \cite{joshi-formal-groups}}] This paper establishes the key relationship between perfectoid geometry and Mochizuki's multiplicative monoidal point of view and is a recommended reading even though its results are not directly used in the present series of papers.
\item[{\href{https://arxiv.org/pdf/2003.01890}{Anabelomorphy} \cite{joshi-anabelomorphy}}] I will freely use the notion of anabelomorphy and related set of ideas introduced therein, but no other results from it are required in the present paper.
\item[{\href{https://arxiv.org/pdf/2106.11452}{Constructions I} \cite{joshi-teich}}] This shows that the absolute Grothendieck Conjecture is false in the category of Berkovich spaces (just as it is false in the category of complex analytic spaces and specifically for Riemann surfaces) and uses this to formulate a Theory of Arithmetic Teichmuller Spaces  which parallels classical Teichmuller Theory. 

\item[{\href{https://arxiv.org/pdf/2008.01228}{Abs. Grothendieck Conj. and $p$-adic Hodge Theory} \cite{joshi-gconj}}]
This shows that the Absolute Grothendieck Conjecture is False for Fargues-Fontaine curves. These curves play a central role in my approach to Arithmetic Holomorphic Structures and the main assertion of \cite{joshi-gconj} provides an intrinsic geometric approach to the existence of Mochizuki's Indeterminacy Ind1 in \iut.

\item[{\href{https://arxiv.org/pdf/2210.11635}{Untilts Paper} \cite{joshi-untilts}}] This paper formulates the precise notion of Arithmetic Holomorphic Structures and provides  the basic properties of (local and global) arithmetic Teichmuller spaces arising from a quasi-projective variety over a $p$-adic field or a number field.
\item[{\href{https://arxiv.org/pdf/2303.01662}{Constructions II} \cite{joshi-teich-estimates}}] Formulates and proves a  local prototype of Mochizuki's Corollary 3.12. The key discovery of this paper is a canonical geometric description of Mochizuki's $\flog$-Links and Mochizuki's $\Theta_{gau}$-Links using the theory of untilts of perfectoid fields. Notably, this approach establishes certain local correspondences and the set of such correspondences on Arithmetic Teichmuller spaces. 
\item[{\href{https://arxiv.org/pdf/2305.10398}{Constructions II(1/2)} \cite{joshi-teich-def}}] This deals with the Global theory of Arithmetic Teichmuller space of a fixed number field, and notably, establishes  the existence of (global) deformations of a fixed number field and parameter spaces which parameterize such deformations. 
\item[{\href{https://arxiv.org/pdf/2401.13508}{Constructions III} \cite{joshi-teich-rosetta}}] Provides a formulation and a  proof of Mochizuki's Corollary 3.12 and also provides a `Rosetta Stone' to facilitate a parallel reading of my work and Mochizuki's work. Notably, it provides a canonical construction of the global version of the local correspondences of \cite{joshi-teich-estimates} using the global theory of \cite{joshi-teich-def} and provides a new approach (with proof) to \moccor\ (a discussion of this is given below).
\item[{\href{https://arxiv.org/pdf/2403.10430}{Constructions IV (This paper)}}] Deals with the proof of the \abc\ as asserted in \cite{mochizuki-iut4}.
\item[{\href{link to be added}{Worked Example} \cite{joshi-expository}}] [\textit{in preparation}] Illustrates the theory of the present series of papers by means of a worked example.
\item[Recommended Reading] {\bf(1)}
\href{https://math.arizona.edu/~kirti/report-on-scholze-stix-mochizuki-controversy.pdf}{The status of the Scholze-Stix Report and an analysis of the Mochizuki-Scholze-Stix Controversy}, {\bf(2)} \href{https://math.arizona.edu/~kirti/joshi-teich-quest.pdf}{Mochizuki’s Corollary 3.12 and my quest for its proof}, and {(\bf3)} \href{https://arxiv.org/pdf/2111.06771}{Comments on Arithmetic Teichmuller Spaces}.
\end{description}
It will be useful to understand the logical dependency of various theorems of \cite{mochizuki-iut3} and \cite{mochizuki-iut4}. This is not the  order in which these results appear in \cite{mochizuki-iut3}, \cite{mochizuki-iut4}. But in my opinion, this is the order in which they should be assimilated.
\begin{enumerate}[label={(\bf\arabic{*})}]
\item{} \cite[Corollary 2.2]{mochizuki-iut4} = \Cref{th:existence} \\
Existence of Initial Theta Data. These data are required for construction  of the set $\thetam,\thetaj$ (discussed below).
\item{} \moccor\ = \constrthr{Corollary }{cor:cor312} \\ This requires the previous result [IUT 4, Corollary 2.2] for constructing  $\thetam,\thetaj$ and provides lower bound on the volume  ${\rm Vol}(\Theta)$ for $\Theta\in\{\thetam,\thetaj\}$.
\item{} \cite[Theorem 1.10]{mochizuki-iut4} = \Cref{th:main-bound} (using the upper bound given by  \Cref{th:theta-upper-bound}). 
\item{} \cite[Corollary 2.3]{mochizuki-iut4} = \Cref{th:main-dioph-thm}\\
Provides the main result (Vojta's Inequality).
\item{} \Cref{th:main-thm} 
The $abc$ and the Szpiro Conjecture.
\end{enumerate}
\subsubsection*{Note}
For  ease of comparison between my papers and Mochizuki's, in [\cite{joshi-teich-rosetta}, \cite{joshi-teich-abc}] I have preserved Mochizuki's original content-wise appearance by placing results of \cite{mochizuki-iut3} in \cite{joshi-teich-rosetta} and results of \cite{mochizuki-iut4} in this paper. But in my opinion, the above ordering is more logical and this is the ordering in which these results  \textit{should} be read.

\subsection{A brief outline of the proof}\label{ss:outline}\nwss
This paper (in \Cref{th:main-dioph-thm}) completes  the remarkable proof of the \abc\ asserted by Mochizuki in \iut. 

The central novel idea of the proof of the \abc\ claimed by Mochizuki in \iut\ is that arithmetic of integers (generally encapsulated in the notion of a Number Field)  is not quite as rigid  as has been previously believed and one may prove the \abc\ by averaging over a suitably chosen set of deformations of a fixed number field. In \cite{joshi-teich-def}, a deformation of the arithmetic of a number field $L$  (assumed to have no real embeddings for simplicity) is captured in the notion of an \textit{arithmeticoid} \constrtwoh{Definition }{def:arithmeticoid} and I have robustly demonstrated (in \cite{joshi-teich-def}) that there is in fact a metrisable parameter space $\bsY_L$ of such deformations of arithmetic (of a fixed number field). Notably, it makes perfect sense to talk of  deformations of a fixed Number Field $L$, as is claimed in \iut.

A \textit{holomorphoid}, $\hol{}{\by}{X/L}$ of a quasi-projective variety $X/L$ \constrtwoh{Definition }{def:global-holomorphoids} (also \constrthr{Definition }{def:holomorphoids}) consists, among other data which I ignore for this discussion, of a deformation of arithmetic i.e. an arithmeticoid $\arithl_{\by}$ (\constrtwoh{Definition }{def:arithmeticoid}) given by a point $\by\in \bsY_L$.  By \constrthr{Definition }{def:holomorphoids} a holomorphoid $\hol{}{\by}{X/L}$ also includes various Berkovich analytic spaces associated to the variety. In general (for $X$ projective) two holomorphoids $\hol{}{\by}{X/L}$ and $\hol{}{\by'}{X/L}$ have isomorphic tempered fundamental groups (and isomorphic geometric tempered fundamental subgroups) \constrone{Theorem }{thm:main} but need not provide isomorphic analytic spaces. However, the schemes underlying these holomorphoids are isomorphic to $X/L$ (\constrthr{Theorem }{th:exist-cats-of-schemes-over-L}, \constrthr{Remark }{re:cat-sch-holmorphoid}). 

\textit{In particular, one can think of a holomorphoid $\hol{}{\by}{X/L}$ as representing a quasi-conformal equivalence class of $X/L$} (see \constrtwoh{\ssep}{sss:classical-analogies}).  Hence the term holomorphoid serves as a reminder that working with the category of holomorphoids provides an analog of classical Teichmuller Theory.

Now let me come to the principal reason why this theory is of interest in Diophantine Geometry. It has been well-understood in Diophantine Geometry that if one had a global Frobenius morphism for a number field, then many conjectures in Diophantine Geometry could be settled using it (see  the survey \cite[\ssep 2.4, Page 15]{mochizuki-gaussian} for an accessible discussion of this point). But, unfortunately, no global Frobenius morphism  exists for a number field. 

However, the central point of \iut\ (\cite[\ssep 2.4]{mochizuki-gaussian}) is that Mochizuki's $\flog$-Link serves as a proxy  for global Frobenius morphism on number fields. In \cite{joshi-teich-def}, I demonstrate that, in fact, there is a global Frobenius morphism (which has all the properties of Mochizuki's $\flog$-Link (and more)) on the natural parameter spaces of deformations of a fixed number fields (in particular, my work demonstrates the existence and the canonical nature of Mochizuki's $\flog$-Links)  and this can be exploited in Diophantine Geometry (as Mochizuki's work suggests).  

Specifically, the  parameter spaces $\yadl,\yadlmax, \yadl^\R$ (constructed in \cite{joshi-teich-def}) are each  equipped with a global Frobenius morphism $\bvphi:\bsY_L\to \bsY_L$ (and similarly for $\yadlmax,\yadl^\R$) and also equipped with other global and local symmetries which topologically transform a given deformation of arithmetic. [As established in \constrthr{\ssep}{ss:log-link-indentified}, $\bvphi$ corresponds to a $\flog$-Link in \cite{mochizuki-iut3}, and hence iterates of Mochizuki's $\flog$-Link correspond to global Frobenius iterates $\{\bvphi^n\}_{n\in\Z}$.]

An important consequence of \constrtwoh{Theorem }{th:consequence-inequivalent-arithmeticoids} is that arithmetic (and geometry) provided by the holomorphoids $\hol{}{\by}{X/L}$ and $\hol{}{\by'}{X/L}$ (in general), and  $\hol{}{\by}{X/L}$ and $\hol{}{\bvphi(\by)}{X/L}$ (in particular), occur  at intrinsically different valuation scales--this scaling is analogous to the  scaling behavior seen in quasi-conformal mappings between complex domains or Riemann surfaces (where this scaling behavior is encapsulated in the theory of the  Beltrami differential or the  Beltrami coefficient \cite[Section 1.4.1]{imayoshi-book}). [Mochizuki's discussion of the importance of this valuation scaling aspect  for his theory is in \cite[Remark 1.5.3, Page 461]{mochizuki-iut3}.] It is precisely the existence of such intrinsic differences in valuation scales which the averaging techniques considered  here (and in \iut) exploit. 

In the context of Szpiro's conjecture, each holomorphoid $\hol{}{\by}{X/L}$ of any elliptic curve $X/L$ provides a Tate parameter at a prime of bad semi-stable reduction and so the Tate parameter may be considered as a function (which should be elevated to a sheaf using Grothendieck's function sheaf dictionary but this is left for other writers) on the parameter space $\bsY_L$. 

Roughly speaking the idea of the proof of the $abc$-conjecture is to average over a suitable subset of holomorphoids of $X/L$. One should think of such averagings as the arithmetic analogs of the volume computations which are quite familiar in classical Teichmuller Theory--the closest analog of what one does here occurs in  \cite{mirzakhani07} (see the survey \cite{wright19}). An example of such an averaging (or integration) is discussed in the genus one case in \cite[Section 3]{wright19}, and consists of integrating functions of geodesic length (a manifestly metric dependent quantity) over an appropriate classical Teichmuller Space. Likewise, in \iut\ and \present\ one wants to average the Tate parameter (a $p$-adic metric dependent quantity) over a configuration space of a suitable Arithmetic Teichmuller Space.  [In fact, precisely because of the analogy with \cite{mirzakhani07}, one is optimistic that the methods of the \present\ should yield higher genus analogs of the genus one estimates of \cite{mochizuki-iut4}.]

In order to apply Arithmetic Teichmuller Theory developed in \present, one needs to prove the existence of elliptic curves equipped with Initial Theta Data \inithtdata\ (this corresponds to \cite[\ssep 3.1]{mochizuki-iut1}). This is a global assertion (i.e. requiring working over number fields) which is established in \Cref{th:existence}.  

[\Cref{th:existence} corresponds to \cite[Corollary 2.2(i,ii)]{mochizuki-iut4}. Note that I have reorganized \cite[Corollary 2.2]{mochizuki-iut4} into two separate portions consisting of \cite[Corollary 2.2(i,ii)]{mochizuki-iut4} and \cite[Corollary 2.2(iii)]{mochizuki-iut4} and in my approach  \cite[Corollary 2.2(iii)]{mochizuki-iut4} is subsumed in the proof of the main theorem, \Cref{th:main-dioph-thm}.]
 
The existence of Initial Theta Data (\Cref{th:existence}) allows one to proceed, via \cite{joshi-teich,joshi-untilts,joshi-teich-def,joshi-teich-estimates} to the two core constructions of \cite{joshi-teich-rosetta} which play a central role in the proof of the $abc$-conjecture given here.  

The first construction, in  \constrthr{\ssep}{ss:mochizuki-adelic-ansatz}, is the construction of a set, denoted there by $\tSigma_{L'}$, of tuples of holomorphoids $X/L$ which is called \textit{Mochizuki's Adelic Ansatz} (here $L'/L$ is a certain finite extension whose precise description, given in \constrthr{\ssep}{ss:theta-data-fixing2}, is not relevant at the moment). This set $\tSigma_{L'}$ may be understood as a graph of a certain correspondence or as a certain configuration space of points of $\yadlp^\ells$ and each point of $\tSigma_{L'}$ is a $\Theta_{gau}$-Link in the sense of \cite{mochizuki-iut3} (this is proved in \constrthr{\ssep}{sss:ansatz-theta-link}, \constrthr{\ssep}{ss:Mochizuki-theta-gau-link}). [It may be useful for readers to first understand the local version of this construction given in \cite{joshi-teich-estimates}.] Mochizuki's Adelic Ansatz is stable under action of the global Frobenius $\bvphi$ i.e. stable under $\flog$-Links and also under other symmetries acting on $\yadlp$ (\constrthr{Theorem }{th:adelic-theta-link}). [Mochizuki does not construct this set, but in my view the precise definition of $\tSigma_{L'}$ provides greater transparency and clarity to the subsequent construction.]

The second core construction   is the construction, in \constrthr{\ssep}{se:mochizuki-construction-thetam}, of an adelic set $$\thetam=\prod_{p\in\vq^{non}} \thetamp,$$
which collates the tuples of Tate parameters arising from holomorphoids given by points of  Mochizuki's Adelic Ansatz $\tSigma_{L'}$. By construction, this set is stable under the iterates of the global Frobenius morphism $\bvphi$ constructed in \constrtwoh{Definition }{def:global-Frobenius-def} and \constrtwoh{Corollary }{cor:global-frobenius}. As the global Frobenius is identified with Mochizuki's $\flog$-Link (\constrthr{\ssep}{ss:log-link-indentified}), this stability property of $\thetam$ is equivalent to the stability under the iterates of $\flog$-Links in \cite{mochizuki-iut3}.
 [Mochizuki's construction of this set is essentially the content of \cite[Theorem 3.11]{mochizuki-iut3}.]  There is, in fact, a more naturally defined set  $$\thetaj=\prod_{p\in\vq^{non}} \thetajp,$$  constructed in  \constrthr{\ssep}{se:construction-of-thetaj-and-thetam},
which is also stable under iterates of the global Frobenius morphism and which provides greater insight into the nature of $\thetam$; for instance using $\thetaj$ one obtains, in \constrthr{Theorem }{th:moccor}, an intrinsic proof of \moccor\ formulated for this set. [\textit{One may think of $\thetam$ as the Galois cohomological version of $\thetaj$.}] However,  making upper bound estimates directly with $\thetaj$ seems a bit difficult at this juncture and hence this paper, following \cite{mochizuki-iut3}, works with $\thetam$. [The collation of Tate parameters is a bit tricky but is carefully detailed in \constrthr{\ssep}{se:construction-of-thetaj-and-thetam} for $\thetaj$ and \constrthr{\ssep}{se:mochizuki-construction-thetam} for $\thetam$.]

To understand why such adelic (i.e. apparently local collection of) constructions encode global properties, one must understand the following:
\benumlab
\item the (logarithmic) height of an algebraic number $x\in \bL^*$ is a function of the choice of a (normalized) arithmeticoid $\by$ (i.e. height of $x$ depends on the chosen deformation of arithmetic) and I will write this dependence as $h_{\by}(x)$  and note that $h_{\by}(x)$ is a sum of all local contributions (except for explicating the dependence on the arithmeticoid, the definition of height given in  \constrtwoh{Definition }{def:heights} is the same as the one given in \cite[\ssep 1.5.1]{bombieri-gubler}). More precisely, one should not think of height (of non-zero algebraic numbers) as a  function $$\bL^*\to \R$$
but instead think of height  as a function
$$\bL^*\times \yadl\to \R$$
given by $(x,\by)\mapsto h_{\by}(x)$.
\item  Simplest way to understand the dependence of $h_{\by}(x)$ on $\by$, is to understand the effect that the global Frobenius $\bvphi$ has on the height i.e. understand $h_{\by}(x)$ and $h_{\bvphi(\by)}(x)$. This is discussed in \constrtwoh{\ssep}{se:heights}. Let me briefly explain this. Let $v\in\vl$ be any non-archimedean prime in the support of the fractional ideal $(x)$.  The key point of \constrtwo{Theorem }{th:flog-kummer-correspondence}{\bf(3,4)} is that the contribution of $v$ to the height function $h_{\by}(x)$  changes (in general) under application of global Frobenius  i.e. under passage  $\by\mapsto\bvphi(\by)$. [In the specific context of  elliptic curves, the precise assertion  is that the absolute values of Tate parameters rise under global Frobenius $\bvphi$ (this assertion is an immediate consequence of \constrtwo{Theorem }{th:flog-kummer-correspondence}{\bf(3,4)}). Such an  assertion is also contained in \cite[Remark 3.6.2]{mochizuki-topics3} and is the reason why  $\flog$-Link (and their iterates) are important in \cite{mochizuki-iut3}.]
\item[\bf{($\textbf{2}'$)}] There is a subtler variant of {\bf(2)}: there is an action  $L^*\act\yadl$ (\constrtwoh{Theorem }{th:galois-action-on-adelic-ff}), in which $z\in L^*$ operates on $\by\in\yadl$ through the powers of the (local) Frobenius morphism of the primes in the support of $z$ and this raises the absolute $\abs{x}_v$ of $x\in \bL^*$, but only at the primes $v$ in the support of $z$. This allows us to  compare the two height functions $h_{\by}(x)$ and $h_{z\cdot\by}(x)$ (for a worked example of this see \constrtwoh{\ssep}{re:normalization2}). Moreover for this action one has $z\cdot\by=\by$ if and only if $z$ is a root of unity contained in $L$ (\constrtwoh{Theorem }{th:fundamental-property-frobenius}). Mochizuki considers this action in \cite[\ssep 3]{mochizuki-iut3}--his discussion of this is in \cite[Remark 1.2.3(ii)]{mochizuki-iut3}.
\item Since the local non-archimedean contributions to $h_{\by}(x)$ are bound together with archimedean contributions by means of the (global) product formula (for the arithmeticoid $\arithl_{\by}$), one deduces that the global Frobenius, in fact,  changes global heights in a subtle and complicated way in order to ensure that the product formula holds  for the arithmeticoids corresponding to $\by$ and $\bvphi(\by)$. 
\item This difference between the height functions $h_{\by}$ and $h_{\bvphi(y)}$ is also a key point in the proofs of the geometric Szpiro Inequality given in \cite{bogomolov00}, \cite{zhang01}. This is detailed in  \constrthr{\ssep}{ss:geom-case-height-frobenius}. [See \cref{ss:comp-with-the-geom-case} for more on this.]
\item The product formula for an arithmeticoid  $\arithl_{\by}$ provides a hyperplane (in an infinite dimensional $\R$-vector space) $$H_\by\subset \P(\oplus_{v\in\vl} \R_v) \qquad \text{where } \R_v=\R\text{ for each }v\in \vl$$ which depends on $\by$, and hence, in fact, gives rise to a non-trivial period morphism  $$\yadl\to \P(\oplus_{v\in\vl} \R_v)$$ given by  $\by\mapsto H_\by$ (see \constrtwoh{Theorem }{th:hyperplane} and \constrthr{Theorem }{th:hyperplane-theta-lnk})  and these results provide a more natural way to understand the above mentioned properties.
\item To understand how these considerations come into play in the context of the \abc, let $(a,b,c)$ be a primitive $abc$-triple. Then one wants to estimate  $h_\by(a\cdot b\cdot c)$. By the above considerations, this height may be compared with $h_{\bvphi(\by)}(a\cdot b\cdot c)$ and also with the heights  $\left\{ h_{\bvphi^n(\by)}(a\cdot b\cdot c)\right\}_{n\in\Z_{\geq0}}$ of all global Frobenius iterates of $\by$. By the relationship between the $p$-adic norms  corresponding to $\by$ and $\bvphi(\by)$ (see \constrtwo{Theorem }{th:flog-kummer-correspondence}), if $p|(a\cdot b\cdot c)$, then one sees that iteratively applying $\bvphi$ (to $\by$) shaves off powers of $p$ from the exponent of $p$ in $\abs{a\cdot b\cdot c}_v$ for any prime $v|p$. Therefore, it appears  natural to replace $a\cdot b\cdot c$ by its conductor i.e. the conductor appears naturally in the context of the \abc\ because of the properties of the (global) Frobenius morphism! Further, one may  even bound $h_\by(a\cdot b\cdot c)$ by the supremum 
$$h_\by(a\cdot b\cdot c)\leq \sup_{(a,b,c)}\left\{ h_{\bvphi^n(\by)}(a\cdot b\cdot c): n\in\Z_{\geq0} \right\}$$
with the supremum over heights of all Frobenius iterates of all primitive $abc$-triples and try and obtain an upper bound on this supremum. \textit{Importantly,} as discussed in \constrthr{\ssep }{se:appendix-geom-case}, this is precisely the sort of strategy for bounding heights deployed using the global Frobenius morphism $\vphi_\infty$ (\constrthr{Definition }{def:global-frob-classical-case}) in the proof of the Geometric Szpiro Inequality due to \cite{bogomolov00}, \cite{zhang01}. Of course, one does not know how to directly bound the above supremum (and it may be infinite), and so, roughly speaking, the path taken in the proof of the \abc\ uses the reduction \cite{mochizuki-general-pos} to bound a similar sort of supremum obtained using  Tate parameters of holomorphoids of Frey elliptic curves, corresponding to  $abc$-triples, whose $j$-invariants lie in compactly bounded subsets.
\eenum
The sets $\thetam$  (resp. $\thetaj$), constructed in \constrthr{\ssep}{se:mochizuki-construction-thetam} (resp. \constrthr{\ssep}{se:construction-of-thetaj-and-thetam}),  capture this natural variation of heights by virtue of stability of $\thetam$ under iterates  $\by\mapsto \bvphi^n(\by)$ of the global Frobenius morphism. [In \cite{mochizuki-iut3}, this stability property corresponds to stability under iterates of $\flog$-Link.] In particular, the absolute values of the Tate parameters grow in the set $\thetam$ (\constrtwo{Theorem }{th:flog-kummer-correspondence}). This is the (global) Diophantine significance of these sets. However, the structure of $\thetam$ (and  $\thetaj$) is quite complicated and at the moment our understanding of these sets is far from complete. 

Mochizuki's central assertion in \cite{mochizuki-iut3,mochizuki-iut4} is that upper and lower bounds on $\thetam$ leads to a proof of the \abc\ via a proof of Vojta's Inequality (the equivalance of the $abc$-Conjecture and Vojta's Inequality was established in \cite{frankenhuysen2002}). \textit{However,} the central problem which has engulfed \iut\ after the appearance of \cite{scholze-stix,scholze-review} is whether or not the set $\thetam$ and the structures required to construct it exist at all (see \cite{joshi-report}, \cite{joshi-teich-quest}, \cite[Introduction]{joshi-untilts} for a thorough discussion of these issues). Thanks to the \present\ one now has a robust demonstration of the existence of the structures required for the construction of $\thetam$.

The proof of the \abc\ described below is based on the one given in \cite{mochizuki-iut4}. However, I have organized the proof from my point of view for greater transparency,  supplying additional proofs which I felt were  necessary for completeness.

Before proceeding further let me introduce some more notation. Let $U=(\P^1-\{0,1,\infty\})(\bL)$ and $U_d=(\P^1-\{0,1,\infty\})(\bL)^{\leq d}\subseteq U$ be the subset of $U$ consisting of $\bL$-points of degree at most $d$ and view $\P^1_L-\{0,1,\infty\}$ as the $j$-line for the Legendre family of elliptic curves $X_\lambda:y^2=x(x-1)(x-\lambda)$. Write $C_\lambda$ for the projective elliptic curve corresponding to  $X_\lambda$ and for simplicity of notation, write $C=C_\lambda$.

In order to apply the above theory to the proof of the \abc, one uses Mochizuki's reduction \cite{mochizuki-general-pos} which reduces the proof of the $abc$-conjecture to establishing Vojta's Height Inequality (see \cref{se:vojta-ineq} for the statement) for compactly bounded subsets of $U_d$  (see \cref{ss:compact-bnded} for the definition).

The quantity  interest here is the Tate divisor (see \cref{ss:tate-divisor}): for an elliptic curve $C$ over $M\supset L$  (for certain finite extensions $M/L$), let
$$
\fT(C/M)= \frac{1}{[M:\Q]}\cdot\sum_{w\in \vnon_L} \ord_w(q_w)\cdot \log w,
$$
where $q_w$ is a Tate parameter at $w\in\vmnon$,
be the normalized degree of the Tate divisor of $C/M$, considered as a real valued function $$\fT:\sM_{ell}(\bL)\to \R.$$ \textit{Strictly speaking,} an important point of the theory is that $\fT(C/M)$ depends on the choice of a holomorphoid $\hol{}{\by}{C/L}$  and one should write $\fT(C/M)_{\by}$ to indicate this choice and $\fT(C/M)$ (as above) in fact denotes the $\fT(C/M)_{\by_0}$ for a holomorphoid $\hol{}{\by_{0}}{C/L}$ with a specific choice  of the standard arithmeticoid $\by_0\in\yadl$ which is declared in \cref{ss:basic-setup}. 

[One expects that the averaging which is carried out  in the present paper to be replaced by the evaluation of an integral $$\int_{\Sigma}\fT(C/M)_{\by}d\boldsymbol{\mu},$$  for a suitable subset $\Sigma\subset \yadl\times\sM_{ell}(\bL)$ equipped with a suitable measure $d\boldsymbol{\mu}$. But presently one is far from being able to establish this expectation. Realizing such  integrals would bring one closer to the analogy with the volume integrals of classical Teichmuller Theory which occur in \cite{mirzakhani07} (\cite{wright19}). A similar expectation has been  tacitly voiced by Mochizuki in \cite[Table on Page 115]{mochizuki-iut1} where he suggests that his log-shell is the canonical unit of volume analogous to the canonical Kahler volume form.]

At any rate, on a given compactly bounded subset satisfying a natural $2$-adic restriction, $\fT(C/M)$ be replaced, up to a bounded error,  by a variant, denoted by $\fT_2(C/M)$, supported on primes away from the prime $2$. The divisors $\fT(C/M),\fT_2(C/M)$ are not  stable under passage to field extensions (as primes of additive reduction may become primes of semi-stable reduction), so one replaces $\fT(C/M)_2$ by the normalized arithmetic degree of another divisor $\fq=\fq(C/M)$  (\cref{ss:tate-divisor}), which is  similar to the divisor $\fT$, but which  receives contributions only from primes of $M$ lying over a prime (not dividing $p=2$) of semi-stable reduction in the field of definition $\lmod$ of the moduli point associated to $C$. By \Cref{pr:tate-divisor}, this new divisor $\fq$ is stable under base-field extensions. Let me remark  that $\fq$ also appears in the fundamental estimate  \constrthr{Corollary }{cor:cor312} (\moccor).

The key bound for $\fq$ is established in \Cref{th:main-bound} (which corresponds to \cite[Theorem 1.10]{mochizuki-iut4}). A crucial step in the proof of \Cref{th:main-bound} is \Cref{th:theta-upper-bound}.

Let me discuss the proof of \Cref{th:theta-upper-bound}. This step requires  lower and upper  bounds for $\thetam$ and so this is one of the central points in the proof of the \abc\ where the theory of Arithmetic Teichmuller Spaces developed in the \present\ come into play. The lower bound in \Cref{th:theta-upper-bound} is given by \constrthr{Corollary }{cor:cor312}. So let me now discuss the upper bound estimate.

On a fixed compactly bounded set $\sZ\subset U_d$, and with $C=C_\lambda\in \sZ$, invoking the Existence of Initial Theta-Data (given by \Cref{th:existence}) also allows one to talk about $\thetam$ (constructed in \constrthr{\ssep}{se:mochizuki-construction-thetam}) and one can proceed to the problem of bounding a suitably defined volume, $\mathit{Vol}(\thetam)$,  of $\thetam$. Under \Cref{th:existence}, one also has the conclusion of \constrthr{Corollary }{cor:cor312} which provides a lower bound on $\mathit{Vol}(\thetam)$.  To avoid dealing with Mochizuki's notational and sign conventions so as to keep the notation uncluttered for clarity in this discussion, I will simply pretend that $\log\mathit{Vol}(\thetam)>0$. [In \cite{mochizuki-iut3,mochizuki-iut4}, Mochizuki works with the specific notational convention of writing  $-\abs{\log\mathit{Vol}(\thetam)}$ (\textit{strictly speaking}, I should be writing $-\frac{1}{\ells}\abs{\log\mathit{Vol}(\thetam)}$ in place of $-\abs{\log\mathit{Vol}(\thetam)}$ but this is not important for the present discussion) which makes keeping track of correct signs a tedious task  but this has been carried out  here in \cref{se:first-main-bnd} and also in \cite{joshi-teich-rosetta}.]

Let me comment on an important and subtle point at this juncture. It is important to remember that  by the virtue of stability under the global Frobenius morphism, one has the freedom of the choice of an arithmeticoid in the computing the lower bound on the log-volumes of $\thetam$: the lower bound on $\thetam$ of \constrthr{Theorem~}{th:moccor-Mochizuki-form2} is independent of the choice of the arithmeticoid one makes at this point. Specifically in \Cref{th:theta-upper-bound}, the lower bound will be computed using the arithmeticoid $\by_0'=\bvphi(\by_{0})$ while the upper bound  will be computed using $\by_{0}$ (the upper bound is not independent of this choice).  Each arithmeticoid $\by_0'=\bvphi(\by_{0}),\by_0$ provides an isomorphic copy of our number field and hence one has chosen the same global normalization for both the copies, the it does not matter which arithmeticoid is used for global computations. Hence one seeks a bound of the form:
\begin{equation}\label{eq:0}
0<A \leq  \log\mathit{Vol}(\thetam^{\bvphi(\by_0)})= \log\mathit{Vol}(\thetam^{\by_0}) \leq A'.
\end{equation}

In \cite[Remark 3.3.3]{mochizuki-iut3}, Mochizuki adopts a similar strategy (calls it $\flog$-shifting i.e. Frobenius-shifting), available in his case by the virtue of the stability under iterates of $\flog$-Links in his constructions.[\textit{The previous release of this paper had omitted this Frobenius-shifting.}] As Mochizuki remarks in \cite[Remark 4.8.2(iii)]{mochizuki-iut2}, this strategy (of using $\by_0'=\bvphi(\by_0)$ and $\by_0$) is best understood in terms of vector bundles equipped with a filtration and is tantamount to comparing slopes of $(V,Fil(V))$ and $(V,\bvphi(Fil(V)))$ (see \cite{joshi16} where slope inequalities are established in a similar bundle situation on curves in presence of Frobenius; the inequalities of \cite{joshi16} are strongly reminiscent of the geometric proofs of the Szpiro inequality in \cite{kim97}, \cite{beauville02}). [For additional details on this important point see \Cref{re:bundle-rmk}.]

Since each of the sets $\thetam^{\by_0'}, \thetam^{\by_0}$ is an adelic set by construction and most primes contribute zero to the log-volume, one can compute $\log\mathit{Vol}(\thetam^{\by_0'})$ by taking the sum  $\log\mathit{Vol}(\thetamp^{\by_0'})$ of the local contributions at all primes (and similarly for the respective quantities for $\by_0$). [Mochizuki's  upper bound in \Cref{th:theta-upper-bound} is obtained by summing over all the local contributions.]

Locally, however, at each prime $p$, $\by_0'$ presents quite distinct arithmetic and geometry from that of $\by_0$ (because of the relationship $\by_0'=\bvphi(\by_0)$). Hence one cannot make local term by term comparisons at each prime as one may naively hope to do without violating concurrent global choice of normalizations for the arithmeticoids $\arithl_{\by_0'}=\arithl_{\bvphi(\by_0)}$ and $\arithl_{\by_0}$ given by $\by_{0}',\by_{0}$.
At any rate, one therefore has for $\thetam^{\bvphi(\by_0)}$:
\begin{equation}\label{eq:1}
0<A = \sum_p \log\mathit{Vol}_p(\thetamp^{\bvphi(\by_0)}) = \log\mathit{Vol}(\thetam^{\bvphi(\by_0)})
\end{equation}
and for $\thetam^{\by_0}$:
\begin{equation}\label{eq:2}
\log\mathit{Vol}(\thetam^{\by_0})=\sum_p \log\mathit{Vol}_p(\thetamp^{\by_0}) \leq A'.
\end{equation}
Now the calculation proceeds as follows. For each prime $p$, one independently picks  real numbers $z_p,z_p'\in\R_{\geq0}$, such that 
\be\label{eq:lower} A=\sum_p z_p\leq \sum_p \log\mathit{Vol}_p(\thetamp^{\by_{0}'})\ee
and \be\label{eq:upper} \sum_p \log\mathit{Vol}_p(\thetamp^{\by_0}) \leq \sum_p z'_p = A'\ee
with $z_p,z_p'=0$ for all but a finite number of primes while $z_p>0$ for a finite, non-empty set of primes. 
\textcolor{red}{Note} that the important requirement on $z_p,z_p'$ is that \cref{eq:lower} and \cref{eq:upper} hold. However as indicated earlier, $z_p$ and $z_p'$ reference $\by_0'$ and $\by_0$ respectively and one cannot claim that $z_p\leq z_p'$ holds for all primes $p$ without violating the middle equality of \eqref{eq:0} because the local data for $\by_0'$ and $\by_0$ differ by Frobenius at each prime number $p$.
Then one has 
\be\label{eq:vol-ineq}
A=\sum_{p}z_p \leq \log\mathit{Vol}(\thetam^{\by_0'}) = \log\mathit{Vol}(\thetam^{\by_0}) \leq \sum_{p}z_p'=A'
\ee

Now the proof of Vojta's Inequality for $\sZ$ is obtained as follows.  Calculating $A'$ explicitly, by finding $z'_p$ for each $p$,  Mochizuki  observes, in \cite[Theorem 1.10]{mochizuki-iut4}, that 
\benumlab
\item  $A'=C\cdot A$. This is a global equation, which only appears after one sums over all primes, and hence \Cref{eq:vol-ineq} may be written as  (for the precise assertion see \Cref{th:theta-upper-bound})
$$
0<A \leq \log\mathit{Vol}(\thetam)\leq A'= C\cdot A
$$
and hence $C\geq 1$, 
\item  and secondly that $C\geq 1$ implies Vojta's Inequality for $\sZ$ (for the precise assertion see \Cref{th:main-bound}).
\eenum

Finally, the main theorem  i.e. the Vojta Inequality is  established here in \Cref{th:main-dioph-thm} (this corresponds to \cite[Theorem A/Corollary 2.3]{mochizuki-iut4}). This also proves \Cref{con:abc-vojta} and hence proves the \abc\ (\Cref{con:abc}).

\subsection{What is new in my approach?}
Let me briefly explain  the new ideas in the \present\ and how they relate to \iut. A detailed discussion of these points, with references to precise results in both my papers and Mochizuki's papers, is given in the `Rosetta Stone' \constrthr{\ssep }{se:intro-rosetta-stone}.
\begin{description}
\item[1: Arithmetic Holomorphic Structures I] My work provides a new, natural and robust definition of the notion of arithmetic holomorphic structures in \cite{joshi-untilts}, \cite{joshi-teich}; one obtains Mochizuki's definition of arithmetic holomorphic structures by applying the tempered fundamental group functor (and variants) to the arithmetic and geometric objects (called holomorphoids) of my theory. The objects of my theory include an algebraically closed perfectoid field and its tilting data for each prime and naturally provides \textit{Mochizuki's \'etale pictures} (i.e. group theoretic data) plus \textit{Mochizuki's Frobenius Pictures.} (i.e. realified Frobenioid data). [A tabular summary of the key objects of the two theories is given in \cref{tab:comptable}.] This approach provides the cleanest geometric proof that arithmetic holomorphic structures and the sort of Teichmuller Theory asserted by Mochizuki exists--not just in the genus one case as asserted by Mochizuki, but in all genera, and even in higher dimensions. [That such a theory cannot exist is the core objection \cref{the scholze-stix-report}{\bf(1)} of \cite{scholze-stix}.]
\item[2: Mochizuki's $\Theta_{gau}$-Links and $\flog$-Links--local version] In \cite{joshi-teich-estimates}, I provide a natural (but local) construction of Mochizuki's $\Theta$-Link as a geometric correspondence on a suitable Fargues-Fontaine curve and also a natural demonstration of the non-triviality of its key valuation scaling property (asserted, but not proved, by Mochizuki). [For a worked example of the $\Theta_{gau}$-Link see \cite{joshi-expository}.] The second important discovery of \cite{joshi-teich-estimates} is that the Frobenius morphism of suitable Fargues-Fontaine curves has all the properties of Mochizuki's (local) $\flog$-Link. This not only provides a natural reason why the $\flog$-Link is important in \iut, but also provides a geometric and hence a quantitatively more useful version of Mochizuki's (local) $\flog$-Link. For a worked example of the $\flog$-Link see \cite{joshi-expository}. [That such objects cannot be non-trivial are the core objections \cref{the scholze-stix-report}{\bf(2, 3)} of \cite{scholze-stix}.]
\item[3: Mochizuki's $\Theta_{gau}$-Links and $\flog$-Links--global version] Mochizuki's (global) $\flog$-Link is a proxy for the Frobenius morphism and its existence is construction of the global Frobenius morphism of \cite{joshi-teich-def}--the $\flog$-Link aspect of this is detailed in \cite{joshi-teich-rosetta}  (but also see the next point {\bf(4)}). The constructions of the global versions of $\Theta_{gau}$ is given in \cite{joshi-teich-rosetta}. [The existence of these global objects  contradicts the core assertions \cref{the scholze-stix-report}{\bf(2, 3)} of \cite{scholze-stix}{\bf(2, 3)}.] 
\item[4: Arithmetic Holomorphic Structures II] An important and completely new part of the global theory is established in \cite{joshi-teich-def}--this paper is devoted to the demonstration of the existence of a Teichmuller Theory of a fixed Number Field (this allows one to talk about Arithmetic Holomorphic Structures on a fixed Number Field). The existence of such a theory is asserted by Mochizuki  \cite[Page 25]{mochizuki-iut1} and its existence underpins \iut. My work provides the most intrinsic formulation of such a theory. [\cite{scholze-stix} is silent on this but tacitly imply that arithmetic of a number field is rigid i.e. no deformations can exist.] Important points to note about \cite{joshi-teich-def} are:
\begin{enumerate}[label={\bf(\alph{*})}]
	\item  this theory (\cite{joshi-teich-def}) includes archimedean contributions 
	\item a global Frobenius morphism (which appeared in points {\bf(1), (2)} above). [This morphism, together with the product formula changes global heights.]
	\item The product formula as a period morphism: the construction of a global period morphism, which assigns to a given arithmetic holomorphic structure (i.e. an arithmeticoid), the hyperplane provided by the logarithm of the (global) product formula. [Mochizuki asserted that the product formula varies in his theory (this point is now demonstrated robustly by the construction of this period morphism). [The third core objection (\cref{the scholze-stix-report})  of \cite{scholze-stix} is that this phenomenon cannot exist.]
	\item I provide an intrinsic construction of Mochizuki's adelic log-shell in \cite{joshi-teich-rosetta} which explains why the log-shell is a natural object.
	\item In \cite{joshi-teich-rosetta}, an intrinsic construction of Mochizuki's Theta-value locus $\thetam$ is given, as well as a more natural version of $\thetam$, denoted $\thetaj$, is also constructed. As shown in \cite{joshi-teich-rosetta}, Mochizuki's Corollary 3.12 is, essentially, a tautology for $\thetaj$. This demystifies the said corollary and it is the reason why I have asserted in \cite{joshi-teich-quest} that this corollary is evidently a natural assertion.
	\item Mochizuki's method of establishing his lower bound (\moccor) for $\thetam$ relies on \topics\ while I provide a natural approach based on the fact that the set $\thetam$ is adelic and so amenable to local computations (along the lines of \cite{joshi-teich-estimates}). This is analogous to the fact that while heights of algebraic numbers are global, they can certainly be computed by summing local terms.
	\item From a philosophical and mathematical point of view, the familiar invariants such as  the discriminant, the different and the Swan conductor are  sensitive to the additive structure and this is how anabelomorphy comes into play. This was the key discovery of \cite{joshi-anabelomorphy}. In the present paper, I have made this view (and mathematical fact) completely manifest in my approach to \cite{mochizuki-iut4}. Thus in some sense, my work makes quite natural Mochizuki's idea of proving Szpiro's conjecture (about discriminants) by considering variations of rings structures.
	\item One of the important points in my approach to \iut\ is that I have made the role of $p$-adic Hodge Theory completely transparent. While Mochizuki makes it clear that collections of Hodge-Theaters should be viewed as a proxy for variations of Hodge Structures, a lack of quantification (in \iut) of this point has been a major stumbling block in understanding the claims of \iut.
\end{enumerate}
\item[5: Canonical geometric description of Mochizuki's Indeterminacies] My work on quantifying Arithmetic Holomorphic Structures (see {\bf(1,4)}) provides  a canonical geometric description of Mochizuki's Three Fundamental Indeterminacies Ind1, Ind2, Ind3 \cite[Theorem 3.11]{mochizuki-iut3}. This is detailed in \constrthr{\ssep}{ss:Mochizuki-indeterminacies}. Mochizuki's Indeterminacy Ind3 is subtlest of the three and Mochizuki's quantification of it, in \cite{mochizuki-iut3}, is woefully inadequate. On the other hand it is central to the formulation of the Mochizuki's Theory. In \constrthr{\ssep}{ss:Mochizuki-indeterminacies} provides a natural demonstration of the existence of Ind3 and of its geometric origin and its relationship to arithmetic.

\item[6: Construction of $\thetam$ and $\thetaj$] My approach to constructing Mochizuki's Theta-values locus $\thetam$ is slightly different from that of \cite{mochizuki-iut3}. As a precursor to constructing $\thetam$, I provide the construction of the theta-values locus $\thetaj$ which is more natural and makes transparent the fact that absolute values of Tate parameters arising from distinct arithmetic holomorphic structures can, in fact, be compared in the fixed adelic period ring $\mathbb{B}_{L'}$ (here $L'/L$ is a finite extension given by the choice of Initial $\Theta$-data). This ring is essentially the product of certain Fargues-Fontaine rings for all primes $v\in\vlp$. Moreover, this is independent of the choice of arithmetic holomorphic structures. Mochizuki's $\thetam$, on the other hand, lives in an adelic vector space constructed using $L'$ (and necessarily involves fixing one arithmetic holomorphic structure with respect to which all others are compared). Since this is only a vector space, no natural multiplicative operations exists, forcing Mochizuki to work with tensor products. From my point of view, because $\mathbb{B}_{L'}$  is a ring one can multiply (in it)  any finite collection of theta-values even if these values arise from distinct arithmetic holomorphic structures; and multiplication, by virtue of bilinearity (or multilinearity), provides natural passage to tensor products. Thus this provides a natural reason why tensor products of vector spaces arise in \cite[Section 3]{mochizuki-iut3}. All in all, $\thetaj$ is quite natural and the expected relationship between the $\thetam$ and $\thetaj$ is that $\thetam$ is the image of $\thetaj$ under a suitable connecting homomorphism in Galois cohomology sequence for Fargues-Fontaine rings but such a relationship is not presently established.
\item[7: A worked example] [\textit{in preparation}] The example of $y^2=x(x-1)(x-p)$  discussed in \cite{joshi-expository}   will help clarify the reasons for the existence of Arithmetic Teichmuller Spaces asserted by the \present.
\end{description}

\newpage
\nwss
\begin{boxedcontent}[label={tab:comptable},grow to left by=2cm, grow to right by=1cm]{Comparison of Primary Objects and their purpose in the two theories}{}
\centering
	\renewcommand{\arraystretch}{1.5}
	{\small 
		\begin{tabular}{|p{2in}|p{2in}|p{2in}|}
			\hline 
Object$\downarrow$ \ \ Purpose$\to$  & \iut\ & \cite{joshi-teich,joshi-teich-estimates,joshi-teich-def,joshi-teich-rosetta} \\ 
			\hline 
\'etale Picture (of a variety) at each prime of a number field & Provides an instance of group theoretic  data (i.e. the data of tempered fundamental group, absolute Galois groups etc) & \multirow{2}{2in}{Holomorphoid $$\hol{}{}{X/L}$$ of a variety over a number field. Each holomorphoid canonically provides both an \'etale-picture and a Frobenius-picture for every prime  of $L$.}\\ 
			& & \\
			\cline{1-2}
Frobenius Picture (of a variety) & Provides instance of Realified Frobenius data for computing local arithmetic degrees (global arith. degree is the sum of local degrees)&  \\ 
			\hline 
			Global object & A member of the isom. class of the (realified) Frobenioid of a Number Field & An arithmeticoid $$\arithl$$ of the Number field; each holomorphoid provides a canonical arithmeticoid and each arithmeticoid provides a canonical global Frobenioid.\\
			\hline
			local+global object & A Hodge Theater given by an (\'etale, Frobenius)-Picture for each prime + Global Frobenioid of the Number Field & Each holomorphoid has both local (for all primes) and global components. Each holomorphoid naturally gives rise to the associated  Hodge Theater (in Mochizuki's sense)\\
			\hline
		\end{tabular}
	} 
	\tcblower
	{\small 
	\textcolor{red}{\bf Note:}
	\begin{enumerate}[label={\bf(\arabic{*})}]
	\item A holomorphoid (resp. an arithmeticoid) is a quantitatively more precise object and working with it provides mathematical advantage over working with Mochizuki's (\'etale, Frobenius)-pictures and Frobenioids.
\item  \cite{joshi-teich-def} provides a metrisable topological space of arithmeticoids and a space (i.e. a category with a list of symmetries) of holomorphoids. In particular, in \constrthr{\ssep}{se:frobenioids},  one sees that Mochizuki's Hodge Theaters are parameterized by the Arithmetic Teichmuller Space of \cite{joshi-untilts}.
\end{enumerate}
}
\end{boxedcontent}
\newpage

\subsection{Where is Anabelian Geometry in all this?}
Let me clarify this important point. Let $L$ be a number field, let $v\in\vl$ be a prime of $L$, let $L_v$ be the completion of $L$ at $v$, and let $X/L$ be a smooth quasi-projective variety   over  $L$ (in Mochizuki's context $X/L$ is a smooth, hyperbolic curve over a number field $L$).  

One of my main observations is that for each non-archimedean $v$, the structure of the tempered fundamental groupoid of $X/L_v$, which underlies Mochizuki's approach,  is quite similar to that of  the Diamond $X^\Diamond/L_v$ associated to $X/L_v$ in \cite{scholze-diamonds}. This was pointed out explicitly in  \cite{joshi-blogpost},  \construntilts{Remark }{re:diamonds} (or \constrone{Proposition }{re:diamonds}) and more detailed version appears in \constrthr{Remark }{rem:diamonds-elaborated}. [The   definition of Arithmetic Teichmuller Space in \constrone{Definition }{def:arith-hol-space-local} is a bit more general than that of the Diamond associated to $X/L_v$ as it  glues together many Diamonds.] 

In the context of the discussion of holomorphoids considered in \cref{ss:outline}, in \cite{joshi-teich}, \cite{joshi-teich-rosetta} it is shown that each  holomorphoids $\hol{}{\by}{X/L}$ provides, for each non-archimedean $v$,  a group and the inclusion of a normal subgroup:
$$ \text{temp. fund. group of } X^{an}/L_v \supset \text{geom. temp. fund. group of }  X^{an}/L_v $$ 
and hence  a group and the inclusion of a normal subgroup:
$$ \text{\'etale fund. group of } X/L_v \supset \text{geom. \'etale fund. group of }  X/L_v$$
and for each archimedean $v$  a group and the inclusion of a normal subgroup
$$ \text{\'etale fund. group of } X/L_v\supset \text{geom. \'etale fund. group of }  X/L_v.$$ The construction of this group-theoretic data, is up to isomorphism, independent of the choice of the holomorphoid $\hol{}{\by}{X/L}$ i.e. if one has two holomorphoids  $\hol{}{\by}{X/L}$, $\hol{}{\by'}{X/L}$, then there exists a topological isomorphism between the corresponding group-theoretic data provided by these two holomorphoids. This is analogous to the fact that the Teichmuller space parameterizes Riemann surfaces with isomorphic fundamental groups. In \iut, Mochizuki works with the inclusion (for each $v\in\vl$): $$ \text{temp. fund. group of } X/L_v\supset \text{geom. temp. fund. group of }  X/L_v.$$

So anabelian geometry is built into the very foundation of my theory and it strongly resembles classical Teichmuller Theory and that is why I have asserted that my conclusions regarding the relationship between \iutthr, the \present\ and the theory of diamonds $X^\Diamond/L_v$ for each prime $v\in\vl$, are quite robust (see \constrone{Proposition }{re:diamonds}, \constrthr{Remark }{rem:diamonds-elaborated}).

My second central observation is that for each $v\in\vlnon$, the existence of 
\benumlab
\item many topologically distinct classes  of algebraically closed, complete rank-one valued fields, with residue fields $\bar{\mathbb{F}}_{p_v}$, and all containing an isometrically embedded  $\Q_{p_v}$, and
\item the existence of many inequivalent tilting data (as in \cite{scholze12-perfectoid-ihes}) for each field provided by {\bf(1)}
\eenum
given to us, with increasing precision, by the fundamental theorems of \cite{schmidt33}, \cite{kaplansky42}, \cite{matignon84},  \cite{poonen93}, \cite{fargues-fontaine},  \cite{kedlaya18},  lies at the foundation of the variation of Arithmetic Holomorphic Structures described in the \present\ and in \iut. In other words, the theory is birthed by  the vast richness of arithmetic itself and not by the group theory surrounding the Galois and the tempered (or \'etale) fundamental groups. These groups appear as symmetries of the objects of the theory, and to be sure, these groups play many important roles in the theory, but the genesis of the theory  lies  in the infinite  richness of Arithmetic itself.

\subsection{Comparison of the arithmetic and the geometric cases}\label{ss:comp-with-the-geom-case}
\newcommand{\syi}{\sY_{\infty}}
\newcommand{\sxi}{\sX_\infty}
\newcommand{\slr}{{\rm SL}_2(\R)}
\newcommand{\tslr}{\widetilde{{\rm SL}_2(\R)}}
\newcommand{\fH}{\mathfrak{H}}
There is a rich history of Diophantine conjectures being first established in the geometric case and the arithmetic case being subsequently established by imitating the strategy used in the geometric case.  
The geometric cases of the Szpiro Conjecture and the \abc\ have been known for quite some time.   During my sabbatical visit to RIMS, Kyoto, Mochizuki lectured to me on the geometric case indicating the close parallel between the geometric proofs of \cite{bogomolov00}, \cite{zhang01} and his own proof given in \iut. 
The importance of the \cite{bogomolov00}, \cite{zhang01} proofs is that these proofs are Teichmuller Theoretic in nature i.e. the proofs exploit the existence of Teichmuller Spaces and their rich geometry.
My own reflections on the geometric case are detailed in \constrthr{\ssep}{se:appendix-geom-case}. I have supplemented \cite{mochizuki-bogomolov} by demonstrating  the existence of $\Theta_{gau}$-Links and $\flog$-Links (aka global Frobenius) in the geometric case and by establishing the geometric case of Mochizuki's Corollary 3.12. 

Here I provide a comparison table (\Cref{tab:geomcomptable}) for these ideas and how they are related in the two cases.

For the geometric case, recall some notations from \constrthr{\ssep}{se:appendix-geom-case}. Let  $C/\C$ be a  smooth, projective curve of genus $g\geq 1$ over $\C$; let  $X\to C$ be a proper, generically smooth morphism whose geometric generic fiber is a smooth proper curve of genus one. The proof of \cite{bogomolov00}, \cite{zhang01} provides a Teichmuller Theoretic proof of the geometric Szpiro Inequality for $X\to C$. Let $\fH\subset \C$ be the upper half-plane.  Let  $\tslr$ be the universal cover of $\slr$. By \cite{bogomolov00}, \cite{zhang01} one has a central extension 
\be\label{eq:slr-cover}
1\to \langle z^2\rangle \to\tslr \to \slr \to 1,
\ee
where $\langle z^2\rangle\subset \tslr$ is the cyclic subgroup generated by a specific element defined in \cite{bogomolov00}, \cite{zhang01}. One also has the standard mapping $\slr\to \fH$ and therefore the composite mapping $$\tslr\to\slr\to\fH.$$ The upper half-plane $\fH$ is the classical Teichmuller Space in genus one.

Following \constrthr{Definition }{def:global-frob-classical-case}, write  and 
\begin{align*}
\syi&=\tslr\\
\vphi_\infty&=z^2\\
\sxi&=\slr=\syi/\vphi_\infty^\Z
\end{align*}
(more precisely, $\vphi_\infty(g)=z^\cdot g$ for all $g\in \tslr$).
Then $\syi,\sxi$ are the (global) Fargues-Fontaine ``curves'' in the geometric case.  Moreover, $\vphi_\infty$ plays the role of the (global) Frobenius (in the geometric case).  The $p$-adic analogs are the Fargues-Fontaine curves $\syQp$ and $\sxqp$ (\cite{fargues-fontaine}). The global arithmetic analogs are the curves are $\yadl$ and $\xadl$ constructed in \cite{joshi-teich-def}.  \Cref{eq:slr-cover} asserts that $\sxi=\syi/\vphi_\infty^\Z$ which is analogous to $\xadl=\yadl/\bvphi^\Z$ in the arithmetic case and  $\sxqp=\syQp/\vphi^\Z$ in the $p$-adic case.

The function $h:\syi\to\R$, denoted by $\tilde g \mapsto \ell(\tilde{g})$ in \cite{zhang01}, is the analog of the logarithmic height function. \Cref{tab:geomcomptable} provides a comparison of the geometric case including   Mochizuki's approach and my  approach in the arithmetic case.
\begin{boxedcontent}[label={tab:geomcomptable},grow to left by=2cm, grow to right by=1cm]{Comparison of the Geometric and the Arithmetic case}{}
\centering
	\renewcommand{\arraystretch}{1.5}
\begin{tabular}{|p{2in}|p{2in}|p{2in}|}
			\hline 
Geometric Case (genus one)    & Arithmetic case (Joshi) & Arithmetic Case (Mochizuki) \\ 
			\hline 
$\syi=\tslr$ a parameter space for hol. strs	& $\yadl$ parameter space for arith. hol. strs. &  (Hodge-Theaters and $\flog$-Links between them) \\
\hline 
global Frobenius $\vphi_\infty:\syi\to\syi$, $\vphi_\infty=\text{mult. by }z^2$ on $\tslr$ & global Frobenius $\bvphi:\yadl\to\yadl$ & proxy for global Frobenius = $\flog$-Link between two Hodge Theaters\\ 
			\hline 
			$\sxi=\syi/\vphi_\infty^\Z$ & $\yadl/\bvphi^\Z$ & ``vertically coric objects'' \\
			\hline
			hol. str. $\tilde g\in\syi$ & arith. hol.  $\by\in\yadl$ & arith. hol. str. \\
			\hline
			$\flog$-Link = Frobenius $\bvphi_\infty$  (Joshi)
			$$\tilde{g} \mapsto \vphi_\infty(\tilde{g})=z^2\cdot\tilde{g}$$ & $$\by\mapsto \bvphi(\by)$$ & $\flog$-Link = proxy for Frobenius (Mochizuki) \\
			\hline
			$\Theta_{gau}$-Link (Joshi) & $\Theta_{gau}$-Link (Joshi) & $\Theta_{gau}$-Link (Mochizuki) \\
			\hline
			height $h: \syi \to \R$ as a function of hol. str. & height $\by\mapsto h_{\by}$ as a function of arith. hol. str. & does not directly define this\\
			\hline
			$\vphi_\infty$ modifies height function and the proof exploits this property of $\vphi_\infty$ & $\bvphi$ modifies height function  and the proof exploits this property of $\bvphi$ &  $\flog$-Link modifies height function and the proof exploits this property of $\flog$-Links\\
			\hline
			Specifically $\vphi_\infty$ shaves off $$\text{multiples of }\pi$$  from the logarithmic height function  & $\bvphi$ shaves off $$\sum_p \ast_p \cdot\log p$$ (for suitable $\ast_p\in\R$) from the logarithmic height function & $\flog$-Link shaves off $$\sum_p \ast_p \cdot \log p$$ (for suitable $\ast_p\in\R$) from the logarithmic height function   \\
			\hline
			global constraint equation
			$$\prod_{i=1}^g[a_i,b_i]\cdot \prod_j c_j =1$$ & global constraint equation: 
			$$\prod_{v\in\vl} \abs{x}_v=1$$ &global constraint equation: 
			$$\prod_{v\in\vl} \abs{x}_v=1$$\\
			\hline
		\end{tabular}
	\tcblower
	My discussion of the geometric case appears in \constrthr{\ssep}{se:appendix-geom-case} (and for readers convenience also in \constrtwoh{\ssep}{se:appendix}). Mochizuki's discussion of the geometric case (and comparison with the arithmetic case) appears in \cite{mochizuki-bogomolov}.
\end{boxedcontent}
\newpage

\subsection{The Scholze-Stix Report}\label{the scholze-stix-report}
I have provided a detailed mathematical analysis of the Mochizuki-Scholze-Stix Controversy (\cite{scholze-stix,scholze-review}, \cite{mochizuki-scholze-stix-report}) in my reports \cite{joshi-final},  \cite{joshi-report,joshi-teich-quest,joshi-teich-summary-comments}.

For reasons of completeness and for the reader's reference, I recall the  core objections of \cite{scholze-stix}:
\benumlab
\item The existence of (distinct) Arithmetic Holomorphic Structures is not demonstrated in \iut\  (this is the disambiguation problem discussed in \cref{ss:why-this-is-needed}) and cannot possibly exist because of \cite[Theorem 1.9]{mochizuki-topics3}.
\item The existence and the non-triviality of the valuation scaling property of the $\Theta_{gau}$-Link is not adequately demonstrated in \iut--but it cannot be non-trivial in any case (\cite[Section 2.1.4, Section 2.1.9, Section 2.2]{scholze-stix}). 
\item In \cite[Section 2.1.3 especially Footnote 8]{scholze-stix}  it is argued that the $\flog$-Link is essentially irrelevant in \iutthr.
\item In \cite[Section 2.2]{scholze-stix} it is argued that Mochizuki's Three Indeterminacies must be irrelevant.
\item \cite[Section 2.1.4]{scholze-stix} asserts the impossibility  of the existence of the variation in arithmetic (supposed to be encapsulated by items {\bf(1)--(4)}) of the hyperplane given by the (global) product formula. 
\item According to \cite{scholze-stix}, points {\bf(2)--(4)} above lead to the core conclusion drawn in \cite[Section 2.2]{scholze-stix} and points {\bf(1)--(5)} lead to the conclusions of \cite{scholze-stix,scholze-review}.
\item Additionally, there is no mention of the geometric case in \cite{scholze-stix}. \textit{However, the geometric case is central to understanding the claims of \iut.} For a comparison of the arithmetic and the geometric case see \cref{ss:comp-with-the-geom-case}.
\eenum

The paper \cite{scholze-stix} extrapolates the absence of adequate mathematical evidence for Mochizuki's claims in \iutthr\   into a claim of mathematical impossibility of obtaining such evidence i.e. Scholze-Stix assert that no one could provide such evidence for Mochizuki's claims because it is mathematically impossible to do so.
As is noted in \cite{joshi-report}, Mochizuki's claims underlying {\bf(1)--(5)} have been independently and robustly established in the \present. As the \cref{tab:rebuttal} below shows, every ``impossibility'' assertion of \cite{scholze-stix}, \cite{scholze-review} stands invalidated by my work. Using the `Rosetta Stone' established in \constrthr{\ssep}{se:intro-rosetta-stone}, and the geometric objects provided by my theory, Mochizuki's claims may now also be verified by his methods. 
\newpage
\begin{boxedcontent}[label=tab:rebuttal,grow to right by=1cm,grow to left by=1.5cm]{{Status of various assertions of \cite{scholze-stix}, \cite{scholze-review}}}
\centering{\small
		\begin{tabular}{|p{2in}|p{2.9in}|p{1.5in}|}\hline 
			Claims of section in \cite{scholze-stix}	& Disproved by Joshi in the \present	& Mochizuki \iut\ (for comparison) \\ 
			\hline 
			Sect. 2.1.2 (and also \cite{scholze-review})
			asserts that distinct
			\benumlab 
			\item Hodge Theaters
			\item \'Etale Pictures
			\item Frobenius Pictures
			\eenum 
			do not exist & \benumlab 
			\item distinct Hodge-Theaters exist \cite[Theorem 10.11.5.1]{joshi-teich-rosetta};
			\item distinct \'etale Pictures exist \cite[Proposition 8.3.1.1]{joshi-teich-rosetta}; 
			\item distinct Frobenius Pictures exist \cite[Proposition 8.3.1.2]{joshi-teich-rosetta}
			\eenum	& asserts  the existence of all three (without proof)\\ 
			\hline 
			Remark 9	&  existence of isomorphs \cite[Theorem 2.9.1]{joshi-teich}, \cite[Theorem 4.8]{joshi-untilts} and distinct arithmetic holomorphic structures \cite[Definition 5.1]{joshi-untilts}	 & asserts existence (without proof) \\ 
			\hline 
			Section 2.1.4 & False by \cite[Theorem 5.10.1]{joshi-teich-def} & is tacitly asserted (without proof) \\
			\hline
			Section 2.1.5	& distinct prime-strips exist \cite[\ssep 8.7.1 and Theorem 8.8.3]{joshi-teich-rosetta} & asserts the existence (without proof) \\
			\hline
			Section 2.1.6 (argues $\mathfrak{log}$-Link must be irrelevant)	& $\mathfrak{log}$-Link = global Frobenius and it is constructed in \cite[\ssep 5]{joshi-teich-def}  (and is non-trivial), $\mathfrak{log}$-Link aspect is detailed in \cite[\ssep 8.9]{joshi-teich-rosetta}  & asserts properties in \cite{mochizuki-topics3} (without proof) \\
			\hline
			Section 2.1.7, 2.1.8		 & Each holomorphoid provides $q$ and $\Theta$-Pilot object \cite[Remark 6.4.2.2 and \ssep 6.10.1]{joshi-teich-rosetta}   & asserts existence (without proof) \\
			\hline
			Section 2.1.9 ($\Theta_{gau}$-Link cannot be non-trivial)		 & Construction and properties \cite[Theorem 4.2.2.1 and \ssep 4.2.3 and \ssep 8.10 ]{joshi-teich-rosetta}  & asserts non-triviality (without proof) \\
			\hline
			Section 2.2 (main disproof argument)	 & This argument of \cite{scholze-stix} fails because the $\Theta_{gau}$-Link does come with a non-trivial $j^2$-scaling factor for $j=1,\ldots,\frac{\ell-1}{2}$ and also is equipped with a non-trivial Galois action  \cite[Theorem 4.2.2.1 and \ssep 4.2.3 and \ssep 8.10 ]{joshi-teich-rosetta}  &  \\
			\hline
		\end{tabular} 
	}
	\tcblower
	{\small
		\noindent\textcolor{red}{{\bf NOTE}}
		This covers all  sections of \cite{scholze-stix} except for notations and generalities (and also covers the main assertion of \cite{scholze-review}). This table is taken from \cite[October 2024 Version]{joshi-report}.
	}
\end{boxedcontent}
\newpage
\begin{boxedcontent}[label={tab:sumtable},grow to left by=2cm, grow to right by=1cm]{Summary of the present status of \cite{scholze-stix}}{}
To summarize the above discussion of \cite{scholze-stix}, \cite{scholze-review} (for additional details see \cite{joshi-report}):
\benumlab
\item  I have examined the claims of \cite{scholze-stix} and \cite{scholze-review} in meticulous detail.  
\item At this point, all the objections regarding Mochizuki's Inter-Universal Teichmuller Theory (\iut) voiced in these reports stand completely dismantled by my work (see Table \ref{tab:rebuttal}, also \cite[\ssep 1.2]{joshi-report}).
\item Hence, the logical mathematical conclusion one may draw regarding \iut\  \textit{at this point in time}
is that the \abc\ stands established by the methods of \cite{mochizuki-iut1,mochizuki-iut2,mochizuki-iut3,mochizuki-iut4,joshi-untilts,joshi-teich-estimates,joshi-teich-rosetta,joshi-teich-abc-conj}.
\eenum
\end{boxedcontent}

\newcommand{\sV}{\mathscr{V}}
\newcommand{\bbW}{\mathbb{W}}
\newcommand{\bbWl}{\bbW_L}

\subsection{Acknowledgments}
My foray into this subject was made possible by my conversations with a few individuals I wish to acknowledge here: 

During my stay at RIMS, in Spring 2018: Many thanks to Shinichi Mochizuki who lectured to me on his view of the geometric case \cite{mochizuki-bogomolov}--his work continues to be a source of inspiration; Yuichiro Hoshi and Machiel van Frankenhuysen who answered my  questions about \topics\ and \cite{frankenhuysen2002} (and the \abc) respectively.  

Thanks to Taylor Dupuy and Anton Hilado for some  early conversations on \iut. In the Summer of 2021, and recently  in March 2024, I lectured to Taylor on my papers in great details. I thank him for the opportunity to explain some important points of Mochizuki's work (namely Mochizuki's Key Principle of Inter-Universality \cite[\ssep I3]{mochizuki-iut1}) missed by many and whose relevance to Teichmuller Theory (claimed in \iutthr) is proved using my work. 

Special thanks to Brian Conrad who encouraged me to provide a detailed proof of the main theorem of \cite[{Theorem 2.15.1, \ssep\ref{I-se:grothendieck-conj}}]{joshi-teich}. Brian's advice and help is deeply appreciated. 

Thanks to Kiran Kedlaya for his brief comments and encouragement on reading \cite{joshi-untilts-2020}. 

An early version of \cite{joshi-teich-rosetta} was sketched in Zoom lectures (Spring 2022) to my former colleagues Minhyong Kim and Dinesh Thakur and I thank them for their interest, comments and encouragement. Over the duration of this project (2018--to date), Dinesh has provided feedback, questions, suggestions and encouragement and I acknowledge my gratitude here. Thanks to Nicholas Christofferson for pointing out some typos. I thank Peter Scholze for some questions.

I met  \href{https://www.ams.org/journals/notices/202110/rnoti-p1763.pdf}{Lucien Szpiro}  for the first time in 1994 during my visit to the Max Planck Institute (Bonn). Occasionally, I would encounter him on the tram ride to our respective apartments. When Szpiro found out that my PhD advisor was N. Mohan Kumar, he gleefully exclaimed ``So, you are my grand-student!'' (Mohan Kumar wrote notes for \cite{szpiro-book1979}). While I have never entertained  any ambition of proving Szpiro's Conjecture, by some twist of mathematical fate I have been handed a role in its resolution, and I am happy for it.

\section{The conjectures and the theorems}\label{se:cons-and-thms} \nwss
\subsection{The \abc}
The following tantalizing assertion was conjectured by David Masser and Joseph Oesterle (see \cite{oesterle1988}) and is known as the $abc$-conjecture:
\bcon[The \abc]\label{con:abc} For each $\varepsilon>0$, there exists an absolute constant $C(\varepsilon)>0$, such that for all primitive triples of $a,b,c$ integers (i.e. triples of integers with $\text{gcd}(a,b,c)=1$) satisfying 
$$a+b=c,$$
one has $$\max\{\abs{a},\abs{b},\abs{c}\}\leq C(\varepsilon)\cdot \prod_{p|a\cdot b\cdot c}p^{1+\varepsilon}$$
where the product is over all the prime number $p$ dividing $a\cdot b\cdot c$.
\econ

\newcommand{\rad}[1]{{\rm rad}(#1)}

\brem\label{re:masser1} For any $abc$-triple as above, write $N=\rad{abc}=\prod_{p|a\cdot b\cdot c}p$. It is known, \cite{stewart1986} (also \cite{masser1990}) that the exponent of $N$ is optimal and cannot be weakened. More precisely, in \cite{stewart1986}, it is shown that for any $0<\delta\in \R$ and any $N_0\in \N$, there exists primitive $abc$-triples with $N=\rad{abc}\geq N_0$ such that 
$$\max\{\abs{a},\abs{b},\abs{c}\} \geq N\cdot exp{((4-\delta)\cdot \log(N)^{1/2}/\log\log(N))}.$$
\erem

\subsection{Arithmetic Szpiro Conjecture}
For an elliptic curve $E/\Q$, write $\Delta_E$ for its discriminant and $N_E$ for its conductor \cite{silverman-arithmetic}. The following was conjectured in  \cite{szpiro81}.
\bcon[The Arithmetic Szpiro Conjecture]\label{con:szpiro} Let  $0<\varepsilon\in\R$, then there exists an absolute constant $C(\varepsilon)>0$ such that for all elliptic curves $E/\Q$, one has  $$\abs{\Delta_E}_\R\ll C(\varepsilon)\cdot N_E^{(6+\varepsilon)}.$$
\econ

\brem\label{re:masser2} 
As was shown in \cite{masser1990}, the exponent of $N_E$ is optimal: for any $0<\delta\in \R$ and any $N_0\in \N$, there exists an elliptic curve over $\Q$ with $N_E\geq N_0$ such that 
$$\abs{\Delta_E} \geq N_E^6\cdot exp{((24-\delta)\cdot \log(N_E)^{1/2}/\log\log(N_E))}.$$ 
\erem

\brem 
Szpiro conjectured this inequality on the basis of his discovery of  the geometric case of his inequality in \cite{szpiro1991}. Since then a number of proofs of the Geometric Szpiro Inequality have appeared \cite{kim97}, \cite{bogomolov00}, \cite{zhang01}, \cite{beauville02}. The proofs given in \cite{bogomolov00}, \cite{zhang01} are Teichmuller Theoretic in nature and are the ones which are closest to the proof given here. Mochizuki's discussion of the geometric case is in \cite{mochizuki-bogomolov}. My discussion of the geometric case, along with a demonstration of how the structures of the sort considered in \iutthr\ and \cite{joshi-teich-def} appear in \cite[{\ssep\ref{II5-se:appendix}}]{joshi-teich-def}. Reader may find the geometric proofs useful in understanding the arithmetic case.
\erem

\subsection{Functions of bounded discrepancy class}
Let $X$ be a set (in the main context of this paper $X=\P^1-\{ 0,1,\infty\}-(\bQ)$). Let $f,g:X\to \R$ be two real valued functions. Then I will say that $f,g$ have the same  \textit{bounded discrepancy class}, and write $f\approx g$, if $\abs{f-g}:X\to \R$ is a bounded function (equivalently $f(x)=g(x)+O(1)$ holds on $X$). It is easy to see that $\approx$ is an equivalence relation on the set of real valued functions on $X$. I will  say that $f(x)\lesssim g(x)$ if $f(x)-g(x)\leq C$ for some constant $C$ (as in pointed out in \cite[Remark 2.3.1(ii)]{mochizuki-iut4}, readers should beware that the definition of $f(x)\lesssim g(x)$ given in \cite[Definition 1.2(ii)]{mochizuki-general-pos} suffers from an obvious misprint). The definition of $f\gtrsim g$ is then self-evident.

\newcommand{\logdiff}{\textrm{log-diff}}
\newcommand{\logcon}{\textrm{log-con}}
\newcommand{\hite}{{\rm ht}}
\subsection{The functions $\logdiff$ and $\logcon$}
Let $X$ be a geometrically connected, smooth projective variety over a number field $L$ and let $D\subset X$ be an effective divisor then there are two naturally defined real valued functions $\logdiff_X$ and $\logcon_D$. The definition is as follows: if $P\in X(\bQ)$ has a minimal field of definition $L(P)$ then 
$$\logdiff_X:X(\bQ)\to \R$$
is the function $$X(\bQ)\ni P\mapsto \frac{1}{[L(P):L]}\deg(\fd_{L(P)/L})$$
which associates to $P$  the normalized arithmetic degree of the arithmetic divisor given by the different $\fd_{L(P)/L}$ of $L(P)/L$. The quantity $\logdiff_X(P)$ will be called the \textit{log-different of $P$}. This depends on $X$.

Now let $D$ be an effective (Cartier) divisor on $X$. Any point $P\in (X-D)(\bQ)$, with minimal field of definition $L(P)/L$ defines a morphism $\iota_{P}:\Spec(\O_{L(P)})\to X$ ($X$ is proper). The divisor $D$ may be pulled-back along $\iota_P$. Then $$\logcon_D(P)=\frac{1}{[L(P):L]}\deg(\iota_{P}(D)_{red})$$
is the normalized arithmetic degree of the reduced divisor $(\iota_{P}^*(D))_{red}$ and $\logcon_D(P)$ will be called the \textit{log-conductor of $P$}. This obviously depends on $(X,D)$.

In this paper one considers the bounded discrepancy class of these two functions.

\subsection{Vojta's Height Inequality}
For a smooth, projective curve over a number field, let $\omega_X$ denote its sheaf of differentials. In \cite{vojta1998} Paul Vojta formulated an inequality which is now called Vojta's Height Inequality for curves. For any open subset $U\subset X$, and any integer $d\geq 1$, let $U(\bQ)^{\leq d}\subset U(\bQ)$ be the subset of points of degree $\leq d$. The following version of Vojta's Height Inequality is due to Mochizuki \cite{mochizuki-general-pos}:
\bcon[Vojta's Height Inequality for curves: Mochizuki's formulation]\label{con:vojta-inequality}
Let $L$ be a number field. Let $(X,D)/L$ be a pair consisting of a  geometrically connected, smooth, projective curve  over $L$ and $D\subset X$ a reduced divisor such that $U=X-D$ is a (punctured) hyperbolic curve over $L$. Let $d\geq 1$ be an integer. Then for every $\varepsilon>0$, the following inequality 
$$\hite_{\omega_X(D)}\lesssim (1+\varepsilon)\left(\logdiff_X+\logcon_D\right)$$
holds on $U(\bQ)^{\leq d}=(X-D)(\bQ)^{\leq d}$.
\econ
\subsection{The strong \abc}
In \cite[Conjecture 14.4.12]{bombieri-gubler}, the following conjecture is called the strong \abc:
\bcon[The strong \abc: Mochizuki's formulation]\label{con:abc-vojta}
Let $(X,D)=(\P^1,\{0,1,\infty\})$ be the \textit{fundamental hyperbolic tripod} over $\Q$. Let $d\geq 1$ be an integer. Then for every $\varepsilon>0$ the following inequality 
$$\hite_{\omega_{\P^1}(D)}\lesssim (1+\varepsilon)\left(\logdiff_{\P^1}+\logcon_D\right)$$
holds on $U_{tpd}(\bQ)^{\leq d}=(\P^1-D)(\bQ)^{\leq d}=\left(\P^1-\{0,1,\infty\}\right)(\bQ)^{\leq d}$.
\econ

\bpro
The strong \abc\ (\Cref{con:abc-vojta}) implies the \abc\ (\Cref{con:abc}).
\epro
\bp 
Taking $d=1$ one obtains the assertion.
\ep

\subsection{van Frankenhuysen's Theorem}
In \cite{vojta1998}, it had been shown that Vojta's Height Inequality for Curves (\Cref{con:vojta-inequality}) implies the \abc\ (\Cref{con:abc}). The converse of this assertion is the following theorem
of \cite{frankenhuysen2002} which is of importance to this paper:
\bthm\label{th:frankenhuysen}
The following conjectures are equivalent:
\benumlab
\item The strong \abc\ (\Cref{con:abc-vojta}).
\item Vojta's Height Inequality for Curves (\Cref{con:vojta-inequality}).
\eenum
\ethm

\subsection{Mochizuki's Refinement of van Frankenhuysen's Theorem}
In \cite[Theorem 2.1]{mochizuki-general-pos}, Mochizuki proved the following important refinement of \cite{frankenhuysen2002} which asserts that a weaker version of the strong \abc\ implies Vojta's Height Inequality for Curves (\Cref{con:vojta-inequality}) and hence the \abc\ (\Cref{con:abc-vojta}). More precisely, the validity of Vojta's Inequality  (\Cref{con:vojta-inequality}) on compactly bounded subsets of $\P^1-\{0,1,\infty \}(\bQ)$ implies (\Cref{con:vojta-inequality}). For compactly bounded subsets see \cref{ss:compact-bnded}.

\bthm[\cite{mochizuki-general-pos}] 
Let $S$ be a finite set of primes including archimedean primes. Let $d\geq1$ be an integer. Then the following are equivalent:
\benumlab 
\item Vojta's Height Inequality for Curves (\Cref{con:vojta-inequality}) holds. 
\item \Cref{con:abc-vojta} holds for any $S$-supported compactly bounded subset of $U(\bQ)^{\leq d}=\P^1-\{0,1,\infty \}(\bQ)^{\leq d}$. 
\eenum
\ethm

\brem\label{re:wada}
In the context of the above refinement, let me remark that there is an analog of examples discussed in \cref{re:masser1} for compactly bounded subsets i.e. the exponent of the conductor in the \abc\ cannot be improved even on compactly bounded subsets. This is established in \cite[Theorem 0.5]{wada2016} and this also implies the existence (\cite[Theorem 2.7]{wada2016}) of examples of Frey curves with $j$-invariants lying in compactly bounded subsets and satisfying an estimate similar to that given in  \cref{re:masser2} i.e. the exponent of the conductor in the Arithmetic Szpiro's conjecture cannot be improved even on compactly bounded subsets.
\erem

\subsection{A brief history}
\benumlab
\item In \cite{szpiro81}, Szpiro proved the geometric Szpiro Inequality and this led him to conjecture the inequality now know as the Arithmetic Szpiro Inequality (\Cref{con:szpiro}). Several proofs of the geometric Szpiro Inequality now exist \cite{szpiro81}, \cite{kim97},  \cite{bogomolov00}, \cite{zhang01}, \cite{beauville02}.
\item Mordell's Conjecture was proved by Faltings in \cite{faltings1983}. 
\item In \cite{oesterle1988}, it was observed that the $abc$-conjecture (\Cref{con:abc}) is equivalent to the arithmetic Szpiro Conjecture (\Cref{con:szpiro}). 
\item In \cite{elkies1991} it was shown that the \abc\ (\Cref{con:abc}) implies Mordell's Conjecture.
\item In \cite{vojta1998} conjectured a general inequality, which is now known as Vojta's Height Inequality  for Curves (\Cref{con:vojta-inequality}), and showed that this inequality implies the $abc$-conjecture (\Cref{con:abc}). 
\item The equivalence of the two forms of the \abc\  given by \Cref{con:abc} and \Cref{con:abc-vojta} is quite elementary. 
\item In an important development, in \cite{frankenhuysen2002} it was shown that the \abc\ is equivalent to  Vojta's Height Inequality (\Cref{con:vojta-inequality}).
\item In \cite[Theorem 2.1]{mochizuki-general-pos}, Mochizuki refined the method of \cite{frankenhuysen2002} and showed that  Vojta's Height Inequality for compactly bounded subsets of $(\P^1-\{0,1,\infty\})(\bQ)$ (\Cref{con:abc-vojta}) implies Vojta's Height Inequality (\Cref{con:vojta-inequality}) and by van Frankenhuysen's Theorem he deduced the \abc\ (\Cref{con:abc}) from \abc\ for compactly bounded subsets of a fixed curve. 
\eenum

\subsection{The Main Theorem}
The following assertion of \cite{mochizuki-iut4}, together with the \present\ or \iutthr\ together with the corrections and changes suggested in the `Rosetta Stone' of \cite[{\ssep\ref{III-se:intro-rosetta-stone}}]{joshi-teich-rosetta},  establish:
\bthm\label{th:main-thm}
Let $(X,D)/L$ be a geometrically connected, smooth, projective curve over a number field $L$ and let $D\subset X$ be a reduced divisor over $L$. Let $d\geq 1$ be an integer. Then for every $\varepsilon>0$, Vojta's Height Inequality $$\hite_{\omega_X(D)}\lesssim (1+\varepsilon)\left(\logdiff_X+\logcon_D\right)$$
holds on $U(\bQ)^{\leq d}=(X-D)(\bQ)^{\leq d}$.

As a consequence the following assertions are also true:
\benumlab
\item The strong  \abc\  (\Cref{con:abc-vojta}) (and \abc).
\item The Arithmetic Szpiro Inequality (\Cref{con:szpiro}) for elliptic curves over number fields.
\item Effective form of Mordell's Conjecture is true for $X/L$.
\eenum
\ethm
\brem 
The assertion that the effective \abc\ implies the effective Mordell's conjecture  is established in \cite{elkies1991} (also see \cite{frankenhuysen2002}). The strategy of the proof of \Cref{th:main-thm} is to prove \Cref{th:main-thm}{\bf(2)} for a given compactly bounded subset supported on a suitable finite set of primes.
\erem

\section{Averaging over Arithmetic Teichmuller Space of a Number Field}\label{se:averaging}\nwss

\subsection{Anabelomorphisms}
\emph{Throughout this paper a $p$-adic field means a finite extension of $\Q_p$ for some prime number $p$. All valuations on any valued field considered here will be rank one valuations.} 

Throughout the paper I will freely use the notion  of \emph{anabelomorphisms of schemes, fields etc.}  and the related notion of \emph{amphoric quantities, properties and structures} associated to anabelomorphic schemes, fields introduced in \cite[{\ssep\ref{0-se:anabelomorphy}}]{joshi-anabelomorphy}. Some basic results can be found in \cite[{\ssep\ref{0-se:five-fundamental}}]{joshi-anabelomorphy}. 
\newcommand{\sC}{\mathcal{C}}

\subsection{Averaging over additive structures of a number field}
It is asserted in \cite{mochizuki-iut4} that  proof of the \abc\ (and Arithmetic Szpiro Conjecture) presented therein rests on averaging over distinct additive structures i.e. one averages over distinct avatars of the arithmetic of a fixed number field. In \iutthr, the existence of distinct avatars is asserted without adequate and convincing proof while its theory rests, centrally on a precise quantification of this notion. 

The Theory of Arithmetic Teichmuller Spaces of the \present, allows me to demonstrate this central missing point  of \iut\ in exquisite detail, and with far greater precision than that can be achieved by the methods of \iut. 

My work provides a topological space  $\yadl$ (constructed in \cite[{\ssep\ref{II5-se:adelic-ff-curves}}]{joshi-teich-def}) which parameterizes  the many inequivalent avatars of a fixed number field. The central idea which emerges is that the arithmetic structure of a number field is bendable while its multiplicative structure of the number field remains fixed, but the additive structures and the intertwining between the additive and multiplicative structures may differ in a precise topological sense.  The differences between two such structures move in a manner compatible with the differences between two such structures at  completions at  each prime of the number field \cite[{\ssep\ref{II5-se:adelic-ff-curves}}]{joshi-teich-def}.

The question of great importance,  which arises is this: 

\begin{verse}
	How does one measure the differences between two additive structures of a number fields, $p$-adic fields?
\end{verse}

Understanding the answer to this question, holds the key to the proof of \cite{mochizuki-iut4},  (\iut\ does not explicitly provide a direct quantitative answer to this question).

\subsection{Differents, discriminants, and Measuring differences between two additive structures} 
A direct quantitative answer to the above question was provided in \cite[{\ssep\ref{0-se:disc-unamph}}]{joshi-anabelomorphy} for strict anabelomorphisms of $p$-adic fields and for any hyperbolic curve viewed over two such  $p$-adic over-fields. The answer is quite elegant:
 \begin{verse}
 \textit{Discriminants, differents, Swan conductors are measures of additive structure of a Number Field and of Varieties over a Number Field!}
 \end{verse}
 
A quantitative version of the answer for $\yadl$ is the theorem given below.

\newcommand{\frakd}{\mathfrak{d}}
\bthm\label{th:diff-add-var} Let $L$ be a number field with no real embeddings. Let $X/L$ be a geometrically connected, smooth hyperbolic curves over $L$. For $\by_1,\by_2\in\yadl$, let $\arith{L}_{\by_1}^{nor}, \arith{L}_{\by_2}^{nor}$ be two normalized arithmeticoids of $L$ (\cite[{\ssep\ref{II5-se:arithmeticoid-adeloid-frobenioid}}]{joshi-teich-def}). Let $\bL_{\by_1}, \bL_{\by_2}$ be the algebraic closures of $L$ provided by the two arithmeticoids. Let $M/L$ be a finite, Galois extension of $L$ viewed in the algebraic closures provided by the two arithmeticoids $M_{\by_i}\subset\bL_{\by_i}$ (for $i=1,2$). Let $v\in\V_L$ and $w|v$ be a prime of $M$ lying over $v\in\V_L$. Then 
\benumlab
\item The Differents $\frakd_{M_w/L_v;\by_i}$, Discriminants $d_{M_w/L_v;\by_i}$, Artin and Swan conductors of ${M_w/L_v;\by_i}$, and the discriminants, the Artin and the Swan conductors  of hyperbolic curves $X/\arith{L}_{\by_i}$ are dependent on the additive structure of $\arith{L}_{\by_j}$ and hence are measures of the relative differences between the additive structures of two arithmeticoids $\arith{L}_{\by_j}$.
\item For $v\in\vl$, the values $d_{M_{j,w}/L_v}, \fd_{M_{j,w}/L_v}\in L_{v,j}^*\subset K_{v,j}$ may not be directly compared as the two arithmeticoids may not be comparable (in the sense of \cite[{\Cref{II5-th:anabelomorphic-deforms}}]{joshi-teich-def}). 
\item A choice of an anabelomorphism $G_{L_v,K_{y_{v,1}}}\isom G_{L_v,K_{y_{v,1}}}$ provides an isomorphism of topological groups 
$$L_{v,1}^*\isom  L_{v,2}^*$$ and an isomorphism of the unit groups $$\O_{L_{v,1}}^*\isom  \O_{L_{v,2}}^*$$ 
\item Using this isomorphism, the respective discriminants and differents $$d_{M_{j,w}/L_v}\text{ and }\fd_{M_{j,w}/L_v}\in L_{v,j}^*$$ resp. $$\frac{d_{M_{j,w}/L_v}}{p_v^{ord_v(d_{M_{j,w}/L_v})}} \text{ and }\frac{\fd_{M_{j,w}/L_v}}{p_v^{ord_v(\fd_{M_{j,w}/L_v})}}\in \O_{L_{v,j}}^*$$ may be compared.
\eenum
In particular $$\fd_{M/L;\by}=(\fd_{M_w/L_v,\by})_{v\in\vl}\in \prod_{p\in\vl}K_{y_v}^*$$ depends on the arithmeticoid $\arith{L}_{\by}$ for $\by\in\yadl$, and $\by\mapsto \fd_{M/L;\by}$ is a non-constant function on $\yadl$.
\ethm
\bp This was established in \cite[{Theorem \ref{0-th:discriminant-is-unamphoric}}]{joshi-anabelomorphy} for the case of distinct anabelomorphic $p$-adic fields and for curves over such  fields. However, this principal underlying the proof in that case is equally valid for distinct arithmeticoids for the same $p$-adic field $L_v$ being considered here. 

The central point is the following. Let $G_i=G_{L_v;K_{y_{v,i}}}$ be the absolute Galois groups of $L_v$ computed using the algebraic closures of $L_v$ provided by the two arithmeticoids.  Let $G_i^\mydot=G_{L_v;K_{y_{v,i}}}\subset G_i$ be the upper numbering ramification filtration. As pointed out in \cite[{\ssep\ref{0-se:intro}}]{joshi-anabelomorphy}, this filtration is a stand-in for the additive structure of $L_v$ (thanks to \cite{mochizuki-local-gro}). 

Since the two arithmeticoids are distinct in general, there is no natural isomorphism of the filtered groups $(G_i,G_i^\mydot)$ and the normalization of  the arithmeticoid  $\arith{L}_{\by_1}^{nor}$ does not carry over to $\arith{L}_{\by_2}^{nor}$ i.e. both may not be simultaneously normalizable.

On the other hand, the two discriminants, differents and  Swan conductors are dependent on the ramification filtration. Hence these invariants, computed for one of the arithmeticoids, can be numerically distinct from the corresponding invariants for the other arithmeticoid.  The simplest way of understanding this assertion is to recognize that the function on Fargues-Fontaine curve given by $\syQp\ni y\mapsto \abs{p}
_{K_y}$ is a highly non-trivial function.
\ep

\brem 
Let me remark that \Cref{th:diff-add-var} is not useful in practical computations or estimates. Nevertheless it is of importance in understanding how the choices of additive structures of a number field impacts the geometric and numerical invariants of interest.
\erem

\newcommand{\bz}{{\bf z}}

\newcommand{\bsi}{\boldsymbol{\mathscrbf{I}}}
\newcommand{\bsiq}{\boldsymbol{\mathscrbf{I}^{\Q}}}
\newcommand{\blogsh}[1]{\bsi(#1)}
\newcommand{\blogshq}[1]{\bsi^{\Q}(#1)}

\newcommand{\sImj}{\mathscrbf{I}_{\scriptscriptstyle{Joshi}}}
\newcommand{\sImjq}{\mathscrbf{I}_{\scriptscriptstyle{Joshi}}^{\Q}}
\newcommand{\tsIm}{\widetilde{\mathscrbf{I}}_{\scriptscriptstyle{Mochizuki}}}
\newcommand{\tsImj}{\widetilde{\mathscrbf{I}}_{\scriptscriptstyle{Joshi}}}
\newcommand{\tsImq}{\widetilde{\mathscrbf{I}}_{\scriptscriptstyle{Mochizuki}}^{\Q}}
\newcommand{\tsImjq}{\widetilde{\mathscrbf{I}}_{\scriptscriptstyle{Joshi}}^{\Q}}
\newcommand{\bsIm}{\boldsymbol{\mathscrbf{I}}_{\scriptscriptstyle{Mochizuki}}}
\newcommand{\bsImj}{\boldsymbol{\mathscrbf{I}}_{\scriptscriptstyle{Joshi}}}
\newcommand{\bsImq}{\boldsymbol{\mathscrbf{I}}_{\scriptscriptstyle{Mochizuki}}^{\Q}}
\newcommand{\bsImjq}{\boldsymbol{\mathscrbf{I}}_{\scriptscriptstyle{Joshi}}^{\Q}}
\newcommand{\Vol}{{\rm Vol}}

\section{Differents and Tate Divisors}\label{se:diff-tate-div}\nwss
\subsection{The basic set up}\label{ss:basic-setup}
Recall the theory of \cite{joshi-teich-def}. Let $\arith{L}_{\by}$ be the standard arithmeticoid of a number field $L$. This arithmeticoid is defined using the chosen standard point $\by_0=(y_{v})_{v\in\vl}\in\yadl'$  defined in \cite[{\ssep\ref{II5-se:arith-and-galois-cohomology}}]{joshi-teich-def}. Let $C/\arith{L}_{\by}$ be a holomorphoid of an elliptic curve $C/L$. This choice provides for us an algebraic closure of $\bL$ of $L$ and an algebraic closure of $\bL_v$ for each $v\in\vl$ and thus provides us a way of computing Galois groups, modules of roots of unity, Galois cohomology, arithmetic invariants of elliptic curve etc. 

The assumptions \inithtdata\ and \assumptions\ are valid for each of these holomorphoids $C/\arith{L}_{\by}$. As was established in \cite[{\ssep\ref{II5-se:arithmeticoid-adeloid-frobenioid}}]{joshi-teich-def}, each arithmeticoid $\arith{L}_{\by}$ represents a deformation of the arithmetic of $L$. 

Let $\bL$  be the algebraic closures of $L$ provided by  the  arithmeticoid of $\arith{L}_{\by}$. And let $M/L$ be a finite extension of $L$ and let $L\subset M \subset\bL$ be this extension viewed in the two arithmeticoids. Let $v\in\vl$ be a prime and let $w|v$ be a prime of $M$ lying over $v$. I work with this choice of arithmeticoid for all the calculations. As has been pointed out, the theory of the \present\ works simultaneously with a large collection of arithmeticoids and one averages Tate parameters and other quantities associated with respect to these arithemeticoids.

Notations and assumptions of \initassumptions\ are strictly in force.
\subsubsection{}\label{sss:degrees}
The following assumptions and notations are added to \initassumptions:
\benumlab
\item $C/L$ has good reduction at $v\in\vl^{good}\cap\vlnon\cap\{v\in \vl:v \text{ does not divide }(2\cdot\ell) \}$.
\item Let $\dmod=[\lmod:\Q]$,
\item let $\dmods=2^{12}\cdot 3^3\cdot 5\cdot\dmod$.
\item Assume that $$\emod=\max_{v\in\V_{\lmod}^{non}}\left\{ e_v \right\}\leq \dmod,$$ where $e_v$ is the absolute ramification index of $v\in\V_{\lmod}^{non}$. 
\item and let $\emods=2^{12}\cdot 3^3\cdot 5\cdot \emod\leq\dmods$.
\eenum

\subsubsection{}\label{sss:field-assump}
In addition to \cref{sss:degrees}, also assume the following:
\benumlabresume
\item $\ell\geq 5$ is a prime (to be chosen to be sufficiently large later on).
\item $C_L(L)\supset C_L[15](\bL)$, and
\item $\lmod\subseteq \ltpd=\lmod(C_{\lmod}[2](\bL))\subseteq L$, and
\item $L=\ltpd(\sqrt{-1},C_{\ltpd}[3\cdot5])$ this implies that $C_L$ is extended from $C_{\ltpd}$.
\eenum
The following consequence will be used throughout:
\bpro\label{pr:basic-ramification-data}
The assumptions  \initassumptions\ and \assumptions\ are in force. Let $v\in\vnon_L$. Then
\benumlab
\item $C/L$ has semistable reduction at $v$ if and only if $v\in\ubblvoss$.
\item $C/L$ has good or additive reduction if $v\in \{w\in\vl: w|(2\ell)\}$.
\item $C/L$ has good reduction at $v\not\in\left(\ubblvoss\cup \{w\in\vl: w|(2\ell)\}\right)$.
\eenum
\epro
\bp 
This is clear from \assumptions\ and \initassumptions.
\ep
\blem 
The extension $L/\lmod$ is a composite extension of Galois extensions,, it is Galois with  $$\gal(L/\lmod)\into \Z/2\times\glt{\F_2}\times\glt{\F_3}\times\glt{\F_5}.$$
\elem
\bp 
By \assumptions, $L/\lmod$ is a composite of Galois extension and the assertions follows from \cite[Chap VI, Theorem 1.14]{lang-algebra}.
\ep
\brem
Thus $C/L$ is special elliptic curve. The existence of such an elliptic curve is a delicate argument of \cite{mochizuki-iut4} using \cite{mochizuki-general-pos}.
\erem

\blem\label{le:tame-lemma} 
\benumlab 
\item The extension $L/\ltpd$ is unramified outside $$\{v\in \V_{\ltpd}: v\text{ divides } 2\cdot 3\cdot 5\cdot\ell\}\cup \supp(\fq_{\ltpd})$$ and  $L/\ltpd$ is tamely ramified outside $$\{v\in \V_{\ltpd}: v\text{ divides } 2\cdot3\cdot5\}.$$
\item The extension $L'/L$ is unramified outside $$\{v\in \V_{\ltpd}: v\text{ divides } 2\cdot\ell\}\cup \supp(\fq_{\ltpd})$$ and  $L'/L$ is tamely ramified outside $$\{v\in \V_{\ltpd}: v\text{ divides } 2\cdot\ell\}.$$
\eenum
\elem
\bp 
This is immediate from \Cref{pr:basic-ramification-data}.
\ep
\newcommand{\fdiff}{\mathfrak{diff}}
\newcommand{\fdisc}{\mathfrak{disc}}

\renewcommand{\fa}{\mathfrak{a}}

\subsection{Some more notations}
Let $L'\supseteq M\supseteq\lmod$ be a finite Galois extension. Let $v\in\vmnon$ be a prime. Let $w|v$ with $w\in\vnon_M$, write $f_w,f_v$ for the absolute residue  field degree of $M_w, \lmodv$ and $f_{w|v}$ for the relative residue  field degree of $M_w/\lmodv$. Then one has $f_w=f_{w|v}\cdot f_v$. Write $\log w=\log(p_v^{f_w})$ (resp. $\log v=\log(p_v^{f_v})$) for the logarithm of the cardinality of the residue field of $M_w$ (resp. $\lmodv$). Then one has $\log w=f_{w|v}\cdot \log v$.  Let $e_w$ be the absolute ramification degree of $M_w/\Q_p$. Let $e_{w|v}$ be the relative ramification index of $M_w/\lmodv$.

\subsubsection{} Assumptions \initassumptions, \cref{sss:degrees}, \cref{sss:field-assump} are strictly in force. Let $$L'\supset M\supset L$$ be the number field, where $L'$ is the number field determined by \initassumptions. Let $$\fa_M=\sum_{v\in\V_M} c_v\cdot v,$$ where $c_v\in \R$ be an \textit{arithmetic divisor on $M$} (this means  $c_v=0$ for all but finitely many $v\in\V_M$). Write 
\be\label{def:deg} \deg(\fa_M)=\sum_{v\in\vnon_M} c_v\cdot \log v+\sum_{v\in\varc_M} c_v\in\R,\ee
where $\log v$ denotes the cardinality of the residue field of $M_v$. One calls
$\deg(\fa_M)$ the \textit{arithmetic degree} 
and then
\be\label{def:logdeg} \log(\fa_M)=\frac{\deg(\fa_M)}{[M:\Q]},\ee
is called the \textit{normalized arithmetic degree} of $\fa_M$.

For $M=L$, I will simply write $\fa,\deg(\fa),\log(\fa)$ instead of $\fa_L,\deg(\fa_M),\log(\fa_M)$ etc. 

\subsection{The Different and its arithmetic degree}\label{ss:diff-and-log-diff}
Let $M/\Q$ be a finite extension. Then 
\be 
\fd_M=\sum_{w\in\V_M} {\rm ord}_w(\fd_{M,w})\cdot w,
\ee
where \be \fd_{M,w} ={\rm ord}_w(\fdiff_{\O_{L,w}/\Z_p})\ee is the order of any generator of the different ideal $\fdiff_{\O_{L,w}/\Z_p}$ of $\O_{L,w}/\Z_p$. 
Let $$\log(\fd_M)=\frac{1}{[M:\Q]}\sum_{w\in \V_M}\fd_{M,w}\cdot \log w.$$ be its \textit{normalized arithmetic degree.}
Let $$\fdisc_{M_w/\Z_p}=Norm_{M_w/\Q_p}(\fdiff_{\O_{L,w}/\Z_p})$$
be the discriminant ideal.
Then the discriminant ideal of $M/\Q$ is given by the arithmetic divisor  on $\Q$
\be \fdisc_{M/\Q}=\sum_p{\rm ord}_p(\fdisc_{M_w/\Z_p})\cdot p
\ee
and let
\be 
\log\fdisc_{M/\Q}=\frac{1}{[M:\Q]}\sum_p{\rm ord}_p(\fdisc_{M_w/\Z_p})\cdot\log(p)
\ee
its normalized arithmetic degree.
\bpro\label{re:diff-rmk}
Let $M/\Q$ be a finite extension.  Then one has the equality
$$\log\fd_M=\log\fdisc_{M/\Q}.$$
\epro
\bp 
This is immediate from the relationship between  differents and discriminants of fields \cite[Proposition B.1.19]{bombieri-gubler}  and the proof of this assertion is spelled out in \cite[14.3.8, Proposition B.1.19]{bombieri-gubler}.
\ep

\subsection{The Tate divisor and the reduced Tate divisor}\label{ss:tate-divisor}\nwss
Let the notations and assumptions \assumptions, \initassumptions\ be in force. Let $L'\supseteq M\supseteq\lmod$ be a finite Galois extension. Let $v\in\vmnon$ be a prime. Let $w|v$ with $w\in\vnon_M$, write $f_w,f_v$ for the absolute residue  field degree of $M_w, \lmodv$ and $f_{w|v}$ for the relative residue  field degree of $M_w/\lmodv$. Then $f_w=f_{w|v}\cdot f_v$; write $\log w=\log(p_v^{f_w})$ (resp. $\log v=\log(p_v^{f_v})$) for the logarithm of the cardinality of the residue field of $M_w$ (resp. $\lmodv$); then $\log w=f_{w|v}\cdot \log v$.  Let $e_w$ be the absolute ramification degree of $M_w/\Q_p$. Let $e_{w|v}$ be the relative ramification index of $M_w/\lmodv$. 

Let $C_M=C\times_LM$ be the base extension of $C$ to $M$. Let $\ubblvoss_M$ be the inverse image of the set of primes   $\ubblv_{\lmod}\subset \V_{\lmod}$  of $\lmod$ (of odd residue characteristics by assumption) over which $C_{\lmod}$ has semi-stable reduction. 

\brem Beware that since additive reduction is unstable under finite extensions and may become multiplicative reduction \cite{silverman-arithmetic},  some prime $w\in\vnon_{M}$ of semi-stable reduction of $C_{M}$ may lie over a prime $v\in\vnon_{\lmod}$ not contained in $\ubblv_{\lmod}$. Such a $w\not\in\ubblvoss_M$.
\erem

\begin{defn}
If $w\in\ubblvoss_M$ and let  $q_w$ at $C_{M_w}$ be a Tate parameter (well defined up to multiplication by a local unit and its order at $w$). The arithmetic divisor 
Let 
$$
\fq_M= \sum_{w\in \ubblvoss_M} \ord_w(q_w)\cdot w.
$$
will be called the \textit{Tate divisor of $C_M$} and 
$$
\ff_M= (\fq_M)_{red}=\sum_{w\in \ubblvoss_M} 1\cdot w
$$
will be called the \textit{reduced Tate divisor of $C_M$}. One has the arithmetic degree of the Tate and reduced Tate divisors:
$$
\log(\fq_M)= \frac{1}{[M:\Q]}\cdot\deg(\fq_{M}),
$$
of $\fq_M$ and 
$$
\log(\ff_M)= \frac{1}{[M:\Q]}\cdot\deg(\ff_{M})
$$
of $\ff_M$ respectively.
\end{defn}

\brem 
In \cite[Definition 3.3]{mochizuki-general-pos}, the number $\ord_w(q_w)$ is called the \textit{local height} of the elliptic curve at $w$.
\erem

The following fundamental property of the Tate divisor is asserted without proof in \cite{mochizuki-iut4} (no proof is given in \cite{fucheng} or \cite{yamashita}).

\bpro\label{pr:tate-divisor} 
If $L'\supseteq M\supseteq\lmod$ is a finite extension, then
$$\log(\fq_M)=\log(\fq_{\lmod}).$$ 
\epro
\bp 
One has by definition
\be \log(\fq_M)=\frac{1}{[M:\Q]}\cdot\left(\sum_{w\in \ubblvoss_M} \ord_w(q_w)\cdot \log w\right)\ee 
Let $\pi_w$ and $\pi_v$ be the uniformisers of $\O_{M_w}$ and $\O_{\lmodv}$ respectively. 
If  $q_v$ is a Tate parameter of $C_{\lmodv}$, then for some $\alpha\in\Q$ one has $q_v=\pi_v^\alpha$.
Let $e_{w|v}$ be the relative ramification index of $M_w/\lmodv$  then $q_w=q_v$ is also the Tate parameter over $\O_{M_w}$. Hence as $\pi_w=\pi_v^{e_{w|v}}$ one sees that $$\ord_w(q_w)=\ord_w(q_v)=\ord_w(\pi_v^\alpha)=\ord_w(\pi_w^{e_{w|v}\cdot\alpha})=e_{w|v}\cdot \ord_v(q_v).$$

If $f_{w|v}$ is the relative residue field degree of $M_w/\lmodv$ then
$$\log w=f_{w|v}\cdot \log v.$$
Hence  
\begin{align}
 \log(\fq_M) & = \frac{1}{[M:\Q]}\cdot\left(\sum_{v\in \ubblvoss_{\lmod}}\left(\sum_{w|v}  \ord_w(q_w)\cdot \log w\right)\right)\\
 &=\frac{1}{[M:\Q]}\cdot\left(\sum_{v\in \ubblvoss_{\lmod}}\sum_{w|v}  e_{w|v}\cdot\ord_v(q_v)\cdot f_{w|v}\cdot\log v\right)\\
 &=\frac{1}{[M:\Q]}\cdot\left(\sum_{v\in \ubblvoss_{\lmod}} \ord_v(q_v)\cdot \log v\cdot \left(\sum_{w|v}  e_{w|v}\cdot f_{w|v}\right)\right)\\
 &=\frac{1}{[\lmod:\Q]}\cdot\left(\sum_{v\in \ubblvoss_{\lmod}} \ord_v(q_v)\cdot \log v\right),
 \end{align}
 as $$\sum_{w|v} e_{w|v}\cdot f_{w|v}=[M:\lmod] \text{ and } [M:\Q]=[M:\lmod]\cdot[\lmod:\Q].$$ This completes the proof of the proposition.
\ep

\subsection{Idelic Tate divisors}
The key idea in the proof of \Cref{th:main-dioph-thm} is to average over  Tate divisors (more generally a similar divisor defined using theta-values) of an elliptic curve
$C/L$ with respect to many distinct arithmetic holomorphic structures or holomorphoids $C/\arith{L}$.  For this it is more convenient to work with the idele given by a choice of Tate parameters at each $v\in\vlnon$:
\newcommand{\ideles}{\prod_{v\in\vl} L_v^*}
\newcommand{\tideles}{\prod_{v\in\vl} \widetilde{L}_v^*}
\newcommand{\TI}{\mathscr{T\kern-2mmI}}
\newcommand{\tTI}{\widetilde{\TI}}
\be\TI_{C/\arith{L}}=(q_v)_{v\in\vlnon}\in \ideles\ee
where 
\be q_v = \begin{cases}
 1 & \text{ if } v\in\vl-\ubblv,\\
 q_v & \text{ a choice of the Tate parameter at } v \text{ if }v\in\ubblv.
\end{cases}  
\ee
One may similarly define, for $j=1,\ldots,\ells=\frac{\ell-1}{2}$, the tuple of Tate divisors
\be\TI^j_{C/\arith{L}}=(q^{j^2}_v)_{v\in\vlnon}\in \ideles\ee
where 
\be q_v = \begin{cases}
	1 & \text{ if } v\in\vl-\ubblv,\\
	q_v^{j^2} & q_v \text{ a choice of the Tate parameter at } v \text{ if }v\in\ubblv.
\end{cases}  
\ee
I will continue to write $\TI=\TI^1$ interchangeably.

Then $\TI_{C/\arith{L}}$ is an idele of $L$ (for $\arith{L}$) and it  provides an arithmetic divisor
\be{\rm Divisor}(\TI_{C/\arith{L}})=\sum_{v\in\vlnon}\ord_v(q_v)\cdot v=\fq_{C/L}.\ee
In other words, the arithmetic divisor associated to $\TI_{C/\arith{L}}$ is none other than the Tate divisor $\fq_{C/\arith{L}}$ of $C/\arith{L}$ defined in \cref{ss:tate-divisor} (i.e. the Tate divisor of $C$ for the arithmeticoid $\arith{L}$. There are many choices of  $\TI_{C/\arith{L}}$ (with $\arith{L}$ fixed), but different choices of $\TI_{C/\arith{L}}$ provide the same Tate divisor $\fq_{C/\arith{L}}$. Readers should beware that the dependence of on the choice of $\arith{L}$ needs to be emphasized. 

Similarly the sum
\be\sum_{j=1}^\ells{\rm Divisor}(\TI^j_{C/\arith{L}})=\sum_{j=1}^\ells\left(\sum_{v\in\vlnon}\ord_v(q^{j^2}_v)\cdot v\right).\ee
is the divisor associated to the tuple of Theta-values in \cite{mochizuki-iut3}. One is interested in bounding its degree. This is not directly possible at the moment. But it is possible to work with these divisors in Galois cohomology. This is worked out below for $\TI$ (the construction of $\thetami$ given in) \cite[{\ssep\ref{III-se:mochizuki-construction-thetam}}]{joshi-teich-rosetta} provides a more detailed view of this).

For passage to Galois cohomology it will be convenient to work with a variant $\tideles$ of $\ideles$ where \be\widetilde{L}_v^*=\invlim_{n} L_v^*/(L_v^*)^{p_v^n},\ee
(where $p_v$ is the residue characteristic of $v$). 
Note that, by Kummer Theory, one has the natural identification 
\be \widetilde{L}_v^*\isom H^1(G_{L_v},\Q_{p_v}(1)).\ee
So in the notation and terminology of \cite[Section 7, \ssep 7.2]{joshi-teich-def},
\be\tideles=H^1(\arith{L},\Q(1)).\ee
A compatible system, $\widetilde{q}_v$, of $p_v^n$-roots of $q_v$ provides an element of the 
galois cohomology of the arithmeticoid
\be\widetilde{\TI}_{C/\arith{L}}=(\widetilde{q}_v)_{v\in\vlnon}\in \tideles=H^1(\arith{L},\Q(1)).\ee

The following proposition  allows us to collate $\tTI_{C/\arith{L}_\bz}$ arising from many different arithmeticoids in the Galois cohomology of the chosen arithmeticoid $\arith{L}_{\by}$, and will clarify the volume computations which are carried out here and in \cite{mochizuki-iut4} and their relationship to the role distinct arithmetic holomorphic structures play in the proof of the \abc\ and \Cref{th:main-dioph-thm}.

\bpro\label{pr:collation-principle} 
Let $C/L$ be an elliptic curve over $L$. Let $\sA\subset\yadl$ be a collection of arithmeticoids and write $$\sA_C=\{C/\arith{L}_\bz:\bz\in \sA\}$$ for the collection of holomorphoids of $C/L$ provided by the arithmeticoids in $\sA$. For each holomorphoid $C/\arith{L}_\bz\in\sA_C$ let $$\tTI_{C/\arith{L}_\bz}\in  H^1({\arith{L}_\bz},\Q(1))$$ be a choice of Tate idele cohomology class arising from $C/\arithl_\bz$ for ${\bz\in\sA}$. Let $\by=\by_0$ be the standard arithmeticoid chosen earlier. Then there exists a subset 
$$\Psi_{C,\sA}\subseteq  H^1({\arith{L}_{\by_0}},\Q(1))$$
which consists of the union,  over $\bz\in\sA$, of the images of all the subsets $\Psi_\bz$ under all the  isomorphisms (of topological groups) of each factor of $$H^1(\arith{L}_\bz,\Q(1))\isom  H^1(\arith{L}_{\by},\Q(1)).$$
\epro

\bp 
This is immediate on taking $\Psi_{C,\sA}=\{ \tTI_{C/\arith{L}_\bz}:\by \in\Psi\}$ in \cite[{Proposition \ref{II5-pr:collation-principle}}]{joshi-teich-def}.
\ep

The following is a variant in which one replaces adelic computation in the above proposition by a direct computation of normalized arithmetic degrees.

\bpro\label{pr:collation-principle2} 
Let $C/L$ be an elliptic curve over $L$. Let $\sA\subset\yadl$ be a collection of arithmeticoids and write $$\sA_C=\{C/\arith{L}_\bz:\bz\in \sA\}$$ for the collection of holomorphoids of $C/L$ provided by the arithmeticoids in $\sA$. For each holomorphoid $C/\arith{L}_\bz\in\sA_C$, let $$\TI_{C/\arith{L}_\bz}$$ be a choice of Tate idele  arising from $C/\arithl_\bz$ for ${\bz\in\sA}$. Let $\arith{L}_{\bz}^{nor}$ be the normalized arithmeticoid of $L$. Then the normalized arithmetic degree $\log(\TI_{C/\arith{L}_\bz})$ is a function of $\bz$. In particular if one takes $\sA=\yadl$, then one obtains a non-constant function ``normalized degree of the Tate divisor'' 
$$\yadl\to \R$$
given by $$\yadl\ni\bz\mapsto \log(\TI_{C/\arith{L}_\bz}).$$
\epro

\bp 
This is clear from the previous proposition.
\ep

\brem\  
\benumlab
\item Of interest to us is the volume of $\Psi_{C,\sA}\subseteq  H^1({\arith{L}_{\by_0}},\Q(1))$ for a suitable measure on the latter vector spaces. A fancier version of this is a sort of calculation carried out  in the proof the main theorem of this paper.
\item It is difficult to understand what the image of a specific element of $\Psi_{C,\sA}\subseteq  H^1({\arith{L}_{\by_0}},\Q(1))$ looks like. If the standard arithmeticoid $\by_0\in\sA$, then one knows one element of $\Psi_{C,\sA}$ namely $ \tTI_{C/\arith{L}_{\by_0}}$. 
\item In practice $H^1(\arith{L}_\bz,\Q(1))$ is a bit difficult to work with for computing volumes needed in this paper (and in \iut) and in fact one works with the Bloch-Kato subspace $H^1_e(\arith{L}_\bz,\Z(1))$ and suitably defined ``tensor product version'' of this cohomology which is identified with  Mochuziki's adelic tensor product  log-shell $\bsIm$ in \cite[{\ssep\ref{III-se:mochizuki-construction-thetam}}]{joshi-teich-rosetta}. Volumes of the set constructed above (and the theta-values locus $\thetami$ which is its variant) are computed by means of \cite[Proposition 1.4]{mochizuki-iut4} and this gets used in the proof of \Cref{th:theta-upper-bound}.
\item In some sense the idea of the proof of \Cref{th:main-dioph-thm} is to bound a specific value $\log(\TI_{C/\arith{L}_\bz})$ by a suitably defined average of function $\yadl\ni\bz\mapsto \log(\TI_{C/\arith{L}_\bz})$ constructed above. Ideally, in the context of \Cref{pr:collation-principle2}, one would like to talk about $$\int_{\yadl} \log(\TI_{C/\arith{L}_\bz}) d\bz$$
for a suitable measure $d\bz$ on $\yadl$. But presently, one is far from being able to do this. So the proof proceeds quite differently. 
\item More precisely, in the context of \cite{mochizuki-iut4}, one would like to talk about $$\int_{\Sigma_L'\subset \yadl^\ells} \left(\sum_{j=1}^\ells\log(\TI_{C/\arith{L}_{\bz_j}})\right) d\bz_1\cdot d\bz_2\cdots d\bz_\ells$$
for a suitable measure $d\bz_1\cdot d\bz_2\cdots d\bz_\ells$ on $\yadl^\ells$ and where $\Sigma_L'$ is Mochizuki's Adelic Ansatz constructed in \cite[{\ssep\ref{III-se:adelic-ansatz}}]{joshi-teich-rosetta}. Again one is far from being able to do this, but this is roughly how one can think of the volume computations of \cite{mochizuki-iut4}. 
\item One should note that similar sort of volume/integral computations occur in Classical Teichmuller Theory--for instance the survey \cite{wright19} of \cite{mirzakhani07}.
\item Mochizuki's log-Volume computations of \cref{ss:log-vol} (which enter the proof of \Cref{th:main-dioph-thm}) are similar to this in spirit.
\eenum
\erem

\subsection{The log-Different+log-Conductor Theorem}\nwss
A central point in the proof of the main theorem of \cite[Theorem 2.1]{mochizuki-general-pos} (\abc\ for compactly bounded subsets implies Vojta's Height Inequality for Curves) and of \cite[Theorem 1.10]{mochizuki-iut4} is a basic assertion \cite[Proposition 1.7]{mochizuki-general-pos} which I call \textit{``The log-Different+log-Conductor Theorem.''}

This theorem, stated as \cite[Proposition 1.7]{mochizuki-general-pos}, is stated in greater generality for arithmetic surfaces. But the proof is not adequately documented (in my opinion) and certainly uses facts about differents of fields which are not established in \cite{serre1979-local-fields}. Moreover, the proof of this point does not appear in \cite{fucheng}, \cite{yamashita} or \cite{dupuy2020probabilistic} (but its consequences are used in \cite[Theorem 2.1]{mochizuki-general-pos} and in \cite{mochizuki-iut4}). 

Also note that a version of \cite[Proposition 1.7]{mochizuki-general-pos} already appears as \cite[Proposition 14.4.6]{bombieri-gubler} where this is proved (for arithmetic surfaces) using theory of local height functions. 

Since I do not use the theory of local heights (and nor do \cite{mochizuki-general-pos} \& \cite{mochizuki-iut4}) I have felt that I should provide a direct proof  of this theorem. 

[I have pointed out in  \cite[{\ssep\ref{II5-se:heights}}]{joshi-teich-def}, that to understand how Arithmetic Teichmuller Spaces impacts the theory of heights in Diophantine Geometry, one must understand how many inequivalent local heights functions arise (at all primes) from many distinct local analytic geometries because of the existence of Arithmetic Teichmuller Spaces.]  

The theorem is  as follows. 

\bthm[The log-different+log-conductor theorem]\label{th:diff-cond-thm}
Assume $C,L$ are as in \inithtdata, \initassumptions. Let $M$ be a finite extension of $L$. Let $\ff_L,\ff_M$ and $\fd_{L},\fd_M$ be as defined in \cref{ss:diff-and-log-diff} and \cref{ss:tate-divisor} respectively.
Then 
\benumlab
\item one has the formula $$(\log\fd_M+\log\ff_M) -(\log\fd_L+\log\ff_L) =\frac{1}{[M:\Q]} \sum_{v\in\vlnon}\sum_{w|v}e_{w|v}\cdot\tau_w\cdot f_{w|v}\cdot\log v$$
where $\tau_{w|v}=0$ if and only if $w|v$ is at worst tamely ramified, and $1\leq \tau_w\leq {\rm ord}_v (e_{w|v})$ if $w|v$ is wildly ramified. 
\item Hence 
$$0\leq (\log\fd_M+\log\ff_M) -(\log\fd_L+\log\ff_L).$$
\item Moreover $$(\log\fd_M+\log\ff_M) -(\log\fd_L+\log\ff_L)$$ depends only on primes $v\in\V_{L}$ which are wildly ramified in $M/L$.
\item Especially $(\log\fd_M+\log\ff_M) = (\log\fd_L+\log\ff_L)$ if and only if $M/L$ is tamely ramified.
\item Let $S^\Q_{wild}$ be the set of all rational primes which lie below primes in $S_{wild}$. Then  one has the bound
$$(\log\fd_M+\log\ff_M) -(\log\fd_L+\log\ff_L)\leq \#S_{wild}^\Q\cdot \log[M:L].$$
\eenum
\ethm
\brem 
The reason this theorem is so named is as follows. By \Cref{re:diff-rmk}, the normalized arithmetic degree of the reduced Tate divisor  $\ff_L=(\fq_L)_{red}$  is the conductor of $C/L$ at its primes of semi-stable reduction (\cref{ss:tate-divisor}). While $\log\fd_L$ is the called  ``log-different'' in \cite{mochizuki-general-pos}.
\erem
\bp 
All the assertions are immediate from {\bf(1)}. So it suffices to prove {\bf(1)}.
By \cite[Proposition 14.3.10]{bombieri-gubler} one has
\be\label{eq:log-diff-eq}  
\log \fd_M-\log \fd_L=\frac{1}{[M:\Q]}\cdot\sum_{v\in\vlnon}\fd_{M_w/L_v}\cdot \log v.
\ee
Now by \cite[Dedekind's Discriminant Theorem, B.2.12]{bombieri-gubler} one has a formula for $\fd_{L_v}$ given by
\begin{align}\label{eq:dedekind-formula}
\fd_{M_w/L_v}&=\sum_{w|v} (e_{w|v}-1+\tau_{w|v})\cdot f_{w|v}.
\end{align}
where there are the following possibilities for $\tau_{w|v}\in\Z_{\geq0}$:
 $$\tau_{w|v}=\begin{cases}
 0 & \text{ if } w|v \text{ is tamely ramified},\\
[1,e_{w|v}\cdot {\rm ord}_v(e_{w|v})]\subset \N   & \text{ if } w|v \text{ is wildly ramified}.
\end{cases}$$
Substituting the formula \eqref{eq:dedekind-formula} into \eqref{eq:log-diff-eq}  one obtains:
\begin{align}\label{eq:log-diff-eq2}
\log \fd_M-\log \fd_L&=\frac{1}{[M:\Q]}\cdot\sum_{v\in\vlnon}\left(\sum_{w|v} (e_{w|v}-1)\cdot f_{w|v}\cdot \log v+\sum_{w|v}\tau_{w|v}\cdot f_{w|v}\cdot \log v\right),
\end{align}

Next let me compute $\log\ff_L-\log\ff_M$
\begin{align}\label{eq:log-cond-eq}
\begin{split}
\log\ff_L-\log\ff_M&=\frac{1}{[L:\Q]}\cdot\sum_{\vlnon} 1\cdot \log v-\frac{1}{[M:\Q]}\cdot\sum_{\vnon_{M}} 1\cdot \log w\\
&=\frac{[M:L]}{[M:\Q]}\cdot\sum_{\vlnon} 1\cdot \log v-\frac{1}{[M:\Q]}\cdot\sum_{\vnon_{M}} 1\cdot \log w\\
&=\frac{1}{[M:\Q]}\cdot\sum_{\vlnon} \left(\sum_{w|v} e_{w|v}\cdot f_{w|v} \right)\cdot 1\cdot \log v-\frac{1}{[M:\Q]}\cdot\sum_{\vnon_{M}} 1\cdot \log w,\\
&=\frac{1}{[M:\Q]}\cdot\sum_{\vlnon} \left(\sum_{w|v} e_{w|v}\cdot f_{w|v} \right)\cdot 1\cdot \log v-\frac{1}{[M:\Q]}\cdot\sum_{\vlnon} \sum_{w|v} 1\cdot \log w,\\
&=\frac{1}{[M:\Q]}\cdot\sum_{\vlnon} \left(\sum_{w|v} (e_{w|v} -1) \cdot  f_{w|v}\cdot \log v\right).
\end{split}
\end{align}
Hence subtracting \eqref{eq:log-cond-eq} from \eqref{eq:log-diff-eq2} one obtains
\be\label{eq:log-diff-eq3}
\log \fd_M-\log \fd_L-\left(\log \ff_L-\log \fd_M\right)=\frac{1}{[M:\Q]}\cdot
\sum_{v\in\vlnon}\left(\sum_{w|v}\tau_{w|v}\cdot f_{w|v}\cdot \log v\right).
\ee
Hence rearranging the left hand side one obtains
\be
\left(\log \fd_M+\log \ff_M\right)-\left(\log \fd_L+\log \ff_L\right)=\frac{1}{[M:\Q]}\cdot\sum_{v\in\vlnon}\left(\sum_{w|v}\tau_{w|v}\cdot f_{w|v}\cdot \log v\right).
\ee
This proves the first assertion of \Cref{th:diff-cond-thm}. Since the terms on the right hand side in {\bf(1)} are all non-negative,  this also proves the second assertion.
As $\tau_{w|v}=0$ unless $w|v$ is wildly ramified, the remaining assertions are immediate from {\bf(1)}. This completes the proof of 
\Cref{th:diff-cond-thm}.

Now let me prove the estimate given in {\bf(5)}.   By  {\bf(1)} one has to bound
\be
\frac{1}{[M:\Q]} \sum_{v\in\vlnon}\sum_{w|v}\tau_{w|v}\cdot f_{w|v}\cdot\log v
\ee
Since $\tau_{w|v}\leq e_{w|v}\cdot {\rm ord}_v(e_{w|v})$ one obtains
\be
\frac{1}{[M:\Q]} \sum_{v\in\vlnon}\sum_{w|v}\tau_{w|v}\cdot f_{w|v}\cdot\log v\leq \frac{1}{[M:\Q]} \sum_{v\in\vlnon}\sum_{w|v}e_{w|v}\cdot {\rm ord}_v(e_{w|v})\cdot f_{w|v}\cdot\log v.
\ee
Firstly, both the sum is over all $v\in S_{wild}$.  Secondly, one has $\sum_{w|v}e_{w|v}\cdot f_{w|v}=[M:L]$, and thirdly one has
 \be {\rm ord}_v(e_{w|v})=e_{v|p}\cdot {\rm ord}_p(e_{w|v}).\ee 
  Now $e_{w|v}\leq [M_w:L_v]$ 
  implies that \be \tau_{w|v}\leq e_{w|v}\cdot e_{v|p}\cdot \frac{\log([M_w:L_v])}{\log(p)}\leq e_{w|v}\cdot e_{v|p}\cdot \frac{\log([M:L])}{\log(p)}.\ee 
 One also has 
 \begin{align}
 \begin{split}
 \log v&=f_{v|p}\cdot\log(p),\\ 
 [M:\Q]&=[M:L]\cdot [L:\Q],\\
[L:\Q]&=\sum_{v|p}e_{v|p}\cdot f_{v|p}.\\ 
 \end{split}
 \end{align}
 So one sees that 
\begin{align}\label{eq:wild-term-ex}
\begin{split}
\frac{1}{[M:\Q]} \sum_{v\in\vlnon}\sum_{w|v}e_{w|v}\cdot\tau_{w|v}\cdot f_{w|v}\cdot\log v&\leq
\frac{1}{[M:\Q]} \sum_{v\in S_{wild}}\sum_{w|v}e_{v|p}\cdot f_{v|p}\cdot e_{w|v}\cdot f_{w|v}\cdot\frac{\log([M:L])}{\log(p)}\cdot\log(p)\\
&\leq
\frac{1}{[L:\Q]} \sum_{p\in S^\Q_{wild}}\sum_{v|p} e_{v|p}\cdot f_{v|p}\cdot{\log([M:L])}\\
&\leq \#S^\Q_{wild}\cdot\log([M:L]).
\end{split}
\end{align}
This completes the proof of the theorem.
\ep

An important corollary of the log-Different+log-conductor Theorem (\Cref{th:diff-cond-thm}) is the following:
\bcor 
In addition to the notations and assumptions of \Cref{th:diff-cond-thm}, assume that one is given  a positive constant $0<d\in\R$, and a subset $S\subset \vlnon$ and a subset $S_{wild}\subset S$ and a family of finite extensions $$\sF_L=\left\{M/L: M\subset \bQ \right\}$$ of $L$ such that 
\benumlab
\item  $[M:L]\leq d$
\item  each $M\in\sF_L$ is  unramified outside $S$, 
\item each $M\in\sF_L$ is tamely ramified outside $S_{wild}\subset S$. 
\eenum 
Then there exists a constant $0<c\in\R$ depending on $C,[M:L],S_{wild}$ such that
$$0\leq \left(\log \fd_M+\log \ff_M\right)-\left(\log \fd_L+\log \ff_L\right)\leq c.$$
\ecor

\brem 
The arithmetic surface case in \cite[Proposition 14.4.6]{bombieri-gubler} and \cite[Proposition 1.7(i)]{mochizuki-general-pos} requires additional arguments to handle vertical ramification but the general principle of the proof is the same as the one given above.
\erem
\bp 
The right hand side of the equation in \Cref{th:diff-cond-thm}, depends only on the set of primes of wild ramification and these are all contained in $S_{wild}$. Since $[M:L]\leq d$ so $[M_w:L_v]\leq C$ for each $M/L$ and each prime $w\in S_{wild}$ of wild ramification. The local extensions all have bounded degrees bounded by $[M:L]\leq d$. Hence there are only a finitely many local extensions of bounded degrees of $L_v$ for each $v\in S$ and as one is dealing with a fixed finitely many primes of wild ramification, so the numbers $\tau_{w|v}$ appearing are absolutely bounded. Hence the right hand side in \Cref{th:diff-cond-thm} is absolutely bounded with the bound $0<c$ depending only on $L,d,[M:L], S_{wild}$. This proves the assertion. 
\ep

\newcommand{\exc}{\mathcal{Exc}}

\section{The Existence of Initial Theta-Data}\label{se:existence}
\subsection{Bounding Domains and Compactly Bounded Subsets}\label{ss:compact-bnded}
I recall a few facts about compactly bounded subsets defined in \cite{mochizuki-general-pos}. My approach to it is different (but equivalent) to the one given in \cite{mochizuki-general-pos}.

Let $L$ be a number field with no real embeddings. Let $X/\arith{L}_\by$, with $\by=(y_v)_{v\in\vl}$, be an arithmeticoid of $L$ given by $\by\in\yadl'$. This means that one is given an algebraic closure $\bL$ of $L$ and for each $v\in\vl$, let $\bQ_{p_v}\subset K_{y_v}$ be the algebraic closure of $\Q_{p_v}$ contained, as a valued field, in the algebraically closed perfectoid field $K_{y_v}$ provided by $\arith{L}_{\by}$ at $v$. Note by our assumption on $L$  for each $v\in\vlarch$, one has $\bQ_{p_v}=\C$.

Let $\vlarch\subset S\subset \vl$ be a finite set of primes of $L$ (so $S\neq \emptyset$). Following the convention in the preceding section, for $v\in \vlnon$ write $p_v$ for the residue characteristic of $v$.   Then consider the product topological space (with $\bQ_{p_v}\subset \hat{\bQ}_{p_v}\subset K_{y_v}$ as above and where $\hat{\bQ}_{p_v}$ is the completed algebraic closure of $\Q_{p_v}$ in $K_{y_v}$).
\be  
\xs =\prod_{v\in S}X(\bQ_{p_v}) \subset \prod_{v\in S}X(\hat{\bQ}_{p_v})\subset \prod_{v\in S}X(K_{y_v}).
\ee
Then $\xs$ is equipped with a natural Galois action of $$G_{L,S}=\prod_{v\in S}\gal(\bQ_{p_v}/L_v)\act \prod_{v\in S}X(\hat{\bQ}_{p_v}).$$ 
A \textit{bounding domain} in $\xs$ is a $G_S$-stable,  compact subset $\emptyset\neq Z=\prod_{v\in S} Z_v \subset \xs $ such that $Z$ is the closure of its interior. The finite set of primes $S$ will be called \textit{the support of the bounding domain $Z$}. 

Let $Z$ be a bounding domain with support $S$. Then $(Z,S)$ determine a subset of $\sZ\subset X(\bL)$ as follows.
\begin{align}
\begin{split}
\sZ &=\left\{ \vphantom{\prod_{v\in S}\bQ_{p_v}} x : \text{ if } x\in X(M) \text{ for some finite extension } L\subseteq M\subset\bL \right.  \\ 
& \qquad \left. \text{ then for each embedding } \iota:M\to \prod_{v\in S}\bQ_{p_v}, \right. \\  & 
\left. \qquad \text{ one has } \iota(x)\in Z\subset \xs
\vphantom{\prod_{v\in S}\bQ_{p_v}} \right\} \subset X(\bQ)
\end{split}
\end{align}
The set $\sZ\subset X(\bQ)$ will be called the \textit{$S$-supported, compactly bounded domain} determined by the bounding domain $Z$. The set $\sZ$  uniquely determines the pair $(Z,S)$ \cite[Page 6]{mochizuki-general-pos}. If I want to emphasize $(Z,S)$ then I will write $\sZ(Z,S)$ in place of $\sZ$.

\subsection{A worked example: $\P^1-\{ 0, 1, \infty\}$} The following example will help understanding the idea of bounding domains and compactly bounded domains. 
Let $X=\P^1-\{ 0, 1, \infty\}/\Q$. Let $S=\{\infty,p\}=\varc_\Q\cup\{ p\}$. Then $\bQ_\infty=\C$  and let $R_0,R_1,R_\infty\subset X(\C)=\P^1-\{ 0, 1, \infty\}(\C)$ (resp. $T_0,T_1,T_\infty\subset X(\bQ_p)= \P^1-\{ 0, 1, \infty\}(\bQ_p)$) be mutually disjoint galois invariant compact subsets not containing the points $0,1,\infty$ respectively and contained in $X(\C)$  (resp. $X(\bQ_p)$ respectively). [One can take $R_0,R_1,R_\infty$ to be closed, mutually disjoint annuli contained in $X(\C)\subset \P^1(\C)$  around the points $0,1,\infty$ but not containing any of these points. Since $\bQ_p$ is not locally compact, a compact Galois stable subset in it is more complicated to describe. Here is the construction of a  typical example of a compact, Galois stable subset of $\bQ_p$ (or $\C_p$): let $t\in\C_p-\bQ_p$ so $t$ is transcendental over $\bQ_p$. Then  by a well-known theorem of Tate, $t$ is not fixed by the Galois action on $\C_p$. The Galois orbit of $t$ in $\C_p$ is compact (by the continuity of the Galois action on $\C_p$ and by the compactness of the Galois group), and the intersection of this compact, Galois stable orbit with $\bQ_p$ is compact and Galois stable. On the other hand, as $\bQ_p$ (resp. $\C_p$) is not locally compact, a subset of these fields defined by finitely many inequalities is bounded but not necessarily compact in general).] Then $$\xs = \left((\P^1-\{ 0, 1, \infty\})(\C)\right)\times \left((\P^1-\{ 0, 1, \infty\})(\bQ_p)\right)$$ as a product topological space.   Let $Z_\infty=(R_0\cup R_1\cup R_\infty)$ and $Z_p=(T_0\cup T_1\cup T_\infty)$. Then $$Z=Z_\infty\times Z_p\subset \xs=\left((\P^1-\{ 0, 1, \infty\})(\C)\right)\times \left((\P^1-\{ 0, 1, \infty\})(\bQ_p)\right).$$ Then $Z$ is a bounding domain with support $S$. The compactly bounded subset $\sZ=\sZ(Z,S)$ determined by $(Z,S)$ consists of all $x\in X(\bQ)$ such that $x\in X(M)$ for some finite extension $\Q\subset M\subset \bQ$ such that under any embedding $\iota:M\into \C\times \bQ_p$, one has $\iota(x)\in Z$. This means  $\iota_\infty(x)\in R_0\cup R_1\cup R_\infty$ and $\iota_p(x)\in T_0\cup T_1\cup T_\infty$.

\subsection{Legendre Elliptic Curves}
Let $L$ be a field and let $\lambda\in L-\{0,1\}$. Then one has the associated elliptic curve $C_\lambda: y^2=x(x-1)(x-\lambda)$. Let $j_\lambda=j(C_\lambda)$ be the $j$-invariant of $C_\lambda/L$. Then by the standard formulae \cite{silverman-arithmetic}, one knows \cite[Chap III, Proposition 1.7]{silverman-arithmetic} that $j_\lambda$ is a rational function of $\lambda$ and by standard moduli theory, the field of definition $\lmod$ of the moduli point $[C_\lambda]\in \sM_{ell}$ is $\lmod=\Q(j_\lambda)$.

\subsection{Two Variants of the Tate Divisor Function}
Let $\arith{L}_\by$ be an arithmeticoid of a number field $L$. Let $C/\arith{L}_\by$ be a holomorphoid of an elliptic curve $C/L$ and let $X/\arith{L}_\by$ be the holomorphoid of $X=C-\{O\}$. One can view the isomorphism class of $[C/L]$ as providing a point of the moduli $\sM_{ell}$ of elliptic curves.  

The following definition is a variant of the definition given in \cref{ss:tate-divisor} (where the discussion was restricted to primes $\ubblvoss_{\lmod}$--but the definition below uses all primes of $L$ at which one has semi-stable reduction):
\begin{defn}
	Let $C/M\in \sM_{ell}(M)$ be an elliptic curve over a number field $M$.	Let $w\in\vnon_M$ and let  $q_w$ at $C_{M_w}$ be a Tate parameter (well defined up to multiplication by a local unit and its order at $w$). Consider the function 
	$$
	\fT(C/M)= \frac{1}{[M:\Q]}\cdot\sum_{w\in \vnon_L} \ord_w(q_w)\cdot \log w.
	$$
	Then the function $$\fT:\sM_{ell}(\bQ)\to \R$$
	given by $C/M\mapsto \fT(C/M)$ with the elliptic curve $C$ being viewed over its field of moduli, say $M$,  is a well defined $\R$-valued function on $\sM_{ell}(\bQ)$.
\end{defn}

Let $\Psi\subset \vnon_{\Q}$ be a fixed finite set of rational primes. The following variant, denoted $\fT_\Psi$, of $\fT$,  will also be useful. Let 
\be
\fT_\Psi(C/M)=\frac{1}{[M:\Q]}\cdot\deg\left(\sum_{w\nmid p\in\Psi}\ord_w(q_w)\right)\in\R
\ee be the logarithmic degree of the portion of $\fT({C/M})$ which is not supported on any of the primes of $M$ lying over the rational prime $p\in\Psi$. In \cite{mochizuki-iut4}, the choice of $\Psi=\{2 \}$ is used and I will write $\fT_2:=\fT_{\{2 \}}$ for this choice of $\Psi$.

\brem 
The following table provides a translation to the notation used in \cite{mochizuki-iut4}:
\begin{center}
	\begin{tabular}{|c|c|c|}
		\hline 
		\cite{mochizuki-iut4}	& $\fq^{\forall}_{(-)}$ & $\fq^{\nmid2}_{(-)}$ \\ 
		\hline 
		This paper	& $\fT(-)$ &  $\fT_2:=\fT_{\{2 \}}(-)$\\ 
		\hline 
	\end{tabular} 
\end{center}
Since some of proofs given in \cite{mochizuki-iut4} use the usual Tate divisor $\fq_{C/M}$ of a fixed elliptic curve (defined in \cref{ss:tate-divisor}), and especially $\fq$ is the shorthand for this Tate divisor and $\log(\fq)$ is its normalized degree, I found Mochizuki's choice of notation $\fq^\forall$ and $\fq^{\nmid2}$ visually and algebraically confusing. 
\erem

\subsection{More notations}\label{ss:notate-setup-exist-thm}
Let $X=\P^1_\Q$, $\O_{\P^1}(1)$ be the standard ample line bundle on $\P^1$. Let $D=[0]+[1]+[\infty]$ be the reduced divisor supported on $\{0,1,\infty \}\subset \P^1_{\Q}$. Since $\omega_{\P^1}=\O_{\P^1}(-2)$, so $$\omega_{\P^1}(D)=\O_{\P^1}(1).$$ Let $h_{\omega_{\P^1}(D)}=h_{\O_{\P^1}(1)}$ be the logarithmic height function on $\P^1$ given by the ample divisor $\omega_{\P^1}(D)$. 

Let $U=\P^1-\{0,1,\infty \}=X-D$ as parameterizing the Legendre family of elliptic curves
\be\label{eq:legendre} U\ni\lambda\mapsto C_\lambda: y^2=x(x-1)(x-\lambda).\ee
This provides the classifying morphism 
$$U\to \sM_{ell},$$
given by $\lambda\mapsto [C_\lambda]$.
Let $$j_{inv}:U\to \sM_{ell}\to \A^1$$
be the $j$-invariant morphism $\lambda\mapsto j(C_\lambda)$. Let $\sM_{ell}\to \sMb_{ell}$ be the usual compactification of $\sM_{ell}$. Then one can pull-back the divisor at infinity of $\sMb_{ell}$ to $U$ by composition with the above morphisms. Let
\be 
h_\infty: U(\bQ)\to \R 
\ee
be the height function on $U$ defined by this pull-back of the divisor at infinity of $\sMb_{ell}$ (see \cite{mochizuki-general-pos} for more details).
\subsection{Height comparisons}
\bpro\label{pr:height-comparisons} 
Suppose $U$ is as above and $\{ 2, \infty\}\subseteq S\subset \V_{\Q}$ is a finite subset of primes.  Suppose $\sZ(Z,S)\subset \sM_{ell}(\bQ)$ is an $S$-supported, compactly bounded subset. Assume that 
\be\label{eq:2-adic-assumption}
Z_2\subset 2^{N_{j_{inv}}}\cdot \O_{\bQ_2}
\ee
for some integer $N_{j_{inv}}\in \Z$ (this implies that if $C_\lambda\in\sZ(Z,S)$ then the image of $C_\lambda$ under the $j$-invariant mapping is contained in $Z_2\subset 2^{N_{j_{inv}}}\cdot \O_{\bQ_2}$ for every embedding of $\Q(\lambda)\into \bQ_2$). Then one has the following equality of bounded discrepancy classes of functions:
$$\frac{1}{6}\cdot \fT_2\approx \frac{1}{6}\cdot \fT \approx \frac{1}{6}\cdot h_\infty \approx h_{\omega_{\P^1}(D)}.$$
\epro
\bp This is standard \cite{mochizuki-general-pos}. The assertion $\frac{1}{6}\cdot h_\infty \approx h_{\omega_{\P^1}(D)}$ is quite elementary: $\omega_{\P^1}(D)$ is a divisor of degree one on $U$, but by the pull-back of the divisor at $\infty$ in $\sMb_{ell}$ to $U$ is of degree $6$. The assertion $\frac{1}{6}\cdot \fT_2\approx \frac{1}{6}\cdot \fT$ or equivalently $\fT_2\approx \fT$ follows from the definition of $\fT$ and $\fT_2$ and \eqref{eq:2-adic-assumption} holds. The rest of the assertions are clear.
\ep

\newcommand{\Rpos}{\R_{>0}}
\newcommand{\crit}{\mathscr{C}}
\newcommand{\epd}{\varepsilon_d}
\newcommand{\hunif}{H_{unif}}

\newcommand{\thmreminder}{With above notations and under the hypothesis of \Cref{th:existence}.}

\subsection{The Existence Theorem}\nwss
In this section the goal is to demonstrate the existence of an elliptic curve $C/\arith{L}_{\by_0}$ over a standard arithmeticoid $\arith{L}_{\by_0}$ of $L$ which is equipped with Initial Theta Data \inithtdata, \assumptions. The following theorem is my formulation of \cite[Corollary 2.2]{mochizuki-iut4}. This theorem is fundamental for the construction of the theta-values locii of \cite{joshi-teich-rosetta} (and hence even for the very statement of \cite[{Corollary \ref{III-cor:cor312}}]{joshi-teich-rosetta}, \moccor). 

\noindent\textcolor{red}{Note}: I have broken up the original assertion of \cite[Corollary 2.2]{mochizuki-iut4} in two parts: what appears below corresponds to  \cite[Corollary 2.2(i,ii)]{mochizuki-iut4} while \cite[Corollary 2.2(iii)]{mochizuki-iut4} is subsumed in the proof of \Cref{th:main-dioph-thm}. This is better in my opinion.

\bthm\label{th:existence}
The notations, conventions of \cref{ss:notate-setup-exist-thm} and assumptions of \Cref{pr:height-comparisons}  hold. This in particular means that one is given an $S$-supported compactly bounded subset $\sZ=\sZ(Z,S)\subset \sM_{ell}(\bQ)$ with a bounding domain $Z\subset \xs$ and the assumption \eqref{eq:2-adic-assumption} holds. Let $d\in\N$, $\delta=\dmods=2^{12}\cdot 3^3\cdot 5\cdot\dmod$. 
Then there exists a finite set $\exc=\exc(\sZ,d)\subset U(\bQ)^{\leq d}$, depending only on $\sZ,d$, with the following properties:

Suppose $x_\lambda\in\sZ\cap U(\bQ)^{\leq d}$ and $x_\lambda\not\in\exc$ and $C_\lambda$ be the Legendre elliptic curve \eqref{eq:legendre} corresponding to $x_\lambda$. Let $\lmod$ be the minimal field of definition of $C_\lambda$, let $L=\lmod(\sqrt{-1},C_\lambda[2\cdot3\cdot5])$ and let $\ltpd=\lmod(C_\lambda[2])$. Let $Q=\fT(C_\lambda/L) \text{ and } Q_2=\fT_2(C_\lambda/L)$. Let $\fq=\fq_{C_\lambda/L}$ be the Tate divisor (defined in \cref{ss:tate-divisor}) of $C_\lambda$. Then there exists a prime $\ell$ satisfying
$$Q^{1/2}\leq \ell \leq 10\cdot\delta\cdot Q^{1/2}\cdot \log(2\cdot\delta\cdot\log(Q))$$ 
such that $C_\lambda$ satisfies Initial-Theta Data \inithtdata\ (also \cite[\ssep 3.1]{mochizuki-iut1}) and \assumptions, with the prime $\ell$ satisfying the above bound.
\ethm
\newcommand{\xip}{\xi_{prm}}
\subsection{Proof of the Existence Theorem (\Cref{th:existence})}\nwss
Let me say that \Cref{th:existence} is a global assertion over a number field which cannot be proved by local means and \Cref{th:existence} is central not just to this paper but also to the main results of \cite{joshi-teich-rosetta}. Logically since the \Cref{th:existence} is needed for the proof of the fundamental estimate \moccor, the theorem should appear in \cite{joshi-teich-rosetta}. However for ease of comparison between my theory and Mochizuki's theory, a result which appears in \cite{mochizuki-iut3} appears in \cite{joshi-teich-rosetta} and similarly with \cite{mochizuki-iut4} and \cite{joshi-teich-rosetta}.

In what follows,  write $$Q=\fT(C_\lambda)=\sum_{v\in \vl} a_v \cdot \log(v)=\sum_{v\in \vl} a_v \cdot f_v\cdot \log(p_v) $$
with $a_v\in\Z_{\geq0}$.

Start with $\exc=\emptyset$. This will be enlarged  step by step as the proof progresses.

Firstly, note that since one is working with the Legendre family, there exists a finite number of points at which one has stacky behavior arising from curves in this family  with automorphism groups of order $>2$ (over $\bQ$). Add to $\exc$ the elliptic curves corresponding to these four $j$-invariants if the points corresponding to these curves are contained in $\sZ$. 

Next I claim that there are finitely many elliptic curves for which $Q^{1/2}\leq \xip$. This says, by \Cref{pr:height-comparisons}, that the height of $C_\lambda$ is bounded and by \cite[Theorem 2.1]{silverman1986}, there are only a finitely many points in $U(\bQ)^{\leq d}$ of bounded degree and of bounded height. Enlarge $\exc$ by adding these curves to it. 

Thus for any point of $\sZ\cap U(\bQ)^{\leq d}$ not in $\exc$, one has $$Q^{1/2}\geq \xip\geq 5.$$

\blem\label{le:prm-num-est}\ 
\thmreminder
\benumlab
\item There exists an $5\leq\xip\in\R$ such that for all $Q\geq \xip$ one has
$$\frac{2}{3}\cdot Q\leq \theta(Q)=\sum_{p\leq{Q}} \log(p)\leq \frac{4}{3}\cdot Q.$$
\item If $\sA$ is a finite set of prime numbers, let 
$$\theta_\sA=\sum_{p\in\sA}\log(p).$$
Then there exists a prime $p\not\in\sA$ such that
$$p\leq 2(\xip+\theta_\sA).$$
\eenum
\elem
\bp 
This is proved in \cite[Lemma 4.1]{mochizuki-general-pos} and the reader is referred to these sources for its proof. 
\ep

\blem\label{le:prm-bnd}\  
\thmreminder\ 
Let $\sA$ be the finite set of primes $p$ such that either
\benumlab
\item $p\leq Q^{1/2}$, or
\item $p|a_v$ for some $v\in\vlnon$, or
\item $p=p_v$ for some $v\in\vlnon$ for which $a_v\geq Q^{1/2}$.
\eenum
Then $$\theta_\sA=\sum_{p\in\sA}\log(p)\leq -\xip+5\cdot\delta\cdot Q^{1/2}\cdot\log(2\cdot\delta\cdot\log(Q)).$$
\elem

\bp 
Write 
$$\theta_\sA=\sum_{p\in\sA,p\leq Q^{1/2}}\log(p)+\sum_{p\in\sA,p|a_v}\log(p)+\sum_{p=p_v, a_v\geq Q^{1/2}}\log(p).$$
and estimate each sum separately. The first sum is estimated by \Cref{le:prm-bnd} to be $\leq\frac{4}{3}\cdot Q^{1/2}$. 

The second sum $\sum_{p\in\sA,p|a_v}\log(p)$ is estimated as follows. \begin{align}
\begin{split}
2\cdot Q^{1/2}\cdot\delta\cdot\log(2\delta\cdot Q)&\geq 2[L:\Q]\cdot Q^{1/2}\cdot\log(2\cdot[L:\Q]\cdot Q)\\
&\geq 2\cdot \frac{Q}{Q^{-1/2}}\cdot \log(2\cdot[L:\Q]\cdot Q)\\
&\geq 2\cdot Q^{-1/2}\sum_{a_v\neq 0}a_v\cdot f_v\cdot \log(p_v)\cdot \log(2\cdot[L:\Q]\cdot Q)\\
&\geq 2\cdot \sum_{a_v\geq Q^{1/2}}(a_v\cdot Q^{-1/2})\cdot \log(p_v)\\
&\geq 2\sum_{a_v\geq Q^{1/2}}\log(p_v).
\end{split}
\end{align}
Hence 
\be 
\sum_{a_v\geq Q^{1/2}}\log(p_v)\leq \delta\cdot Q^{1/2}\cdot\log(2\delta\cdot Q)
\ee
I claim that the third sum is
$$\sum_{p=p_v, a_v\geq Q^{1/2}}\log(p)\leq \delta\cdot Q^{1/2}.$$ This is estimated as follows.
Since $[L:\Q]\leq\delta$, one sees that
\begin{align}
\begin{split}
\delta\cdot Q^{1/2} &\geq [L:\Q]\cdot Q^{1/2}\\
&\geq [L:\Q]\cdot Q\cdot Q^{-1/2},\\
&=Q^{-1/2}\cdot \sum_v a_v\cdot f_v\cdot\log(p_v)\\
&\geq  \sum_v \frac{a_v}{Q^{1/2}}\cdot\log(p_v)\\
&\geq  \sum_{a_v\geq Q^{1/2}} \frac{a_v}{Q^{1/2}}\cdot \log(p_v)\\
&\geq  \sum_{a_v\geq Q^{1/2}} \log(p_v),\\
&=\sum_{p=p_v, a_v\geq Q^{1/2}}\log(p_v).
\end{split}
\end{align}
Hence 
\be 
\sum_{p=p_v, a_v\geq Q^{1/2}}\log(p)\leq \delta\cdot Q^{1/2}.
\ee
This completes the proof of \Cref{le:prm-bnd}.
\ep

\blem\label{le:prm-bnd2}
\thmreminder\ Assume that $Q^{1/2}\geq 5$, the following assertions hold:
\benumlab
\item There exists a rational prime $\ell$ satisfying
$$
Q^{1/2}\leq \ell \leq 10\cdot\delta\cdot Q^{1/2}\cdot\log(2\cdot\delta\cdot Q).
$$
\item $\ell\nmid a_v$ for any $a_v\neq0$ and $v\in\vlnon$,
\item if $\ell=p_v$ for some $v\in\vlnon$, then $a_v<Q^{1/2}$.
\eenum
\elem
\bp 
This is immediate from \Cref{le:prm-bnd} and \Cref{le:prm-num-est}.
\ep

\blem\label{le:isogeny-lem}
\thmreminder\ Let $\sZ\subset U(\bQ)^{\leq d}$. Let $\ell$ be a prime given by \Cref{le:prm-bnd2}. Then there are only finitely many elliptic curves $C_\lambda\in \sZ$ with a subgroup of $C_\lambda[\ell](L)$ of order $\ell$.
\elem
\bp 
This is immediate from \cite[Lemma 3.5]{mochizuki-general-pos}.
\ep

Now enlarge $\exc$ by adding to it the exceptional points given by \Cref{le:isogeny-lem}.

\blem\label{le:exits-ss-prime}
\thmreminder\ 
Suppose  $C_\lambda\in \sZ\cap U(\bQ)^{\leq d}$ (with $\sZ$ satisfying \eqref{eq:2-adic-assumption}) is an elliptic curve such that no prime outside $\{v\in\V_{\lmod}: v\nmid(2\ell)\}$ is a prime of multiplicative reduction of  $C_\lambda$, then $C_\lambda$ is of bounded height.
\elem
\bp
By definition, $Q=\fT(C_\lambda)$. Hence the inequality
\be 
Q=\fT(C_\lambda)\approx \fT_2 (C_\lambda) \leq Q^{1/2}\cdot\log(\ell).
\ee
holds on $\sZ\cap U(\bQ)^{\leq d}$ as  $\sZ$ satisfies \eqref{eq:2-adic-assumption}. By \Cref{le:prm-bnd2}, one has  $$\ell \leq A\cdot Q^{3/4}$$ for a suitable constant $A>0$. So one sees that $Q^{1/2}$ is bounded and hence $Q$ is bounded. So the assertion follows.
\ep

Now enlarge the finite set $\exc$  by adding the curves $C_\lambda$ whose existence is asserted by \Cref{le:exits-ss-prime}.

\newcommand{\brho}{\overline{\rho}}
\blem\label{le:large-image} 
\thmreminder\ 
Let $C_\lambda\in \sZ\subset U(\bQ)^{\leq d}$  and $C_\lambda\not\in\exc$ be an elliptic curve (so $C_\lambda$ has least one prime of multiplicative reduction not lying over $2\cdot \ell$). Let
$$\brho_\ell:G_L\to GL_2(\F_\ell)$$
be the Galois representation on $\ell$-torsion of $C_\lambda$. Then
$$\brho_\ell(G_L)\supset SL_2(\F_\ell).$$ 
\elem
\bp 
This is immediate from \Cref{le:exits-ss-prime} and standard results. For a complete proof see \cite[Theorem 3.8(b)]{mochizuki-general-pos}.
\ep

\bp[Completion of the Proof of \Cref{th:existence}]
The existence of the prime number $\ell$ as asserted in \Cref{th:existence} is given by \Cref{le:prm-bnd2}. The remaining portion of \Cref{th:existence}  may now be completed using \Cref{le:isogeny-lem}, \Cref{le:exits-ss-prime} and \Cref{le:large-image} as follows. Choose any curve $C_\lambda$ in $\sZ$ not in $\exc$, choose the prime $\ell$ to be the prime given by \Cref{le:prm-bnd2}. Then by \Cref{le:exits-ss-prime}, $C_\lambda$ has at least one prime of multiplicative reduction so $\ubblvoss\neq\emptyset$. By \Cref{le:isogeny-lem} and \Cref{le:large-image} one obtains that for $C_\lambda$, the representation $\brho_\ell$ (of \Cref{le:large-image}) has full image. All these results and properties imply that the $C_\lambda$ is equipped with an Initial Theta Data satisfying \inithtdata\ and \initassumptions\ (and hence also the same properties in \cite[Section 3.1]{mochizuki-iut1} and \cite{mochizuki-iut4}). 
\ep

\section{The Main Bounds}\label{se:main-bounds}\nwss
\subsection{The First Main Bound}\label{se:first-main-bnd}
\Cref{th:existence} means that for $\lambda\in \sZ-\exc$,  the curve $(C_\lambda,L,X_\lambda,\ell,\V,\ubblvoss)$ one has the existence of Initial Theta Data satisfying \initassumptions, \assumptions. The following theorem is a cleaner formulation of \cite[Theorem 1.10]{mochizuki-iut4}. The proof is also reorganized for readability:
\newcommand{\bq}{{\bf q}}
\newcommand{\thmremindermnbnd}{Assume that one is given an elliptic curve $C=C_\lambda$ equipped with  the  prime number $\ell$ and \initassumptions, \assumptions\  provided by \Cref{th:existence}. Assume also that the notational conventions of \cref{se:existence} remain in force.}
\bthm\label{th:main-bound}
\thmremindermnbnd\  Let $\fq=\fq_{C/L}$ be the Tate divisor of $C/L$ \cref{ss:tate-divisor}. Then one has
\begin{align*}
\begin{split}
\frac{1}{6}\cdot\log(\fq) & \leq \left(1+\frac{20\cdot\dmod}{\ell}\right)\cdot\left(\log(\dtpd)+\log(\ftpd)\right)+20\left( \emods\cdot\ell+60\right) \\
& \leq \left(1+\frac{20\cdot\dmod}{\ell}\right)\cdot\left(\log(\fd_L)+\log(\ff_L)\right)+20\left( \emods\cdot\ell+60\right).
\end{split}
\end{align*}
\ethm
\brem This is the same as \cite[Theorem 1.10]{mochizuki-iut4} with $\etap=60$ (see \cref{sss:prime-bound}). As will be demonstrated below, one of the important ingredients in the proof is \moccor\ whose proof is completed in \cite[{\ssep\ref{III-cor:cor312}}]{joshi-teich-rosetta}.
\erem 
The proof is broken up into convenient steps. In what follows the assumptions of \Cref{th:main-bound} will be assumed to hold and used in various propositions that follow.
\subsection{Prime number  bounds}\label{sss:prime-bound}\nwss
One has the following explicit estimate:
\bpro For $x\geq \etap:=60$, and let $\pi(x)$ be the number of primes $\leq x$, then one has
$$\pi(x)\leq \frac{4}{3}\cdot\frac{x}{\log(x)}.$$
\epro
\bp 
This is immediate from the explicit estimates for $\pi(x)$ in \cite[Theorem 6.9]{dusart2016} which provides the first inequality given below for $x>1$, while the second inequality below is easily established for $x\geq 60>1$
$$\frac{\pi(x)}{x/\log(x)}\leq 1+\frac{1.2762}{\log(x)}\leq \frac{4}{3}.$$
\ep

\subsection{Easy lemmata}\label{ss:lemmata}\nwss The proofs of the following Lemmata are well-known or trivial.
\blem\label{le:sums}\ 
\benumlab
\item For any integer $n\geq1$, one has
$$\sum_{j=1}^{n}j=\frac{n\cdot (n+1)}{2}.$$
\item For any integer $n\geq1$, one has
$$\sum_{j=1}^{n}j^2=\frac{n\cdot (n+1)(2n+1)}{6}.$$
\eenum
\elem 
\blem\label{le:orders}
For any prime number $p$, one has
$$\abs{\glt{\F_p}}=p\cdot(p+1)\cdot(p-1)^2.$$
Especially
\benumlab
\item $\abs{\glt{\F_2}}=2\cdot3$,
\item $\abs{\glt{\F_3}}=2^4\cdot3$,
\item $\abs{\glt{\F_5}}=2^5\cdot3\cdot5$.
\eenum
\elem
\blem\label{le:deg} Evidently 
$$[\lmod(\sqrt{-1}):\lmod]\leq2.$$
\elem
\blem\label{le:trivial}\ 
\benumlab
\item $0 <\log(2)\leq1$,
\item $1\leq \log(3)<\log(\pi)\leq\log(5)\leq2$.
\eenum
\elem

\newcommand{\mycancelto}[2]{\cancelto{\textcolor{blue}{#1}}{\textcolor{red}{#2}}}
\subsection{Global Discriminant bounds}\nwss
The following lemmas are consequences of the log-Different+log-Conductor Theorem (\Cref{th:diff-cond-thm}). \textit{My estimates  for constants arising from the wild ramification term in \Cref{th:diff-cond-thm}, indicated in \textcolor{blue}{blue},  is different from Mochizuki's and is indicated in \textcolor{red}{red}}. [Difference arises from the factor $\#S^{\Q}_{wild}$ appearing in \eqref{eq:wild-term-ex}.]
\blem\label{le:basic-ineq1} 
Under the notations and assumptions of \Cref{th:main-bound},   one has
$$ \log(\fd_L)+\log(\ff_L)\leq \log(\fd_{\ltpd})+\log(\ff_{\ltpd})+\mycancelto{33}{21}.$$
\elem
\bp 
By \Cref{le:tame-lemma} one knows that $L/\ltpd$ is tamely ramified outside $v|30$. Hence the inequality is a consequence of log-Different+log-Conductor Theorem (\Cref{th:diff-cond-thm}) except for computation of the bound arising from the wild ramification term in \Cref{th:diff-cond-thm}{\bf(5)} and this is straight forward.
Since $[L:\ltpd]\leq 2\cdot(2^4\cdot 3) \cdot(2^5\cdot 3\cdot 5)=2^{10}\cdot 3^2\cdot 5$, and $S^\Q_{wild}=\{2,3,5\}$, one gets $$\#S^\Q_{wild}\cdot\log([L:\ltpd])\leq 3\cdot \log(2^{10}\cdot 3^2\cdot 5)=3\times 10.8\leq 33$$
which is the asserted bound.
\ep
\blem\label{le:basic-ineq2} 
Under the notations and assumptions of \Cref{th:main-bound},
one has 
\benumlab
\item $$\log(\fd_{L'})\leq \log(\fd_{L'})+\log(\ff_{L'})\leq \log(\fd_{\ltpd})+\log(\ff_{\ltpd})+\mycancelto{4}{2}\cdot\log(\ell)+\mycancelto{33}{21}$$ 
and hence,
\item $$\left(1+\frac{4}{\ell}\right)\log(\fd_{L'})\leq \left(1+\frac{4}{\ell}\right) \left(\log(\fd_{\ltpd})+\log(\ff_{\ltpd})\right)+\mycancelto{4}{2}\log(\ell)+\mycancelto{74}{46}.$$
\eenum
\elem
\bp 
The inequality given in {\bf(1)} is also a consequence of the log-Different+log-Conductor Theorem (\Cref{th:diff-cond-thm}) and \Cref{le:tame-lemma}, applied to the extension $L'/L/\ltpd$, except for the explicit computation of the constant arising from wild ramification terms in \Cref{th:diff-cond-thm}. Since $L'/L$ is tamely ramified outside $\{\ell\}$, and as $$[L':L]\leq \ell^4$$ by \Cref{le:orders}, one obtains the bound $$\log([L':L])\leq 4\cdot \log(\ell)$$ for the wild ramification term  at $\ell$ in \Cref{th:diff-cond-thm}, for the extension $L'/L$ and for the extension  $L/\ltpd$ one  uses the bound given by \Cref{le:basic-ineq1}.

The assertion {\bf(2)} follows from {\bf(1)} by multiplying the inequality in {\bf(1)} by $\left(1+\frac{4}{\ell}\right)$ and using the fact that $\log(\ff_{L'})\geq0$ and that for $\ell\geq5$, one has $\frac{\log(\ell)}{\ell}\leq\frac{1}{2}$ and $\left(1+\frac{4}{\ell}\right)\leq 2$.
\ep

\subsection{Analysis of ramification}\nwss
\blem\label{le:unram}
Suppose that the notations and assumptions of \Cref{th:main-bound} are in force. Assume $v\in \vnon_{\ltpd}$ such that $v$ does not divide $2\cdot3\cdot5\cdot\ell$ and $v\not\in \supp(\fq_{\ltpd})$. Then $L'/\ltpd$ is unramified over $v$.
\elem
\bp This is the assertion \cite[(D0) on Page 654]{mochizuki-iut4} and it is an immediate consequence of the assumptions on primes in \cref{sss:field-assump} and \cite[Proposition 1.8(vi,vii)]{mochizuki-iut4}.
\ep

\subsection{}
Assume \inithtdata, \initassumptions\ are in force. Suppose $\lmod \subseteq \ltpd \subseteq M \subseteq L'$ is an extension such that $M/\lmod$ is Galois. Consider a subset of primes  $\vmdst \subset \vmnon$ such that $v\in \vmdst$, if and only if, $v$ extends to a prime of $L'$ which is ramified over $\Q$.

\blem\label{le:vdst}
Under the notations and assumptions of \Cref{th:main-bound},
the following conditions are equivalent 
\benumlab
\item $v\in\vmdst$.
\item $v$ divides $30\cdot\ell$ or $v\in \supp(\fq_M+\fd_{M})$.
\item The image of $v\in\V_{\ltpd}$ is contained in $\V^{dst}_{\ltpd}$.
\eenum
In particular $\vmdst$ (and hence its image in $\V_{\ltpd}$ is finite).
\elem 
\bp  By the definition of $L'$ and assumptions \initassumptions, \cref{sss:field-assump} and the properties of the Weil pairings \cite{silverman-arithmetic} one sees that $L'\supset \Q(\zeta_{60\cdot\ell})$.  If $v_M$ extends  to some prime $v$ of $L'$ which is ramified then $v|\fd_{M}$ or $v|\fq_{M}$ or $v|2\cdot 3\cdot 5\cdot\ell$. This proves the assertion.
\ep

\subsection{Ramification divisors}
\newcommand{\vdst}[1]{\V^{dst}_{#1}}
\newcommand{\fs}{\mathfrak{s}}
\newcommand{\fts}{\mathfrak{s}^\leq}

Define $\vdst{\lmod}\subset \vnon_{\ltpd}$ (resp. $\vdst{\Q}\subset \vnon_{\Q}$) to be the images of $\vdst{\ltpd}$ in $\V_{\lmod}$ (resp. in $\V_{\Q}$). For a rational prime $v_\Q$ and for $L_*\in\{ \lmod,\Q\}$, write $$\V_{L_*,v_\Q}=\V_{L_*}\times_{\V_\Q}\{v_\Q\}.$$

The proof of the next lemma is similar to that of \Cref{le:vdst}:
\blem\label{le:vdst2}
Under the notations and assumptions of \Cref{th:main-bound},
the following conditions are equivalent 
\benumlab
\item $v_\Q\in\vqdst$.
\item $v_\Q$ ramifies in $L'$.
\item $v_\Q$ divides $30\cdot\ell$ or $v_\Q$ is in the image of $\supp(\fq_{\ltpd}+\fd_{\ltpd})$.
\item $v_\Q$ divides $30\cdot\ell$ or $v_\Q$  is in the image of $\supp(\fq_{L}+\fd_{L})$.
\eenum
\elem

Define ramification divisors for primes given by \Cref{le:vdst} as follows. Let
\be
\fs_{L_*}=\sum_{v\in\vdst{L_*}}e_v\cdot v.
\ee
Then \be \deg(\fs_{L_*})\in\R_{\geq0}, \ee
be its degree and let
\be \log(\fs_{L_*})= \frac{\deg(\fs_{L_*})}{[L_*:\Q]}\in \R_{\geq0}\ee
be its normalized degree.

For a prime $v$ of a field, write $p_v$ for its residue characteristic. Let 
\be \iota_v=\begin{cases}
1 & \text{if } p_v\leq \emods\cdot\ell\\
0&  \text{ otherwise}. 
\end{cases}\ee
Let 
\be\label{def:fts} \fts_{L_*}=\sum_{v\in\vdst{L_*}}\frac{\iota_v}{p_v}\cdot v.\ee
Then again one has the degree of $\fts_{L_*}$ and its normalized degree \be \deg(\fts_{L_*}), \log(\fts_{L_*})\in\R_{\geq0}\ee
given by the formulae stated earlier.

\blem\label{le:easy-estimate-fts} 
Under the notations and assumptions of \Cref{th:main-bound}, 
one has $$\emods\cdot \log(\fts) \leq \frac{4}{3}\left(\emods\cdot\ell+\eta_{prm} \right).$$
\elem 
\bp 
Since $\fts$ is supported on primes $\leq \emods\cdot\ell$, this is a straight forward argument using the definition of $\fts$ of \cref{def:fts} and is given in \cite[Page 660]{mochizuki-iut4} using definitions and basic explicit estimate \cref{sss:prime-bound} for the prime counting function.
\ep
\subsection{}
\bpro\label{pr:rampro} Under the notations and assumptions of \Cref{th:main-bound}, one has 
$$ \left(1+\frac{4}{\ell}\right)\cdot\log(\fd_{L'}) +\frac{4}{\ell}\cdot\log(\fs_{\Q}) \leq \left(1+\frac{12\cdot \dmod}{\ell}\right)\cdot \left(\log(\fd_{\ltpd})+\log(\ff_{\ltpd})\right)+\mycancelto{4}{2}\cdot\log(\ell)+\mycancelto{80}{52}.
$$
\epro

\begin{proof}[Proof of \Cref{pr:rampro}]
Let $v\in \vnon_M$ be a prime, write $p_v$ for the characteristic of the residue field on $v$ and write $f_v$ for the degree of the residue field. Then the logarithm of the cardinality of the residue field at $v$ is $\log v=p_v^{f_v}$ and so $\log v=f_v\cdot \log(p_v)$. 

For an arithmetic divisor $\fa_M$ and a rational prime $v_\Q\in\V_{\Q}$, let $$\fa_{M,v_\Q}=\sum_{w|v_\Q} e_{w}\cdot w$$ be the $v_\Q$ component of $\fa_M$. First let me remark that
\be	\frac{\log(\fs_{M,v_\Q})}{p_v}=\frac{\log(p_v)}{p_v}.\ee
This is immediate from the well known equation $$\sum e_i\cdot f_i=n$$ applied to the prime $v_\Q$ and the extension $M/\Q$.  Hence for $v_\Q=\ell$ one has
\be 
\frac{\log(\fs_{M,\ell})}{\ell}=\frac{\log(\ell)}{\ell}\leq \frac{1}{2} \text{ for } \ell\geq 5.
\ee

Next I claim that, for any $p_{v_\Q}\not\in\{2,3,5,\ell\}$ that
\be\label{eq:svq-eq-tame-portion1} 
\frac{\log(\fs_{\lmod})_{v_\Q}}{p_{v_\Q}}\leq \frac{2}{p_{v_\Q}}\cdot(\log(\fd_{\ltpd,v_\Q})+\log(\ff_{\ltpd,v_\Q})),
\ee
for which it is enough to prove , for $p_{v_\Q}\not\in \{2,3,5,\ell \}$, that
\be\label{eq:svq-eq-tame-portion2} 
\log(\fs_{\lmod})_{v_\Q}\leq 2\cdot(\log(\fd_{\ltpd,v_\Q})+\log(\ff_{\ltpd,v_\Q})).
\ee
Mochizuki asserts this as a consequence of \cite[Proposition 1.3(i)]{mochizuki-iut4}, but I will give a different proof using the elementary lower bound on differents given by \cite[Chap. III, Proposition 13]{serre1979-local-fields} applied to the extensions $\lmod/\Q$. For any $v\in\V_{\lmod}$, write $\delta_v=\fd_{\lmod,v}$. Then
$$ e_v-1\leq \delta_v$$ 
equivalently $$e_v\leq\delta_v+1.$$

Two cases need to be considered: $v\in\supp(\fq_{\lmod})$ and $v\not\in\supp(\fq_{\lmod})$. In the first case  one can write the inequality $e_v\leq\delta_v+1$
as the inequality $$e_v\leq\delta_v+\ff_v.$$
On the other hand if $v\not\in\supp(\fq_{\lmod})$ then $\ff_v=0$ and so one can write
\be e_v\leq \delta_v+1\leq 2(\delta_v+0)=2(\delta_v+\ff_v).\ee
Thus in all cases, for primes $v$ not dividing $30\cdot \ell$ one has
\be e_v\leq 2(\fd_{\lmod,v}+\ff_v)\ee and hence
\be e_v\cdot\log v\leq 2(\fd_{\lmod,v}+\ff_v)\cdot\log v.\ee
Summing this over all $v\in\vnon_{\lmod}$ not dividing $30\cdot\ell$ and dividing the resulting sum by $[\lmod:\Q]$ one obtains  the normalized degree of the component for primes $v_\Q\not\in \{2,3,5,\ell\}$:
\be \log(\fs_{\lmod})_{v_\Q}\leq 2(\log\fd_{\lmod}+ \log\ff_{\lmod})_{v_\Q}.\ee
Since $\ltpd/\lmod$ is tamely ramified outside $\{2,3,5,\ell\}$, by the log-Discriminant+log-Conductor Theorem (\Cref{th:diff-cond-thm}), one has \be \log(\fs_{\lmod})_{v_\Q}\leq 2(\log\fd_{\lmod}+ \log\ff_{\lmod})_{v_\Q}= 2(\log\fd_{\ltpd}+ \log\ff_{\ltpd})_{v_\Q}.\ee
This proves \eqref{eq:svq-eq-tame-portion2}  and hence \eqref{eq:svq-eq-tame-portion1} follows.

Next I claim that
\be 
\log(\fs_\Q)\leq 2\cdot\dmod\cdot(\log(\fd_{\ltpd})+\log(\ff_{\ltpd}))+\log(2\cdot 3\cdot5\cdot\ell).
\ee
This is proved as follows. Consider the divisor $\log(\fs_{\Q})$, by definition this is simply
\be \log(\fs_{\Q})=\sum_{p\in\vdst{\Q}} 1\cdot \log(p)=\sum_{p\in\vdst{\Q}-\{2,3,5,\ell\} } 1\cdot \log(p)+\left(\log 2+\log 3+\log 5 +\log \ell\right).\ee 
Now one estimates the main term of the above equation. Consider  the trivial estimate for the field $\lmod$ (this means $e_v,f_v$ refer to $\lmod$)
\be 1 \leq \sum_{v|p} [\lmodv:\Q_p]=\sum_{v|p} e_v\cdot f_v.\ee
Hence multiplying by $\log(p)$ one gets
\be 1\cdot \log(p) \leq \sum_{v|p} [\lmodv:\Q_p]\cdot \log p=\sum_{v|p} e_v\cdot f_v\cdot \log p=\sum_{v|p}e_v\cdot \log v.\ee
The right hand side is $\deg((\fs_{\lmod})_v)$. Summing over all $v|p$ and all $p\not\in \{2,3,5,\ell \}$ one gets
\be \log(\fs_{\Q})\leq \dmod\cdot \frac{\deg((\fs_{\lmod}))}{\dmod}+\log(2\cdot3\cdot5\cdot\ell).
\ee
Since $\frac{\deg((\fs_{\lmod}))}{\dmod}=\log(\fs_{\lmod})$  first term may be bounded by \eqref{eq:svq-eq-tame-portion2}:
\be \log(\fs_{\Q})\leq 2\cdot\dmod\cdot(\log(\fd_{\ltpd})+\log(\ff_{\ltpd}))+\log(2\cdot3\cdot5\cdot\ell),
\ee
Using $\log(2\cdot 3\cdot 5)<5$ one gets
\be 
\log(\fs_\Q)\leq 2\cdot\dmod\cdot(\log(\fd_{\ltpd})+\log(\ff_{\ltpd}))+5+\log(\ell)
\ee
and multiplying this by $\frac{4}{\ell}$ that
\be 
\frac{4}{\ell}\cdot\log(\fs_\Q)\leq \frac{8}{\ell}\cdot\dmod\cdot(\log(\fd_{\ltpd})+\log(\ff_{\ltpd}))+\frac{20}{\ell}+\frac{4\cdot\log(\ell)}{\ell}.
\ee
As $\ell\geq 5$, so $20/\ell\leq 4$ and $\log(\ell)/\ell\leq\frac{1}{2}$ and hence this reduces to
\be 
\frac{4}{\ell}\cdot\log(\fs_\Q)\leq \frac{8}{\ell}\cdot\dmod\cdot(\log(\fd_{\ltpd})+\log(\ff_{\ltpd}))+6.
\ee

Combining this with the inequality provided by \Cref{le:basic-ineq2} gives
\begin{multline}
\left(1+\frac{4}{\ell}\right)\cdot\log(\fd_{L'})+\frac{4}{\ell}\cdot \log(\fs_\Q)\leq  \left(1+\frac{12\cdot\dmod}{\ell}\right)\cdot(\log(\fd_{\ltpd})+\log(\ff_{\ltpd}))\\+\mycancelto{4}{2}\cdot\log(\ell)+\mycancelto{80}{52}.
\end{multline}
This completes the proof of \Cref{pr:rampro}.
\end{proof}

\subsection{Upper bounds for the theta-values locus}
This section uses the notation and the results of \cite[{\ssep\ref{III-se:mochizuki-construction-thetam}}]{joshi-teich-rosetta}. Let $\bsIm$ (resp. $\bsImq$) be the tensor-packet codomains constructed in \cite[{\ssep\ref{III-ss:tensor-packet-codomain}}]{joshi-teich-rosetta} and let $$\thetami\subset \bsIm$$ be the theta-values locus constructed in \cite[{\ssep\ref{III-ss:construction-thetam-thetaj}}]{joshi-teich-rosetta}. 

Of interest to us in this section is an upper bound for the volume (or the log-volume) of 
$$\Vol(\thetami).$$ Readers should note that the fundamental estimate (i.e. \moccor) established in  \cite[{Corollary~\ref{III-cor:cor312}}]{joshi-teich-rosetta} provides a lower bound for $\Vol(\thetami)$. Presently one is seeking upper bounds for $\Vol(\thetami)$. In \cite{mochizuki-iut4}, this upper bound estimate is embedded within the proof of \cite[Theorem 1.10]{mochizuki-iut4}.

\subsection{The Second Main Bound: The log-Volume bound}\label{ss:log-vol}  The construction of $\thetami$ given in \cite[{\ssep\ref{III-se:mochizuki-construction-thetam}}]{joshi-teich-rosetta} is adelic and theta-values are bundled  by rational primes $p$. An important point to observe is that this bundling collects all primes of $L'$ lying over $p$ regardless of the reduction type of the curve. This is precisely the sort phenomenon one encounters in context of the \cite[{Theorem \ref{III-th:hyperplane-theta-lnk}}]{joshi-teich-rosetta} and especially  \cite[{Remark \ref{III-re:hyperplane-theta-lnk}}]{joshi-teich-rosetta} (which assert that altering of valuations at primes of semi-stable reduction implies that the valuations at non-semi-stable reductions must also be altered). For compatibility with \cite{mochizuki-iut4}, I will write $v_{\Q}$ for this rational prime $p$ or write $p_{v_\Q}$ as well (note that $p_v$ is also the characteristic of the residue field of any prime of $L'$ lying over $p$. 

As noted in \cite[{\ssep\ref{III-se:mochizuki-construction-thetam}}]{joshi-teich-rosetta}, the global volume will be the product of the local volumes at each prime $p\in\V_{\Q}$ (including the archimedean prime $p=\infty$) and is finite because it receives contribution of the factor equal to
$$1\in \R$$ at all but finitely many rational primes (in particular, the local log-volume is $\log(1)=0\in\R$ at all but finitely many rational primes). 

\textit{For compatibility with \cite{mochizuki-iut4}, I will follow the notational convention for $\log$-volumes of \cite[{\ssep \ref{III-ss:moccor-proof}}]{joshi-teich-rosetta}.}
\newcommand{\thetamiphibyz}{{\thetam}^{\mathscrbf{I},\bvphi(\by_0)}}
\newcommand{\thetamibyz}{{\thetam}^{\mathscrbf{I},\by_0}}

As indicated in \cref{ss:outline}, I will use two different, arithmeticoids $\by_0'=\bvphi(\by_{0})$ and $\by_0$ for my next assertion--however I will concurrently normalize the (same) number field provided by each of these arithmeticoids. I will indicate the usage of $\by_0'$ or $\by_0$ by a superscript but this convention will be dropped whenever there is no risk of confusion to keep the notation uncluttered. Working with $\by_{0}'=\bvphi(\by_{0})$ and $\by_0$ simultaneously  corresponds in \cite[Remark 3.3.3]{mochizuki-iut3} to working with $\flog$-shifted versions  of various objects. By \constrthr{\ssep}{ss:log-link-indentified}, Mochizuki's $\flog$-Link is identified with the global Frobenius $\bvphi$ morphism of the present theory and hence working with $\flog$-shifted versions is the same as working with Frobenius-shifted versions.  To understand what working with Frobenius-shifted objects means geometrically see \cref{re:bundle-rmk} (given below).

\bthm\label{th:theta-upper-bound}
The preceding notation, convention and assumptions (notably the existence of Initial Theta-Data \inithtdata, \initassumptions\ is given by \Cref{th:existence}) are strictly in force. Let $\by_0$ be the standard arithmeticoid chosen in \cref{ss:basic-setup} and let $\bf{q}^{\by_0}_\ell=\frac{1}{2\ell}\cdot\log(\fq)$ be computed using the arithmeticoid $\by_0$.   Let $\bvphi$ be the global Frobenius morphism of \constrtwoh{Definition }{def:global-Frobenius-def} and let $\by_0'=\bvphi(\by_0)$. Choose the same global normalization for the number field provided by the arithmeticoids $\by_{0}'$ and $\by_0$. With this choice of simultaneous global normalization one has  $$-\abs{\log(\bf{q}^{\bvphi(\by_0)}_\ell)} = -\abs{\log(\bf{q}^{\by_0}_\ell)},$$  and one has
$$-\abs{\log(\bf{q}^{\by_0}_\ell)}\leq -\frac{1}{\ells}\abs{\log\Vol(hull(\thetamiphibyz))}= -\frac{1}{\ells}\abs{\log\Vol(hull(\thetamibyz))} \leq C_\Theta\cdot \abs{\log(\bf{q}^{\by_0}_\ell)},$$
where the first inequality is given by \moccor\ via \cite[{Corollary \ref{III-cor:cor312}}]{joshi-teich-rosetta} and where
\begin{multline*}
C_\Theta=\left\{
\frac{\ell+1}{4 \abs{\log(\bf{q}_\ell)}} \left(
\left(1+\frac{12\cdot\dmod}{\ell}\right)\right.\right.\\
 \left.\left. \cdot(\deg(\fd_{\ltpd})+\deg(\ff_{\ltpd})) + 
 10(\emods\cdot \ell+\etap) - 
 \frac{1}{6} \left(1-\frac{12}{\ell^2}\right)\cdot 
 \log(\fq)\right)-1
\right\}
\end{multline*}
\ethm

\brem
Note however that in \Cref{th:theta-upper-bound}, as arithmeticoids $\by_{0}'\neq \by_0$ as they differ by Frobenius at each prime $v\in\vl$ and this, together with choice of concurrent global normalization for the copy of the (same) number field they provide precludes direct comparison of various local quantities.
\erem

\brem\label{re:bundle-rmk}
Let me explain the consequences of working with $\by_0'$ and $\by_0$ simultaneously. Let $C/L$ be an elliptic curve over a number field $L$. Let $p_v$ be the residue characteristic of $v\in\vlnon$. By \cite{fargues-fontaine} one has a collection of filtered rank two vector bundles $$\{(\sV_v,Fil(\sV_v))\}_{v\in\vlnon}$$ on $\sX_{\cpt,L_v}$ and a similar, but more complicated description, on $\sY_{\cpt,L_v}$ (via \cite{kedlaya-liu15}). [For an elliptic curve $C/L$, $\sV_v$ is the vector bundle on $\sY_{\cpt,L_v}$ associated to $C/L_v$.] If $\varphi_v$ is the local Frobenius of $\sY_{\cpt,L_v}$ then working with $\by_0$ and $\by_0'=\bvphi(\by_0)$ simultaneously in \Cref{th:theta-upper-bound} is analogous to  comparing $$\{(\sV_v,Fil(\sV_v))\}_{v\in\vlnon}\text{ with } \{(\sV_v,\varphi_v(Fil(\sV_v)))\}_{v\in\vlnon}.$$ Mochizuki refers  to this process as a $\flog$-shift (\cite[Remark 4.8.2(iii)]{mochizuki-iut2} and \cite[Remark 3.3.3]{mochizuki-iut3}) and because a $\flog$-Link identified with the global Frobenius of the present theory (see \constrthr{\ssep}{ss:log-link-indentified}), $\flog$-shifting is the same as  (global) Frobenius shifting as above. 

Let me also remark that in the geometric contexts, inequalities for the slopes of ``Frobenius-shifted bundles'' occur in \cite{joshi16}.  Roughly speaking, the inequalities of \cite{joshi16} involve comparison of slopes of this type of vector bundles (for curves in characteristic $p>0$), and similar slope inequalities underlie the proof of the Geometric Szpiro Inequality due to \cite{kim97}, \cite{beauville02}.  Thus Frobenius shifting (or $\flog$-shifting) is natural in the context of the Arithmetic Szpiro Inequality.

To emphasize the analogy between \iut, \cite{joshi16} and \cite{beauville02}, \cite{kim97}, suppose  that $C/L$ has semi-stable reduction, for each prime $v\in\vlnon$ for $C/L$, there exists a connection  $\nabla_v$, such that $\{(\sV_v,Fil(\sV_v),\nabla_v)\}_{v\in\vlnon}$ is an indigenous bundle (up to passage to the projective bundle). Mochizuki had remarked to me (in Spring 2018) that he views an elliptic curve over a number field as an arithmetic analog of an indigenous bundle and that this analogy underlies his theory. The inequalities established in \cite{joshi16} are established using higher rank analogs of indigenous bundles i.e. using higher rank opers. The construction in \cite[Section 4.3]{joshi-formal-groups} of the global arithmetic analog of Fontaine's fundamental ring  $\ainf$  of $p$-adic Hodge Theory, suggests that there is direct construction of such a rank two arithmetic indigenous bundle associated to $C/L$. Because of this, one expects  that (eventually) one may arrive at a proof of the Arithmetic Szipro Inequality along the lines of \cite{joshi16}, \cite{beauville02}, \cite{kim97}.
\erem

\brem 
As mentioned earlier, this theorem and its proof is embedded in the proof of \cite[Theorem 1.10]{mochizuki-iut4} and this makes the proof of the latter theorem difficult to follow.
\erem
\bp 
The equality $$-\abs{\log(\bf{q}^{\bvphi(\by_0)}_\ell)} = -\abs{\log(\bf{q}^{\by_0}_\ell)}$$ is a consequence of the choice of the same global normalization for the copy of the (same) number field provided by $\by_0'$ and $\by_0$. The lower bound inequality  is proved by invoking
\constrthr{Theorem }{th:moccor-Mochizuki-form2} but computed using $\by_0'=\bvphi(\by_{0})$.
The  inequality on the right is computed using the standard arithmeticoid $\by_0$. While a similar strategy is  employed in \cite{mochizuki-iut3} especially see \cite[Remark 4.8.2(iii)]{mochizuki-iut2}, in Mochizuki's formulation of \moccor, the ``Frobenius-shifted'' (``$\flog$-shifted'' in the parlance of \cite[Remark 3.3.3]{mochizuki-iut3}) nature of the inequality is not made explicit.

The middle equality is a tautology using the fact that the number fields provided by $\by_0'$ and $\by_0$ are identically normalized.

So now it remains to prove the inequality on the right. This only involves the arithmeticoid $\by_0$ and hence, in the remainder of the proof, I will drop the notation superscript referencing the arithmeticoid $\by_0$ for notational simplicity.

As is detailed in \constrthr{\ssep}{se:mochizuki-construction-thetam}, the set $\thetami$ is adelic by construction and one may write
$$\thetami=\prod_{p\in\V_{\Q}}\thetampi\subset \bsIm\subset \bsImq.$$
The  volume of $\thetami$ is estimated by estimating volume of component $\thetampi$ for each rational prime  $p$ separately. The tensor product structure of $\bsI$ (and $\thetampi$) is detailed in \cite[{\ssep\ref{III-se:mochizuki-construction-thetam}}]{joshi-teich-rosetta}. The key tool one uses in the calculation of the upper bound for $\log\mathit{Vol}(\thetampi)$ is this tensor structure and \cite[Proposition 1.4, Proposition 1.7]{mochizuki-iut4}. Mochizuki details this computation in Steps (v), (vi) and (vii) of \cite[Proof of Theorem 1.10]{mochizuki-iut4}. These three steps in the \cite[Proof of Theorem 1.10]{mochizuki-iut4} correspond to the following three cases:
\benumlab
\item $v_\Q\in \vdst{\Q}$ (this means $v_\Q$  extends to a ramified prime in $L'$),
\item $v_\Q\in \vnon_{\Q}-\vdst{\Q}$, finally,
\item $v_\Q\in\varc_\Q$ i.e. $v_\Q=\infty$ is the unique archimedean prime of $\Q$.
\eenum
\ep
Mochizuki's computation in this case is easy to follow and I will not repeat the details here (unless required for completeness) as I do not believe any further comments are needed on these proofs.

The proof requires estimates for the numbers $e_v$ appearing in various versions of $\fs$. These estimates are my next lemma. I found the assertions regarding this in \cite[(R1)--(R4), Page 655]{mochizuki-iut4} sufficiently unclear and even confusing (\cite{fucheng} does not prove the assertion and \cite{yamashita} is also unclear (to me)). So I felt the need to write out a proof. 
\blem\ 
Let $(L'\supseteq) M\supseteq L\supseteq \lmod\supseteq \Q$, let  $u\in \V_M$ and let $u|w|v|p$ be the primes lying below $u$ in each of these extensions.
Then one has $$\log(e_u)\leq -3 + 4\cdot\log(\emods\cdot\ell).$$
\elem
\bp 
Since $ L'\supseteq M\supset L\supseteq \lmod\supseteq \Q$, let $u|w|v|p$ be primes lying over $p$ in each of these extensions.  Then $e_u=e_{u|p}^{M/\Q}$. Recall that, for any finite extension $M/L/\Q$  one has the equation $\sum_i e_i\cdot f_i=[M:L]$ and hence one has $e_i\leq [M:L]$. Applying this to the problem of estimating $e_u$ one  has, by \Cref{le:deg} and the definition of $\emod, \emods$ that
\begin{align} 
\begin{split} 
e_u^{M/\Q} & = e_{u|w}\cdot e_{w|v}\cdot e_{v|p}\\
&\leq  [L':\lmod]\cdot \emod\\
&\leq  [L':L] \cdot [L:\lmod] \cdot \emod\\
&\leq \ell^4\cdot 2\cdot \#GL_2(\F_2)\cdot \#GL_3(\F_2)\cdot \#GL_2(\F_2)\cdot \emod\\
&\leq \emods\cdot \ell^4.
\end{split}
\end{align}
Hence 
\be 
\log(e_u)\leq \log(\emods\cdot \ell^4)
\ee
and as $e\leq 3\leq \emods$ so one sees that $-3+4\cdot \log(\emods\cdot \ell)\geq \log(\emods\cdot \ell^4)$ and hence one can replace the above inequality by the weaker inequality
\be 
\log(e_u)\leq -3+4\cdot \log(\emods\cdot \ell).
\ee
This finishes the proof.
\ep

\bpro\label{pr:theta-upper-bound1} 
Let $p=v_\Q\in\vdst{\Q}$, then under the hypothesis of \Cref{th:theta-upper-bound} one has 
\begin{multline*}
-\frac{1}{\ells}\cdot \abs{\log \Vol(\thetampi)}\leq \frac{\ell+1}{4}\left\{ \left(1+\frac{4}{\ell}\right)\log(\fd_{L',v_\Q})-\frac{1}{6}\cdot \log(\fq_{v_\Q}) \right. \\ 
\left. +\frac{4}{\ell}\cdot\log(\fs_{\Q,v_\Q})+\frac{20}{3}\cdot \emods\cdot\log(\fts_{v_\Q}).
\right\}
\end{multline*}
\epro

\bp 
This is the last equation on \cite[{Step (v) on Page 658, Proof of Theorem 1.10}]{mochizuki-iut4} and its proof is all of step (v). 
\ep

\bpro\label{pr:theta-upper-bound2} 
Let $p=v_\Q\in \vnon_{\Q}-\vdst{\Q}$, then under the hypothesis of \Cref{th:theta-upper-bound} one has
$$
-\frac{1}{\ells}\cdot\abs{\log \Vol(\thetampi)}\leq 0.
$$
\epro
\bp 
This is \cite[Step (vi) on Page 659, Proof of Theorem 1.10]{mochizuki-iut4}.
\ep

\brem 
For the archimedean case, considered in the next proposition, I would like to add one conceptually clarifying remark to Mochizuki's computations. Mochizuki treats the theory at any archimedean prime of $L$ (and of $\Q$) using his anabelian reconstruction philosophy. 

However, my remark is that the archimedean case can be treated directly using classical Teichmuller Theory. For $g=1$ the classical Teichmuller Space is $\mathfrak{H}\subset \C$  the upper half plane which is biholomorphic to the open unit disk and this Teichmuller description may be used to understand the archimedean $g=1$ case. This leads to a more direct way of understanding Mochizuki's volume computations at archimedean primes. 
\erem
\bpro\label{pr:theta-upper-bound3}
Let $v_\Q\in\varc_\Q$ i.e. $p=\infty$. Then under the hypothesis of \Cref{th:theta-upper-bound}
$$-\frac{1}{\ells}\cdot\abs{\log \Vol(\thetampi)}\leq (\ell+1).$$
\epro
\bp
This is established in \cite[{Step (vii) on Page 660, Proof of Theorem 1.10}]{mochizuki-iut4}.
\ep

\subsection{Proof of the Second Main Bound (\Cref{th:theta-upper-bound})}
\bp[Proof of \Cref{th:theta-upper-bound}]
Now let me complete the proof of estimate given in \Cref{th:theta-upper-bound}. This is obtained by combining the upper bounds provided by \Cref{pr:theta-upper-bound1,pr:theta-upper-bound2,pr:theta-upper-bound3}.

One has by definition (for the notational conventions see \cite[{\ssep \ref{III-ss:moccor-proof}}]{joshi-teich-rosetta})
\begin{multline}\label{eq:component-vol-added}
-\frac{1}{\ells}\cdot\abs{\log\Vol(hull(\thetami))}=-\frac{1}{\ells}\cdot\left(\sum_{p=\in\vdst{\Q}} \abs{\log\Vol(hull(\thetampi))} \right) \\ 
	+ \abs{\log\Vol(hull(\thetaminfty))},
\end{multline}
with primes $v_\Q\not\in\vdst{\Q}\cup\{\infty\}$ contributing a zero by \Cref{pr:theta-upper-bound2}.
Hence one obtains from \eqref{eq:component-vol-added},  \Cref{pr:theta-upper-bound1,pr:theta-upper-bound2,pr:theta-upper-bound3} that
\begin{align}\label{eq:mess1}
\begin{split}
-\frac{1}{\ells}\cdot\abs{\log\Vol(hull(\thetami))} \ \leq \ 
& \frac{\ell+1}{4}\left\{ \left(1+\frac{4}{\ell}\right)\cdot   \log(\fd_{L'}) \right. \\ 
& \left. -\frac{1}{6}\cdot \log(\fq_{\Q}) +\frac{4}{\ell}\cdot\log(\fs_{\Q})  \right. \\ 
& \left. +\frac{20}{3}\cdot \emods\cdot\log(\fts_{\Q}) \right\} \\
& +(\ell+1).
\end{split} 
\end{align}

Now one uses \Cref{pr:rampro} and \Cref{eq:mess1} and a bit of simplification to get:
\begin{align}\label{eq:mess2}
\begin{split}
-\frac{1}{\ells}\cdot\abs{\log\Vol(hull(\thetami))} \   \leq \ 
& \left(\frac{\ell+1}{4}\right)\cdot\left\{  -\frac{1}{6}\cdot \log(\fq) \right. \\ 
&+ \left. \left(1+\frac{12\cdot\dmod}{\ell}\right)\cdot   \left(\log(\fd_{\ltpd}) 
 + \log(\ff_{\ltpd})\right)  \right. \\ 
& \left. +\frac{20}{3}\cdot \emods\cdot\log(\fts_{\Q}) + \mycancelto{4}{2}\log(\ell) +\mycancelto{80}{56} \right\}.
\end{split}
\end{align}
Next note that one has, by elementary algebra, that
$$ \frac{\ell+1}{4}\cdot \left(-\frac{1}{6}\right)+ \frac{\ell+1}{4}\cdot\frac{1}{6}\cdot \frac{12}{\ell^2}-\frac{1}{2\ell}=\frac{\ell+1}{4}\cdot \left(-\frac{1}{6}\right)+\frac{1}{2\ell^2}$$
and hence
$$
\frac{\ell+1}{4}\cdot \left(-\frac{1}{6}\right)\cdot\left(1-\frac{12}{\ell^2}\right)-\frac{1}{2\ell}
=
\frac{\ell+1}{4}\cdot \left(-\frac{1}{6}\right)+ \frac{\ell+1}{4}\cdot\frac{1}{6}\cdot \frac{12}{\ell^2}-\frac{1}{2\ell} > \frac{\ell+1}{4}\cdot \left(-\frac{1}{6}\right).$$
Hence one may replace the term $-\frac{1}{6}\cdot \log(\fq)$ in \eqref{eq:mess2} by $-\frac{1}{6}\cdot\left(1-\frac{12}{\ell^2}\right)\cdot \log(\fq)$ and still have an excess contribution greater than $\frac{1}{2\ell}\log(\fq)$ i.e
\begin{align}\label{eq:mess3}
\begin{split}
-\frac{1}{\ells}\cdot\abs{\log\Vol(hull(\thetami))} \   \leq \ 
& \left(\frac{\ell+1}{4}\right)\cdot\left\{  -\frac{1}{6}\cdot \left(1-\frac{12}{\ell^2}\right)\cdot \log(\fq) \right. \\ 
&+ \left. \left(1+\frac{12\cdot\dmod}{\ell}\right)\cdot   \left(\log(\fd_{\ltpd}) 
+ \log(\ff_{\ltpd})\right)  \right. \\ 
& \left. +\frac{20}{3}\cdot \emods\cdot\log(\fts_{\Q}) + \mycancelto{4}{2}\log(\ell) +\mycancelto{80}{56} \right\}-\frac{1}{2\ell}\log(\fq).
\end{split}
\end{align}
The term $\frac{20}{3}\cdot\emods\cdot\log(\fts_{\Q})$ in \eqref{eq:mess3} is estimated by \Cref{le:easy-estimate-fts} to be
$$\frac{20}{3}\cdot\emods\cdot\log(\fts_{\Q})\leq \frac{20}{3}\cdot\frac{4}{3}\cdot(\emods\cdot\ell+\etap).$$
Further as $\ell\geq 5$ one has
\be 
\mycancelto{\frac{10}{9}}{\frac{4}{9}}\left(\emods\cdot\ell+\eta_{prm} \right)\geq \mycancelto{\frac{10}{9}}{\frac{4}{9}}\cdot\emods\cdot\ell \geq \mycancelto{10}{4}\cdot 2^{12}\cdot3\cdot 5\cdot \ell  \geq \mycancelto{4}{2}\log(\ell)+\mycancelto{80}{56}.
\ee
Hence 
\begin{align*}
\frac{20}{3}\cdot \emods\cdot\log(\fts_{\Q}) + \mycancelto{4}{2}\log(\ell) +\mycancelto{80}{56}&\leq \frac{20}{3}\cdot\frac{4}{3}\cdot(\emods\cdot\ell+\etap)+\mycancelto{\frac{10}{9}}{\frac{4}{9}}\left(\emods\cdot\ell+\eta_{prm} \right)\\
&\leq\frac{90}{9}\cdot(\emods\cdot\ell+\etap),\\
&\leq\frac{90}{9}\cdot(\emods\cdot\ell+\etap),\\
&\leq 10\cdot(\emods\cdot\ell+\etap).
\end{align*}

Substituting this in \eqref{eq:mess3} and simplifying one obtains
\begin{align}\label{eq:mess4}
\begin{split}
-\frac{1}{\ells}\cdot\abs{\log\Vol(hull(\thetami))} \   \leq \ 
& \left(\frac{\ell+1}{4}\right)\cdot\left\{  -\frac{1}{6}\cdot \left(1-\frac{12}{\ell^2}\right)\cdot \log(\fq) \right. \\ 
&\qquad+ \left. \left(1+\frac{12\cdot\dmod}{\ell}\right)\cdot   \left(\log(\fd_{\ltpd}) 
+ \log(\ff_{\ltpd})\right)  \right. \\ 
&\qquad\qquad \left. + 10\cdot(\emods\cdot\ell+\eta_{prm}) \vphantom{\left(1-\frac{12}{\ell^2}\right)} \right\}-\frac{1}{2\ell}\log(\fq).
\end{split}
\end{align}

For brevity, write
\be\label{eq:rhs-of-cor312} \abs{\log({\bf q}_\ell)}:=\sum_{p,w\in\ubblvossp\neq\emptyset}\log\abs{q_w^{1/2\ell}}_{L_w'}.\ee 
Then the right hand side of \eqref{eq:rhs-of-cor312} is precisely the quantity which appears in number occurring in the right hand side of the fundamental estimate \moccor\ (proved in \cite[{Corollary \ref{III-cor:cor312}}]{joshi-teich-rosetta}). 

Hence observe that (by the definition of $\fq$ in \cref{ss:tate-divisor}) one sees that
$$\frac{1}{2\ell}\log(\fq)=\abs{\log({\bf q}_\ell)}$$
is precisely the number occurring in the fundamental estimate \moccor\ (proved in \cite[{Corollary \ref{III-cor:cor312}}]{joshi-teich-rosetta}). 
So by combining the (global) lower bound on the log-volume provided by \cite[{Corollary \ref{III-cor:cor312}}]{joshi-teich-rosetta} with the (global) upper bound on log-volume provided by \eqref{eq:mess4}
the (global) upper bound in \Cref{th:theta-upper-bound} follows. This proves \Cref{th:theta-upper-bound}.
\ep

\subsection{Proof of the First Main Bound (\Cref{th:main-bound})}
\bp[Proof of \Cref{th:main-bound}] 
Now one can complete the proof of \Cref{th:main-bound}.   By \Cref{th:theta-upper-bound}, 
one obtains the inequality:
\be -\abs{\log(\bq_{\ell})}\leq -\frac{1}{\ells}\cdot\abs{\log\Vol(\thetami)}\leq C_\Theta\cdot \abs{\log(\bq_{\ell})}.
\ee
As the (global) upper bound must be at least as big as the (global) lower bound so one deduces that
$$
C_\Theta\geq -1
$$
which on using the formula for $C_\Theta$ given in \Cref{th:theta-upper-bound} yields (on simplification)
\begin{multline*}
\frac{1}{6}\left(1-\frac{12}{\ell^2}\right)\cdot \log(\fq)\leq 
\left(1+\frac{12\cdot\dmod}{\ell}\right)\\
\cdot(\deg(\fd_{\ltpd})+\deg(\ff_{\ltpd})) + 
10(\emods\cdot \ell+\etap).
\end{multline*}
Finally  one uses $\ell\geq 7$ and the trivial estimates given below
\begin{align*}
\left(1-\frac{12}{\ell^2}\right)^{-1}\leq &\ 2\\
\left(1-\frac{12}{\ell^2}\right)^{-1}\cdot\left(1+\frac{12\cdot \dmod}{\ell}\right)\ \leq\ & \left(1+\frac{20\cdot \dmod}{\ell}\right)
\end{align*}
one obtains
\begin{multline*}
\frac{1}{6}\cdot \log(\fq)\leq 
\left(1+\frac{20\cdot\dmod}{\ell}\right)
\cdot(\deg(\fd_{\ltpd})+\deg(\ff_{\ltpd})) + 
20(\emods\cdot \ell+\etap).
\end{multline*}
to complete the proof of \Cref{th:main-bound}.
\ep

\section{Vojta's Inequality for $U_X(\bQ)^{\leq d}$}\label{se:vojta-ineq}\nwss
\subsection{Diophantine Inequalities}
This is the main theorem of this paper.
\bthm\label{th:main-dioph-thm} 
Let $X$ be a geometrically connected, smooth, projective curve over a number field $L$. Let $D\subset X$ be a reduced divisor, let $U_X=X-D$. Assume $U_X$ is a hyperbolic curve. Let $d\in\N$. Let $\varepsilon>0$. Then one has the inequality
$$ h_{\omega_X(D)}\lesssim (1+\varepsilon)(\logdiff_X+\logcon_D)$$
holds on $U_X(\bQ)^{\leq d}$.
\ethm
\bp 
By the main theorem of \cite{mochizuki-general-pos}, proving this assertion is equivalent to proving this assertion for $X=\P^1_\Q$, $D=\{0,1,\infty \}$ for compactly bounded subsets of $U_{\P^1}(\bQ)^{\leq d}$ satisfying the $2$-adic $j$-invariant condition \eqref{eq:2-adic-assumption}. Then it suffices to prove that
$$ h_{\omega_\P^1(D)}\lesssim (1+\varepsilon)(\logdiff_{\P^1}+\logcon_D)$$
holds on $\sZ\cap U_{\P^1}(\bQ)^{\leq d}$.

Consider an elliptic curve $C_\lambda$ given by \Cref{th:existence} and let $\exc$ be the set constructed in the proof of \Cref{th:existence}. This curve satisfies all the conditions of \Cref{th:main-bound} and the set for $x_\lambda\in\exc$, the height is bounded by some constant depending on $d,\sZ$. In the course of the rest of this proof this set will undergo further enlargement (by finite number of additions) which depend on $d,\varepsilon,\sZ$. For $C_\lambda\not\in \exc$, the estimate provided by \Cref{th:main-bound} holds. Hence
\begin{align*}
\begin{split}
\frac{1}{6}\cdot\log(\fq) & \leq \left(1+\frac{20\cdot\dmod}{\ell}\right)\cdot\left(\log(\dtpd)+\log(\ftpd)\right)+20\left( \emods\cdot\ell+60\right) \\
& \leq \left(1+\frac{20\cdot\dmod}{\ell}\right)\cdot\left(\log(\fd_L)+\log(\ff_L)\right)+20\left( \emods\cdot\ell+60\right).
\end{split}
\end{align*}
Next using $\ell\geq Q^{1/2}$ and $\emods\leq \dmod\leq \delta$, $\etap=60$ and \Cref{le:prm-bnd2} one obtains
\begin{align*}
\begin{split}
\frac{1}{6}\cdot\log(\fq) & \leq \left(1+\frac{20\cdot\dmod}{\ell}\right)\cdot\left(\log(\dtpd)+\log(\ftpd)\right)+20\left( \emods\cdot\ell+60\right) \\
& \leq \left(1+\frac{20\cdot\dmod}{\ell}\right)\cdot\left(\log(\fd_L)+\log(\ff_L)\right)+20\left( \emods\cdot\ell+60\right).\\
&\leq\left(1+\frac{\delta}{Q^{1/2}}\right)\cdot\left(\log(\dtpd)+\log(\ftpd)\right)+200\cdot \delta^2\cdot Q^{1/2}\cdot \log(2\cdot\delta\cdot Q)+20\cdot60.
\end{split}
\end{align*}

Recall that by definition $Q_2=\fT_2(C_\lambda)$ and $Q=\fT(C_\lambda)$   and \Cref{le:prm-bnd2} one has  $$\frac{1}{6}\cdot\fT_2(C_\lambda)-\frac{1}{6}\fq\leq \frac{1}{6}\cdot Q^{1/2}\cdot\log(\ell)\leq Q^{1/2}\log(2\cdot\delta\cdot Q),$$
or equivalently
\be \frac{1}{6}\cdot Q_2-\frac{1}{6}\fq\leq \frac{1}{6}\cdot Q^{1/2}\cdot\log(\ell)\leq Q^{1/2}\log(2\cdot\delta\cdot Q),\ee

From \Cref{pr:height-comparisons} one knows that the difference 
$$\frac{1}{6}\cdot Q-\frac{1}{6}\cdot Q_2\leq A_\sZ$$ for  some a positive real  number $A_{\sZ}$ which depends only on $\sZ$. Hence putting these equations together one obtains, with possibly a new $\sZ$-dependent constant denoted again by $A_{\sZ}$), that
\begin{align}
\begin{split}
\frac{1}{6}\cdot Q&\leq \left(1+\frac{\delta}{Q^{1/2}}\right)\cdot\left(\log(\dtpd)+\log(\ftpd)\right)+(15\cdot\delta)^2
 Q^{1/2}\log(2\cdot\delta\cdot Q)+\frac{1}{2}\cdot A_\sZ.
\end{split}
\end{align}
The last inequality can be written as
\be
\frac{1}{6}\cdot Q\leq
\left(1+\frac{\delta}{Q^{1/2}}\right)\cdot\left(\log(\dtpd)+\log(\ftpd)\right)+ \frac{1}{6}\cdot Q\cdot 
\frac{6\cdot (15\cdot\delta^2)\cdot\log(2\cdot\delta\cdot Q)}{Q^{1/2}}+\frac{1}{2}\cdot A_\sZ.
\ee
Now let 
\be\label{eq:alpha-lambda} \alpha_\lambda=\frac{6\cdot (15\cdot\delta^2)\cdot\log(2\cdot\delta\cdot Q)}{Q^{1/2}}\ee
Then if $1\leq \alpha_\lambda$ then $$Q^{1/2}\leq {6\cdot (15\cdot\delta^2)\cdot\log(2\cdot\delta\cdot Q)}$$
and hence $Q$ is bounded and hence there are only a finitely many elliptic curves in $U(\bQ)^{\leq d}$ where this happens. By again enlarging $\exc$ by adding these curves to it and without affecting other parameters. 
So one may assume that $\alpha_\lambda<1$.

Then the above inequality can be rewritten (with $\alpha_\lambda<1$) as 
\begin{align}
\frac{1}{6}\cdot Q&\leq
\left(1+\frac{\delta}{Q^{1/2}}\right)\cdot\left(\log(\dtpd)+\log(\ftpd)\right)+ \frac{1}{6}\cdot Q\cdot \alpha_\lambda+\frac{1}{2}\cdot A_\sZ.
\end{align}
equivalently
\begin{align}
(1-\alpha_\lambda)\frac{1}{6}\cdot Q&\leq
\left(1+\frac{\delta}{Q^{1/2}}\right)\cdot\left(\log(\dtpd)+\log(\ftpd)\right)+\frac{1}{2}\cdot A_\sZ.
\end{align}
or by dividing both the sides by $(1-\alpha_\lambda)$ one obtains
\begin{align}
\frac{1}{6}\cdot Q&\leq
\left(1+\delta\cdot Q^{-1/2}\right)\cdot(1-\alpha_\lambda)^{-1}\cdot\left(\log(\dtpd)+\log(\ftpd)\right)+(1-\alpha_\lambda)^{-1}\cdot\frac{1}{2}\cdot A_\sZ.
\end{align}

One may complete the rest of the argument along \cite[Page 678]{mochizuki-iut4}, but I will argue slightly different as my argument involves minimal amount of algebraic manipulation which is required in the argument detailed in \cite[Page 678]{mochizuki-iut4}.

The important point is that for $\varepsilon>0$ as given in the statement of the theorem, the condition
$$\frac{1+\delta\cdot Q^{-1/2}}{1-\alpha_\lambda}>1+\varepsilon$$
with $\alpha_\lambda$ given by \eqref{eq:alpha-lambda} 
implies that $Q$ is bounded and hence again there are only a finitely many elliptic curves in $U(\bQ)^{\leq d}$ for a given $\varepsilon$. 

Thus enlarging $\exc$ by adding these curves to it without affecting any other parameters, one can assume that
$$\frac{1+\delta\cdot Q^{-1/2}}{1-\alpha_\lambda}<1+\varepsilon$$

Similarly the condition $(1-\alpha_\lambda)^{-1}>2$ leads to $\alpha_\lambda>\frac{1}{2}$ and again this leads to a bound on $Q$ and so omitting elliptic curves for which $(1-\alpha_\lambda)^{-1}>2$ one may assume that $(1-\alpha_\lambda)^{-1}\leq 2$. Thus using both these inequalities one obtains
\begin{align}
\frac{1}{6}\cdot Q&\leq
\left(1+\varepsilon\right)\cdot\left(\log(\dtpd)+\log(\ftpd)\right)+A_\sZ.
\end{align}
By \Cref{pr:height-comparisons} one has $h_{\omega_{\P^1}(D)}\approx \frac{1}{6}Q$ so one obtains
\be
h_{\omega_{\P^1}(D)}\leq
\left(1+\varepsilon\right)\cdot\left(\log(\dtpd)+\log(\ftpd)\right)+A_\sZ.
\ee
By definitions $\logdiff(x_\lambda)=\log(\dtpd)$ and $\logcon_D(x_\lambda)=\log(\ftpd)$. Hence
\be
h_{\omega_{\P^1}(D)}\leq
\left(1+\varepsilon\right)\cdot\left(\logdiff_{\P^1}+\logcon_D\right)+A_\sZ.
\ee
holds on $\sZ\cap U(\bQ)^{\leq d}$ as claimed. This proves \Cref{th:main-dioph-thm}.
\ep

\brem 
In the notation of \Cref{th:existence}, let $\sZ$ be the given compactly bounded subset. Then in \cite[Corollary 2.2]{mochizuki-iut4}, Mochizuki shows that on $\exc$, the function $\fT$ satisfies
$$\fT \lesssim H_{unif}\cdot \varepsilon^{-3}\cdot d^{4+\varepsilon}+H_{\sZ}$$
for some positive constant $H_{unif}$ which is independent of $\sZ$, while $H_{\sZ}$ is a positive constant dependent on $\sZ$. This conclusion can also be arrived at from the proof given above.
\erem

\subsection{Proof of \Cref{th:main-thm}.} 
From \Cref{th:main-dioph-thm} one obtains \Cref{th:main-thm}. In particular:
\bthm 
The \abc\ (\cref{con:abc}) and the arithmetic Szpiro conjecture over $\Q$ are true.
\ethm
\ENDDOCUMENT